\documentclass[12pt, reqno]{amsart}
\usepackage{etex}
\usepackage{amssymb, mathrsfs}
\usepackage{latexsym}
\usepackage{amsmath, accents}
\usepackage{amsthm}
\usepackage{amsfonts}
\usepackage{dsfont, upgreek}
\usepackage{mathtools}
\usepackage{epsfig}
\usepackage{amscd}
\usepackage{graphicx}
\usepackage{bbm}

\usepackage{enumitem}

\usepackage[left=1.4in,right=1.4in,top=1.2in,bottom=1.2in]{geometry}

\usepackage{tikz}
\usetikzlibrary{decorations.markings}
\usetikzlibrary{arrows.meta}

\DeclareMathSymbol{\widetildesym}{\mathord}{largesymbols}{"65}
\newcommand\lowerwidetildesym{%
	\text{\smash{\raisebox{-1.3ex}{%
				$\widetildesym$}}}}
\newcommand\ntilde[1]{%
	\mathchoice
	{\accentset{\displaystyle\lowerwidetildesym}{#1}}
	{\accentset{\textstyle\lowerwidetildesym}{#1}}
	{\accentset{\scriptstyle\lowerwidetildesym}{#1}}
	{\accentset{\scriptscriptstyle\lowerwidetildesym}{#1}}
}

\newcommand{\dbar}{d\hspace*{-0.08em}\bar{}\hspace*{0.1em}}

\usepackage[colorlinks=true,linkcolor=blue,citecolor=blue]{hyperref}

\usepackage[all,cmtip]{xy}
\theoremstyle{plain}
\newtheorem{theorem}{Theorem}[section]
\newtheorem{proposition}[theorem]{Proposition}
\newtheorem{lemma}[theorem]{Lemma}
\newtheorem{corollary}[theorem]{Corollary}

\theoremstyle{definition} 
\newtheorem{definition}[theorem]{Definition}

\theoremstyle{remark} 
\newtheorem{remark}[theorem]{Remark}

\numberwithin{equation}{section}

\newcommand{\ind}{\textup{Ind}}

\newcommand{\tr}{\mathrm{tr}}
\newcommand{\ch}{\mathrm{ch}}
\newcommand{\Ch}{\mathrm{Ch}}
\newcommand{\sw}{\mathrm{T}}
\newcommand{\m}{\mathrm{M}}
\newcommand{\prop}{\mathit{prop\,}}
\newcommand{\dist}{\mathrm{dist}}

\newcommand{\STR}{\mathrm{STR}}

\begin{document}
	
\title[Delocalized eta, cyclic cohomology and higher rho]{Delocalized eta invariants, cyclic cohomology and higher rho invariants
}

\author{Xiaoman Chen}
\address[Xiaoman Chen]{ School of Mathematical Sciences, Fudan University}
\email{xchen@fudan.edu.cn}
\thanks{The first author is partially supported by NSFC 11420101001.}

\author{Jinmin Wang}
\address[Jinmin Wang]{Shanghai Center for Mathematical Sciences}
\email{wangjm15@fudan.edu.cn}
\thanks{The second author is partially supported by NSFC 11420101001.}

\author{Zhizhang Xie}
\address[Zhizhang Xie]{ Department of Mathematics, Texas A\&M University }
\email{xie@math.tamu.edu}
\thanks{The third author is partially supported by NSF 1500823, NSF 1800737}

\author{Guoliang Yu}
\address[Guoliang Yu]{ Department of
	Mathematics, Texas A\&M University}
\email{guoliangyu@math.tamu.edu}
\thanks{The fourth author is partially supported by NSF 1700021, NSF 1564398}

\date{}	
	
\maketitle
\begin{abstract}
    The first main result of this paper is to prove that  the convergence of Lott's delocalized eta invariant holds for all differential operators with a sufficiently large spectral gap at zero. Furthermore, to each delocalized cyclic cocycle, we define a higher analogue of Lott's delocalized eta invariant and prove its convergence when the delocalized cyclic cocycle has at most exponential growth. As an application, for each cyclic cocycle of at most exponential growth, we prove a formal higher Atiyah-Patodi-Singer index theorem on manifolds with boundary,  under the condition that the operator on the boundary has a sufficiently large spectral gap at zero.

    Our second main result is to  obtain an explicit formula of the delocalized Connes-Chern character of all $C^\ast$-algebraic secondary invariants for word hyperbolic groups. Equivalently, we give an explicit formula for the pairing between $C^\ast$-algebraic secondary invariants and delocalized cyclic cocycles of the group algebra.   When the $C^\ast$-algebraic secondary invariant is a $K$-theoretic higher rho invariant of an invertible differential operator, we show this pairing is precisely the higher analogue of Lott's delocalized eta invariant alluded to above. Our work uses Puschnigg's smooth dense subalgebra for word hyperbolic groups in an essential way. We emphasize that our construction of the delocalized Connes-Chern character is  at $C^\ast$-algebra $K$-theory level.  This is of essential importance for applications to geometry and topology. As a consequence,   we compute the paring between delocalized cyclic cocycles and  $C^\ast$-algebraic Atiyah-Patodi-Singer  index classes for manifolds with boundary, when the fundamental group of the given manifold is hyperbolic. In particular, this improves the formal delocalized higher Atiyah-Patodi-Singer theorem from above and removes the condition that the spectral gap of the operator on the boundary is sufficiently large.

\end{abstract}
\section{Introduction}

Higher index theory is a far-reaching generalization of the classic Fredholm index theory by taking into consideration of the symmetries of the underlying space.  Let $X$ be a complete Riemannian manifold of dimension $n$ with a discrete group $G$ acting on it properly and cocompactly by isometries. Each $G$-equivariant elliptic differential operator $D$ on $X$ gives rise to a higher index class $\ind_G(D)$ in  the $K$-group $K_n(C_r^*(G))$ of the reduced group $C^\ast$-algebra $C_r^\ast(G)$. This higher index is an obstruction to the invertibility of $D$. The higher index theory plays a fundamental role in the studies of many problems in geometry and topology such as the Novikov conjecture, the Baum-Connes conjecture and the Gromov-Lawson-Rosenberg conjecture. Higher index classes are invariant under homotopy and often referred to as primary invariants.

 When the higher index class of an operator is trivial and given a specific trivialization, a secondary index theoretic invariant naturally arises. One such example is the associated Dirac operator on the universal covering $\widetilde M$ of a closed spin manifold $M$ equipped with a positive scalar curvature metric. It follows from the Lichnerowicz formula that the Dirac operator on $\widetilde M$ is invertible. In this case, there is a natural $C^\ast$-algebraic secondary  invariant introduced by Higson and Roe in \cite{Higson1, Higson2, Higson3, Roe}, called the higher rho invariant, which lies in $K_n(C_{L,0}^\ast(\widetilde M)^G)$, where $G$ is  the fundamental group $ \pi_1(M)$ of $M$ and $C_{L,0}^\ast(\widetilde M)^G$ is a certain geometric $C^\ast$-algebra.   The precise definition of  $C_{L,0}^\ast(\widetilde M)^G$ and that of the higher rho invariant are given in Section \ref{preliminaries}. This higher rho invariant is an obstruction to the inverse of the Dirac operator being local and has important applications to geometry and topology. 

Parallel to the $C^\ast$-algebraic approach above, Lott developed a theory of secondary invariants in the framework of noncommutative differential forms  \cite{Lott}. Lott's theory was in turn very much inspired by the work of Bismut and Cheeger on eta forms \cite{JBJC89}, which naturally arise in the index theory for families of manifolds with boundary. 
Despite the fact that Lott's higher eta invariant is defined by an explicit integral formula of noncommutative differential forms, it is difficult to compute in general.  It is only after one pairs Lott's higher eta invariant with cyclic cocycles of $\pi_1(M)$ that it becomes more computable and more applicable to problems in geometry and topology. However, due to certain convergence issues, the question when such a pairing can actually be rigorously defined is often very subtle. We shall devote  the first half of the current paper to these convergence issues. We show that the pairing of Lott's higher eta invariant and a (delocalized) cyclic cocycle is well-defined, under the condition that the operator on $\widetilde M$ has a sufficiently large spectral gap at zero and  the cyclic cocycle has at most exponential growth. In particular, as a special case, if both $\pi_1(M)$ and the cyclic cocycle have sub-exponential growth, then the pairing is always well-defined for all invertible operators on $\widetilde M$.

The second goal of our paper  is to obtain an explicit formula for the delocalized Connes-Chern character of all $C^\ast$-algebraic secondary invariants for hyperbolic groups. Equivalently, this amounts to computing the pairing between $C^\ast$-algebraic secondary invariants and delocalized cyclic cocycles of the group algebra. In the case where the $C^\ast$-algebraic secondary invariant is a $K$-theoretic higher rho invariant, then the pairing is given explicitly in terms of Lott's higher eta invariant\footnote{In the literature, the delocalized part of Lott's noncommutative-differential higher eta invariant sometimes is also referred to as higher rho invariant. To avoid confusion, we shall refer this  noncommutative-differential higher rho invariant as \emph{Lott's higher eta invariant}.}, or rather its periodic version\footnote{See \cite[Section 4.6]{Lott1} for a discussion of the periodic version of Lott's higher eta invariant. As $K$-theory is $2$-periodic, the periodic higher eta invariant is the correct version to used here.}. As mentioned above, one of the main technical difficulties is to resolve various convergence issues. In the case of hyperbolic groups, we overcome these convergence issues with the help of Puschnigg's smooth dense subalgebra. As a consequence,  we compute the paring between delocalized cyclic cocycles and  Atiyah-Patodi-Singer higher index classes for manifolds with boundary in terms of  delocalized higher eta invariants, when the fundamental group of the given manifold is hyperbolic. The details of these results will occupy the second half of the paper.

In the following, we shall give a more precise overview of some of the main results of this paper. Let us recall the definition of Lott's delocalized eta invariant, which shall be thought of (at least formally) as a pairing between Lott's higher eta invariant and traces (i.e. degree zero cyclic cocycles).  Suppose $\left\langle h\right\rangle$ is a nontrivial conjugacy class in $\pi_1(M)$ in the sense that the group element $h$ is not equal to the identity in $\pi_1(M)$. If $D$ is a self-adjoint elliptic differential operator on $M$ and $\widetilde D$ is the lifting of $D$ to $\widetilde M$, then the delocalized eta invariant $\eta_{\left\langle h\right\rangle}(\widetilde D)$ is defined by the formula:
\begin{equation}\label{eq:delocaleta}
\eta_{\left\langle h\right\rangle }(\widetilde D):=
\frac{2}{\sqrt\pi}\int_{0}^\infty \tr_{\left\langle h\right\rangle }(\widetilde De^{-t^2\widetilde D^2})dt.
\end{equation}
Here $\tr_{\left\langle h\right\rangle }$ is the following trace map
$$\tr_{\left\langle h\right\rangle }(A)=\sum_{g\in\left\langle h\right\rangle }\int_{x\in\mathcal F}A(x,gx)dx$$
on $G$-equivariant kernels $A\in C^{\infty}(\widetilde M\times\widetilde M)$, where $\mathcal F$ is a fundamental domain of $\widetilde M$ under the action of $G = \pi_1(M)$. As it stands, the above definition of delocalized eta invariant does not require a choice  of a smooth dense subalgebra of $C^\ast_r(G)$. Of course, in the special event that  $\widetilde De^{-t^2\widetilde D^2}$ lies in an appropriate smooth dense subalgebra   to which the trace map $\tr_{\langle h\rangle}$ continuously extends, this delocalized eta invariant indeed coincides with the  
pairing of $\tr_{\langle h\rangle}$ with Lott's  higher eta invariant.

In \cite{Lott},  the convergence of the above formula is proved by Lott under the assumption that $\left\langle h\right\rangle$ has polynomial growth or is hyperbolic\footnote{In  \cite[Proposition 8]{Lott}, Lott stated that the convergence of the formula $\eqref{delocalizedeta}$ holds for both groups with polynomial growth and  hyperbolic groups. However, his proof for hyperbolic groups contained a technical problem which was later fixed by Puschnigg after constructing a different smooth dense subalgebra of the reduced group $C^\ast$-algebra for hyperbolic groups \cite{Puschnigg}.},  and that $\widetilde D$ is invertible or more generally $\widetilde D$ has a spectral gap at zero. Recall that  $\widetilde D$ is said to have a spectral gap at zero if there exists an open interval $(-\varepsilon, \varepsilon)\subset \mathbb R$ such that $\textup{spectrum}(\widetilde D) \cap (-\varepsilon, \varepsilon)$ is either $\{0\}$ or empty.    In general, the convergence of $\eqref{eq:delocaleta}$ fails. For example, Piazza and Schick gave an explicit example where the convergence of $\eqref{eq:delocaleta}$ fails when $\widetilde D$ does not have a spectral gap at zero.

As the first main result of this paper, we show that if $\widetilde D$ has a sufficiently large spectral gap at zero, then  Lott's delocalized eta invariant $\eta_{\langle h \rangle}(\widetilde D)$ in line $\eqref{eq:delocaleta}$ converges absolutely. We refer to Definition $\ref{def:speclarge}$ in Section $\ref{HigherEtaInvariants}$ for a more precise quantitative explanation of what it means for a spectral gap to be ``sufficiently large". 

\begin{theorem}\label{thm:converge}
    Let $M$ be a closed manifold and $\widetilde M$ the universal covering over $M$. Suppose $D$ is a self-adjoint first-order elliptic differential operator over $M$ and $\widetilde D$ the lift of $D$ to $\widetilde M$. If $\left\langle h\right\rangle$ is a nontrivial conjugacy class of $\pi_1(M)$ and $\widetilde D$ has a sufficiently large spectral gap at zero, then the delocalized eta invariant $\eta_{\left\langle h\right\rangle }(\widetilde D)$  defined in line $\eqref{eq:delocaleta}$ converges absolutely. 
\end{theorem}

We would like to emphasize that the theorem above works for all fundamental groups.  In the special case where the conjugacy class $\langle h\rangle$ has sub-exponential growth, then any nonzero spectral gap is in fact sufficiently large, hence in this case $\eta_{\left\langle h\right\rangle }(\widetilde D)$ converges absolutely as long as $\widetilde D$ is invertible.  

Now a special feature of traces is that they always have uniformly bounded representatives, when viewed as degree zero cyclic cocycles. In fact, the techniques used to prove Theorem $\ref{thm:converge}$ above can be generalized to all delocalized cyclic cocycles of higher degrees, as long as they have at most exponential growth.  Recall that the cyclic cohomology of a group algebra  decomposes into a direct product with respect to the conjugacy classes of the group. A cyclic cocycle in a component of this direct product decomposition that corresponds to a nontrivial conjugacy class  $\langle h\rangle$  will be called a delocalized cyclic cocycle at  $\langle h\rangle$, cf. Definition $\ref{def:delocalcocyc}$. Moreover, see Definition $\ref{def:exgrow}$ in Section $\ref{HigherEtaInvariants}$ for the precise definition of exponential growth for cyclic cocycles.

\begin{theorem}\label{thm:higherconverge}
	Assume the same notation as in Theorem $\ref{thm:converge}$. Let $\varphi$ be a delocalized cyclic cocycle at a nontrivial conjugacy class $\langle h\rangle$. If $\varphi$ has exponential growth and $\widetilde D$ has a sufficiently large spectral gap at zero, then a higher analogue   $\eta_{\varphi}(\widetilde D)$  \textup{(}cf.  Definition $\ref{higherdelocalizedeta}$\textup{)} of the delocalized eta invariant converges absolutely.  
\end{theorem}

    The explicit formula for  $\eta_{\varphi}(\widetilde D)$ is described in terms of  the transgression formula for Connes-Chern character \cite{Connes,ConnesNCG} \cite{JLO88}. It is closely related to the periodic version of Lott's noncommutative-differential higher eta invariant.   In the case where the fundamental group $G$ has polynomial growth, we shall show that our formula for $\eta_{\varphi}(\widetilde D)$ is equivalent to the periodic version of Lott's noncommutative-differential higher eta invariant, cf. Section $\ref{Identification}$. As $\eta_{\varphi}(\widetilde D)$ is an analogue for higher degree cyclic cocycles  of  Lott's delocalized  eta invariant, we shall call $\eta_{\varphi}(\widetilde D)$ a \emph{delocalized higher eta invariant} from now on. 
    Again,  we refer to Definition $\ref{def:specbd}$  for a more precise quantitative explanation of what it means for a spectral gap to be ``sufficiently large" in this context. For now, let us just point out that if both $G$ and $\varphi$ have sub-exponential growth,  then any nonzero spectral gap is in fact sufficiently large, hence in this case $\eta_{\varphi}(\widetilde D)$ converges absolutely as long as $\widetilde D$ is invertible.

    Formally speaking, just as Lott's delocalized eta invariant $\eta_{\left\langle h\right\rangle }(\widetilde D)$ can be interpreted as the pairing between the degree zero cyclic cocycle $\tr_{\langle h\rangle}$ and the \mbox{$K$-theoretic} higher rho invariant $\rho(\widetilde D)$ (or the noncommutative differential higher eta invariant), so can the delocalized higher eta invariant $\eta_{\varphi}(\widetilde D)$ be interpreted as the pairing between the cyclic cocycle $\varphi$ and the \mbox{$K$-theoretic} higher rho invariant $\rho(\widetilde D)$ (or the noncommutative differential higher eta invariant). As pointed out in the discussion above, a key analytic difficulty here is to verify when such a pairing is well-defined, or more ambitiously, to verify when one can extend this pairing to a pairing between the cyclic cohomology  of $\mathbb CG$ and the $K$-theory group $K_\ast(C_{L,0}^\ast(\widetilde M)^G)$. The group $K_\ast(C_{L,0}^\ast(\widetilde M)^G)$ consists of $C^\ast$-algebraic secondary  invariants; in particular, it contains all higher rho invariants from the discussion above. As pointed out above, such an extension of the pairing is important, often necessary, for many interesting applications to geometry and topology (cf. \cite{Piazzarho, MR3590536, Weinberger:2016dq}).

 In a previous paper \cite{Xie}, the third and fourth authors established a pairing between delocalized cyclic cocycles of degree zero (i.e. delocalized traces) and the $K$-theory group $K_\ast(C_{L,0}^\ast(\widetilde M)^G)$,  under the assumption that the relevant conjugacy class has polynomial growth. In this paper,  we shall construct a pairing  between delocalized cyclic cocycles of \emph{all degrees} and the $K$-theory group $K_\ast(C_{L,0}^\ast(\widetilde M)^G)$ for hyperbolic groups. Before we state the theorem, let us recall some notation that will used in the statement of the next theorem.  The cyclic cohomology of a group algebra $\mathbb CG$ has a decomposition respect to the conjugacy classes of $G$ (\cite{MR814144, Nistor}):
 \[ HC^*(\mathbb CG)\cong\prod_{\left\langle h\right\rangle }HC^*(\mathbb CG,\left\langle h\right\rangle), \]
 where 
 $HC^*(\mathbb CG,\left\langle h\right\rangle)$  denotes the component that corresponds to the conjugacy class $\langle h\rangle$, cf. Definition $\ref{def:delocalcocyc}$ below.
 
 \begin{theorem}\label{main}
 	Let $M$ be a closed manifold whose fundamental group $G$ is  hyperbolic. Suppose  $\left\langle h\right\rangle $  is a non-trivial conjugacy class of $G$. Then every element $[\alpha]\in HC^{2k+1-i}(\mathbb CG,\left\langle h\right\rangle ) $ induces a natural map 
 	$$\tau_{[\alpha]} \colon K_i(C^*_{L,0}(\widetilde M)^G)\to \mathbb C$$ such that the following are satisfied. 
 	\begin{enumerate}[label=$(\roman*)$]
 		\item $\tau_{[S\alpha]} = \tau_{[\alpha]}$, where $S$ is Connes' periodicity map \[ S\colon HC^{\ast}(\mathbb CG,\left\langle h\right\rangle )\to HC^{\ast+2}(\mathbb CG,\left\langle h\right\rangle ).\]
 		\item Suppose $D$ is a first-order elliptic differential operator on $M$ such that the lift $\widetilde D$ of $D$ to the universal cover $\widetilde M$ of $M$ is invertible. 
 		Then we have 
 		$$\tau_{[\alpha]}(\rho(\widetilde D))= - \eta_{\alpha}(\widetilde D),  $$
 		where $\rho(\widetilde D)$ is the $C^\ast$-algebraic higher rho invariant of $\widetilde D$ and  $\eta_{\alpha}(\widetilde D)$ is the  delocalized higher eta invariant defined in Definition $\ref{higherdelocalizedeta}$. In particular, in this case, the delocalized higher eta invariant $\eta_{\alpha}(\widetilde D)$   converges absolutely.  
 	\end{enumerate}
 \end{theorem}

The construction of the map $\tau_{[\alpha]}$ in the above theorem uses Puschnigg's smooth dense subalgebra  for hyperbolic groups  \cite{Puschnigg} in an essential way. In more conceptual terms, the above theorem provides an explicit formula to compute the  delocalized Connes-Chern character of $C^\ast$-algebraic secondary invariants. More precisely, the same techniques developed in this paper actually imply\footnote{In fact, even more is true. The same techniques developed in the current paper  imply that if $\mathcal A$ is   smooth dense subalgebra of $C_r^\ast(\Gamma)$ for any group $\Gamma$ (not necessarily hyperbolic) and in addition $\mathcal A$ is a Fr\'echet locally $m$-convex algebra, then there is a well-defined delocalized Connes-Chern character $Ch_{deloc}\colon K_i(C^*_{L,0}(\widetilde M)^\Gamma)\to \overline{HC}^{deloc}_{\ast}(\mathcal A)$. Of course,  in order to pair such a delocalized Connes-Chern character with a cyclic cocycle of $\mathbb C\Gamma$, the key remaining challenge is to  continuously extend this cyclic cocycle of $\mathbb C\Gamma$ to a cyclic cocycle of $\mathcal A$.  } that there is a well-defined delocalized Connes-Chern character 
\[ \ch_{deloc}\colon K_i(C^*_{L,0}(\widetilde M)^\Gamma)\to \overline{HC}^{deloc}_{\ast}(B(\mathbb CG)), \] where $B(\mathbb CG)$ is Puschnigg's smooth dense subalgebra of the reduced group \mbox{$C^\ast$-algebra} of $G$ and $\overline{HC}^{deloc}_{\ast}(B(\mathbb CG))$ is the delocalized part of the cyclic homology\footnote{Here the definition of cyclic homology of $B(\mathbb CG)$ takes the topology of $B(\mathbb CG)$ into account, cf. \cite[Section II.5]{Connes}. } of $B(\mathbb CG)$. Now for Gromov's hyperbolic groups, every cyclic cohomology class of $\mathbb C\Gamma$ continuously extends to  cyclic cohomology class of $B(\mathbb CG)$ (cf. \cite{Puschnigg} for the case of degree zero cyclic cocycles and Section $\ref{SmoothDenseSubalgebras}$ of this paper for the case of higher degree cyclic cocycles). Thus  the map $\tau_{[\alpha]}$ can be viewed as a pairing between  cyclic cohomology and delocalized Connes-Chern characters of $C^\ast$-algebraic secondary invariants.   We point out that, although the spectral gap of $\widetilde D$ is required to be sufficiently large in Theorem $\ref{thm:converge}$ and $\ref{thm:higherconverge}$ in order for $\eta_{\alpha}(\widetilde D)$ to converge, such a requirement is \emph{not} needed in the case of hyperbolic groups. This is again a consequence of some essential properties of Puschnigg's smooth dense subalgebra.

As an application, we use this delocalized Connes-Chern character map to  obtain a delocalized higher Atiyah-Patodi-Singer index theorem for manifolds with boundary. More precisely, let  $W$ be a compact $n$-dimensional spin manifold with boundary $\partial W$. Suppose $W$ is equipped with a Riemannian metric $g_W$ which has product structure near $\partial W$ and in addition has positive scalar curvature on $\partial W$. Let $\widetilde W$ be the universal covering of $W$ and $g_{\widetilde W}$  the   Riemannian metric  on $\widetilde W$ lifted from $g_W$.  With respect to the metric $g_{\widetilde W}$,  the associated Dirac operator $\widetilde D_W$ on $\widetilde  W$ naturally defines a higher index $\ind_G(\widetilde D_W)$ in $K_n(C^\ast(\widetilde W)^G) = K_n(C_r^\ast(G))$, where $G = \pi_1(W)$, cf. \cite[Section 3]{Xiepos}. Since the metric $g_{\widetilde W}$ has positive scalar curvature on $\partial \widetilde W$, it follows from the Lichnerowicz formula that the associated Dirac operator $\widetilde D_\partial$ on $\partial {\widetilde W}$ is invertible, hence naturally defines a higher rho invariant $\rho(\widetilde D_\partial)$ in $K_{n-1}(C_{L, 0}^\ast(\widetilde W)^\Gamma)$. We have the following delocalized higher Atiyah-Patodi-Singer index theorem. 

\begin{theorem}\label{thm:delocalaps}
	 With the same notations as above, if $G=\pi_1(W)$ is hyperbolic and $\left\langle h\right\rangle $ is a nontrivial conjugacy class of $G$, then for any $[\varphi]\in HC^{\ast}(\mathbb CG,\left\langle h\right\rangle )$,
	we have \begin{equation}
	\ch_{[\varphi]}(\ind_G(\widetilde D_W))= \frac{1}{2}\eta_{[\varphi]}(\widetilde D_\partial ), 
	\end{equation}
	where $\ch_{[\varphi]}(\ind_G(\widetilde D_W))$ is the Connes-Chern pairing between the cyclic cohomology class $[\varphi]$ and the $C^\ast$-algebraic index class $\ind_G(\widetilde D_W)$.
\end{theorem}

There have been various versions of higher Atiyah-Patodi-Singer theorem in the literature \cite{ MR1708639, Charlotte, MR3500822}. See the discussion after Theorem $\ref{application}$ for more details on the relations and differences of the above theorem with those existing results.

The methods developed in this paper can also be applied to prove analogues of Theorem $\ref{main}$ and Theorem $\ref{thm:delocalaps}$ above for virtually nilpotent groups. The complete details can be found in the thesis of Sheagan John \cite{Sheagan}.

We would like to point out that our proof of Theorem $\ref{main}$ does \emph{not} rely on the Baum-Connes isomorphism for hyperbolic groups \cite{MR2874956, MR1914618}, although the theorem is closely connected to the Baum-Connes conjecture and the Novikov conjecture. On the other hand, if one is willing to use the full power of the Baum-Connes isomorphism for hyperbolic groups, there is in fact a different, but more indirect, approach to the delocalized Connes-Chern character map. First, observe that the map $\tau_{[\alpha]}$ factors through a map
\[ \tau_{[\alpha]}\colon  (K_i(C^*_{L,0}(\underline{E}G)^G)\otimes \mathbb C)\to \mathbb C\]
where $\underline{E}G$ is the universal space for proper $G$-actions.  Now the Baum-Connes isomorphism $\mu\colon K_\ast^G(\underline{E}G) \xrightarrow{ \ \cong \ }K_\ast(C^\ast_r(G))$ for hyperbolic groups implies that one can identify $K_i(C^*_{L,0}(\underline{E}G)^G)\otimes \mathbb C$  with  $\bigoplus_{\left\langle h\right\rangle \neq 1}HC_\ast(\mathbb CG,\left\langle h\right\rangle)$, where $HC_\ast(\mathbb CG,\left\langle h\right\rangle)$ is the delocalized cyclic homology at $\langle h\rangle$ (a cyclic homology analogue of Definition $\ref{def:delocalcocyc}$) and the direct sum is taken over all nontrivial conjugacy classes. In particular, after this identification, it follows that the map $\tau_{[\alpha]}$ becomes the usual componentwise pairing between cyclic cohomology and cyclic homology. However, for a specific element, e.g. the higher rho invariant $\rho(\widetilde D)$,  in  $K_i(C^*_{L,0}(\underline{E}G)^G)$, its identification with an element in  $\bigoplus_{\left\langle h\right\rangle \neq 1}HC_\ast(\mathbb CG,\left\langle h\right\rangle)$ is rather abstract and implicit. More precisely,  the   computation of the number  $\tau_{[\alpha]}(\rho(\widetilde D))$ essentially amounts to the following process. Observe that if a closed spin manifold $M$ is equipped with a positive scalar curvature metric, then stably it bounds (more precisely,  the universal cover $\widetilde M$ of $M$   becomes the boundary of another $G$-manifold,  after finitely many steps of cobordisms and  vector bundle modifications). In principle,  the number $\tau_{[\alpha]}(\rho(\widetilde D))$ can be derived from a higher Atiyah-Patodi-Singer index theorem for this bounding manifold. The drawback of such an approach is that there is no explicit  formula for $\tau_{[\alpha]}(\rho(\widetilde D))$, since there is no explicit procedure for producing such a bounding manifold.   
In \cite{MR3563244}, Deeley and Goffeng also constructed a delocalized Connes-Chern character for $C^\ast$-algebra secondary invariants. Their approach is in spirit  similar to the indirect method just described above (making use of the Baum-Connes isomorphism for hyperbolic groups),   although their actual technical implementation is different. The key feature of our approach in this paper is that we obtain  an explicit and intrinsic formula for the delocalized Connes-Chern character of $C^\ast$-algebraic secondary invariants.

The paper is organized as follows. In Section \ref{preliminaries}, we review some standard geometric  $C^*$-algebras and a construction of Higson-Roe's $K$-theoretic higher rho invariants. In Section \ref{HigherEtaInvariants}, we prove the convergence of Lott's delocalized eta invariant holds for all operators with a sufficiently large spectral gap at zero. More generally, for each higher degree delocalized cyclic cocycle, we define a higher analogue of Lott's delocalized eta invariant, and prove its convergence for all operators with a sufficiently large spectral gap at zero, provided that the given cyclic cocycle has at  most exponential growth. In Section \ref{SmoothDenseSubalgebras}, we review Puschnigg's construction of smooth dense subalgebras of reduced group $C^*$-algebras for hyperbolic groups. Puschnigg showed that any trace on the group algebra of a hyperbolic group extends continuously onto this smooth dense subalgebra. We shall generalize this result to cyclic cocycles of all degrees, provided that the cyclic cocycles have polynomial growth (cf. Proposition \ref{highertrace} below).  In Section \ref{CyclicCohomologyGrowth}, we show that every cyclic cohomology class of a hyperbolic group  has a representative of polynomial growth. Furthermore, if the cyclic cohomology class has degree $\geqslant 2$, then it admits a uniformly bounded representative.  In Section \ref{DelocalizedPairing},  we give an explicit formula for the pairing between $C^\ast$-algebraic secondary invariants and delocalized cyclic cocycles of the group algebra for word hyperbolic groups.   When the $C^\ast$-algebraic secondary invariant is a $K$-theoretic higher rho invariant of an invertible differential operator, we show this pairing is precisely the higher  delocalized eta invariant of the given operator.   In Section \ref{Application}, we compute the paring between delocalized cyclic cocycles and  $C^\ast$-algebraic Atiyah-Patodi-Singer  index classes for manifolds with boundary, when the fundamental group of the given manifold is hyperbolic. In section \ref{Identification}, we identify our definition of delocalized higher eta invariant with Lott's higher eta invariant. 

We would like to thank Denis Osin for providing us a proof of a useful result on hyperbolic groups (Lemma $\ref{lm:conj}$). 

\section{Preliminaries}\label{preliminaries}

In this section, we review the construction of some geometric $C^*$-algebras and Higson-Roe's higher rho invariants. We refer the reader to \cite{Roecoarse,Yulocalization, Higson1, Higson2, Higson3} for more details.

Let $X$ be a proper metric space, that is, every closed metric ball in $X$ is compact.
An $X$-module is a separable Hilbert space equipped with a $*$-representation of $C_0(X)$, the algebra of all continuous functions on $X$ which vanish at infinity.
An $X$-module is called nondegenerate if the $*$-representation of $C_0(X)$ is nondegenerate.
An $X$-module is said to be standard if no nonzero function in $C_0(X)$ acts as a compact operator.

\begin{definition}\label{def propa local pseudo}
Let $H_X$ be an $X$-module and $T$ a bounded linear operator acting on $H_X$.
\begin{enumerate}
\item The propagation of $T$ is defined to be $\sup\{d(x,y)| (x,y)\in supp(T) \},$ where $supp(T)$ is the complement (in $X\times X$) of the set of points $(x,y) \in X\times X$ for which there exist $f,g\in C_0(X)$ such that $gTf=0$ and $f(x)g(y)\neq 0$.
\item $T$ is said to be locally compact if $fT$ and $Tf$ are compact for all $f\in C_0(X)$.
\item $T$ is said to be pseudo-local if $[T,f]$ is compact for all $f\in C_0(X)$.
\end{enumerate}
\end{definition}

\begin{definition}\label{def roe and localization}
Let $H_X$ be a standard nondegenerate $X$-module and $B(H_X)$ the set of all bounded linear operators on $H_X$.
\begin{enumerate}
\item The Roe algebra of $X$, denoted by $C^*(X)$, is the $C^*$-algebra generated by all locally compact operators with finite propagations in $B(H_X)$.
\item $D^*(X)$ is the $C^*$-algebra generated by all pseudo-local operators and with finite propagations in $B(H_X)$. In particular, $D^*(X)$ is a subalgebra of the multiplier algebra of $C^*(X)$.
\item $C_L^*(X)$ (resp. $D^*_L (X)$) is the $C^*$-algebra generated by all bounded and uniformly-norm continuous functions $f : [0,\infty)\to C^*(X)$ (resp. $f : [0,\infty)\to D^*(X)$) such that 	
\[
\text{propagation of }f(t)\to 0\  \text{~as~}  t \to \infty.
\]
$D^*_L (X)$ is a subalgebra of the multiplier algebra of $C_L^*(X)$.
\item The kernel of the following evaluation map
\begin{eqnarray*}
ev:\;C_L^*(X)\to  C^*(X),\ f\mapsto   f(0)
\end{eqnarray*}
is defined to be $C^*_{L,0}(X)$. In particular, $C^*_{L,0}(X)$ is an ideal of $C_L^*(X)$. Similarly, we define $D^*_{L,0}(X)$ as the kernel of the evaluation map from $D_L^*(X)$ to $D^*(X)$.
\end{enumerate}
\end{definition}

Now in addition we assume that there is a countable discrete group $G$, and acts properly on $X$ by isometries. In particular, if the action of $\Gamma$ is free, then $X$ is simply a $G$-covering of the compact space $X/G$.
Let $H_X$ be an $X$-module equipped with a covariant unitary representation of $G$.
Let the representation of $C_0(X)$ be $\phi$ and let the action of $G$ be $\pi$. We call $(H_X, G, \phi,\pi)$ a covariant system is to say
\[
\pi(g) \phi(f)=\phi(g^*f)\pi(g),
\]
where $g^*f(x)= f(g^{-1}x)$ for any $f\in C_0(X)$, $g \in G$.

\begin{definition}
A covariant system $(H_X, G, \phi)$ is called \emph{admissible} if
\begin{enumerate}
\item The action of $G$ is proper and cocompact;
\item $H_X$ is a nondegenerate standard $X$-module;
\item For each $x\in X$, the stabilizer group $G_x$ acts on $H_X$ regularly in the sense that the action is isomorphic to the obvious action of $G_x$ on $l^2(G_x)\otimes H$ for some infinite dimensional Hilbert space $H$.
Here $G_x$ acts on $\ell^2 (G_x)$ by (left) translations and acts on $H$ trivially.
\end{enumerate}
\end{definition}

We remark that for each locally compact metric space $X$ with a proper and cocompact isometric action of $G$, there exists an admissible covariant system $(H_X, G,\phi)$. Also, we point out that the condition 3 above is automatically satisfied if $G$ acts freely on $X$. If no confusion arises, we will denote an admissible covariant system $(H_X, G,\phi)$ by $H_X$ and call it an admissible $(X, G)$-module.

\begin{definition}
Let $X$ be a locally compact metric space $X$ with a proper and cocompact isometric action of $G$.
If $H_X$ is an admissible $(X, G)$-module, we denote by $\mathbb{C}[X]^G$ to be $*$-algebra of all $G$-invariant locally compact operators with finite propagations in $B(H_X)$.
We define $C^*(X)^G$ to be the completion of $\mathbb{C}[X]^G$ in $B(H_X)$.
\end{definition}

Since the action of $G$ on $X$ is cocompact, it is known that $C^*(X)^G$ is $*$-isomorphic to $\mathcal K\otimes C^*_r(G)$, where $\mathcal K$ is the algebra of all compact operators and $C^*_r(G)$ is the reduced group $C^*$-algebra.

Similarly, we can also define $D^*(X)^G$, $C^*_L (X)^G$, $D^*_L (X)^G$, $C^*_{L,0}(X)^G$ and $D^*_{L,0}(X)^G$.
\begin{remark}
Up to isomorphism, $C^*(X)$ does not depend on the choice of the standard nondegenerate $X$-module $H_X$.
The same holds for $D^*(X)$, $C_L^*(X)$, $D^*_L (X)$, $C_{L,0}^*(X)$, $D^*_{L,0}(X)$ and their $G$-equivariant versions.
\end{remark}

Let $M$ be a closed Riemannian manifold. Let $G$ be a discrete finitely generated countable group. Suppose $\widetilde M$ is a regular $G$-cover of $M$. For example,  $\widetilde M$ is the universal covering of $M$ and $G$ is the fundamental group of $M$. Let $p$ be the associated covering map from $\widetilde M$ to $M$. Suppose $E$ is an Hermitian vector bundle over $M$ and $\widetilde E$ the lifting of $E$ to $\widetilde M$. Write $\mathcal H$ the collection of all $L^2$-sections of $\widetilde E$. The equivariant Roe algebra $C^*(\widetilde M)^G$ is defined to be the operator-norm completion of all $G$-equivariant locally compact operators of finite propagation acting on $\mathcal H$. The localization algebra $C^*_L (\widetilde M)^G$ and $C^*_{L,0}(\widetilde M)^G$ can be defined similarly.

Suppose that $M$ is spin and $E$ is the corresponding spinor bundle over $M$. Let $D$ be the Dirac operator acting on $E$ and $\widetilde D$ its lifting to $\widetilde E$. If the scalar curvature of the metric over $M$ is strictly positive, the Dirac operator naturally defines a $K$-theory class called the higher rho invariant in $K_*(C_{L,0}^*(\widetilde M)^G)$. For simplicity, we will only discuss the case where $M$ is odd dimensional; the even dimensional case is completely similar, cf. \cite{Xiepos}. 

Define the following functions
\begin{equation}\label{higherrho}
    \begin{split}
        F_t(x)=\frac{1}{\sqrt{\pi}}\int_{-\infty}^{x/t} e^{-s^2}ds \ \textup{ and } \
        U_t(x)=\exp(2\pi iF_t(x)).
    \end{split}
\end{equation}
Since the scalar curvature over $\widetilde M$ is uniformly bounded below by a positive number, it follows from the Lichnerowicz formula that the Dirac operator $\widetilde D$ is invertible. This implies that $F_t(\widetilde D)$ converges to $\frac{1}{2}(1+\widetilde D|\widetilde D|^{-1})$ in operator norm, as $t\to 0$. Thus the path $\{U_t(\widetilde D)\}_{0\leq t <\infty}$ lies in $C_{L,0}^*(\widetilde M)^G$.

\begin{definition}[{\cite{Higson1, Higson2, Higson3}}]\label{def:highrho}	The Higson-Roe higher rho invariant $\rho(\widetilde D)$ of $\widetilde D$  is defined to be the $K$-theory class $[\{U_t(\widetilde D)\}_{0\leq t <\infty}]\in K_1(C_{L,0}^*(\widetilde M)^G)$. 
\end{definition}

\section{Higher eta invariants}\label{HigherEtaInvariants}

In this section, we show  that Lott's delocalized eta invariant converges absolutely for operators with a sufficiently large spectral gap at zero. More generally, we prove the convergence of a higher analogue of Lott's delocalized eta invariant for operators with a sufficiently large spectral gap at zero and  higher degree cyclic cocycles that have at most exponential growth. 

\subsection{Convergence of delocalized eta invariants with a large enough spectral gap}

Let $M$ be a closed Riemannian manifold. Let $G$ be a finitely generated discrete group. Suppose $\widetilde M$ is a regular $G$-covering space of $M$. Let $p$ be the associated covering map from $\widetilde M$ to $M$. Choose any fundamental domain $\mathcal F$ of $G$-action on $\widetilde M$. Suppose $D$ is a first-order self-adjoint elliptic differential operator acting on some Hermitian vector bundle $E$ over $M$ and $\widetilde E$ (resp. $\widetilde D$) is the lifting of $E$ (resp. $D$) to $\widetilde M$. Assume $G$ acts on $\widetilde E$ as well. Let $\mathcal H$ be the space of $L^2$-sections of $\widetilde E$. Then the operator
$\widetilde De^{-t^2\widetilde D}$ lies in $B(\mathcal H)$, the algebra of bounded operators on $\mathcal H$. Moreover,  the associated Schwartz kernel $k_t(x,y)$ of $\widetilde De^{-t^2\widetilde D}$ is smooth.

Any conjugacy class $\langle h \rangle$ of $G$ naturally induces a trace map
$$tr_{\left\langle h\right\rangle }\colon \mathbb CG\to\mathbb C, \textup{ by }  \sum_{g\in G}a_g g\mapsto \sum_{g\in\left\langle h\right\rangle }a_g.$$ In particular, when $\left\langle h\right\rangle$ is the trivial conjugacy class, $tr_{\langle e\rangle}$ is the canonical trace on $\mathbb CG$. 
The trace map $tr_{\left\langle h\right\rangle }$ generalizes (formally) to a trace map on  $G$-equivariant integral operators $T$ with a smooth Schwartz kernel $T(x,y)$ as follows: 
\begin{equation}
    \tr_{\langle h\rangle}(T)=\sum_{g\in\langle h\rangle}\int_{x\in\mathcal F} tr(T(x,gx))dx,
\end{equation}
provided that the right hand side converges. 

\begin{definition}[{\cite{Lott}}]\label{def:deloceta}
	For any nontrivial conjugacy class $\left\langle h\right\rangle $  of $G$, Lott's \emph{delocalized eta invariant} 	$\eta_{\left\langle h\right\rangle }(\widetilde D)$ of $\widetilde D$ is defined to be
	\begin{equation}\label{delocalizedeta}
	\eta_{\left\langle h\right\rangle }(\widetilde D)\coloneqq 
	\frac{2}{\sqrt\pi}\int_{0}^\infty \tr_{\left\langle h\right\rangle }(\widetilde De^{-t^2\widetilde D^2})dt.
	\end{equation}
\end{definition}
 The terminology ``delocalized" refers to the fact that $\langle h\rangle$ is a nontrivial conjugacy class. If we were to take the trivial conjugacy class in Definition $\ref{def:deloceta}$, we would recover the $L^2$-eta invariant of Cheeger and Gromov \cite{MR806699}.

Lott proved  the convergence of the integral in line $\eqref{delocalizedeta}$ under the assumption that $G$ has polynomial growth or is hyperbolic, and $\widetilde  D$ is invertible (or more generally has a spectral gap at zero) \cite{Lott}.
Piazza and Schick gave an example where the formula \eqref{delocalizedeta} diverges for non-invertible $\widetilde D$ \cite[Section 3]{Piazza}. They then raised the question of whether a divergent example still exists if one assumes the invertibility of  $\widetilde D$ \cite[Remark 3.2]{Piazza}. 

Our first main result states that the answer to this question of Piazza and Schick is negative, as long as the spectral gap of $\widetilde D$ is sufficiently large. Before we give the precise statement of the theorem, we first fix some notation.

Fix a finite generating set $S$ of $G$. Let  $\ell$ be the corresponding word length function on $G$ determined by $S$. Since $S$ is finite, there exist $C>0$ and $K_{\langle h\rangle}>0$ such that 
\begin{equation}\label{eq:growthofh}
\#\{g\in\langle h\rangle:\ell(g)=n\}\leqslant C e^{K_{\langle h\rangle}\cdot n}.
\end{equation}
We define $\uptau_{\langle h\rangle}$ to be 
\begin{equation}\label{eq:equivalentwithlength}
\uptau_{\langle h\rangle}=\liminf_{\substack{g\in\langle h\rangle \\ \ell(g)\to\infty}} \Big(\inf_{x\in\mathcal F}\frac{\dist(x,gx)}{\ell(g)}\Big).
\end{equation}
Since  the action of $G$ on $\widetilde M$ is free and cocompact, we have $\uptau_{\langle h\rangle} > 0$.

Given a the first-order differential operator $D$ on $M$, we denote the principal symbol of $D$ by $\sigma_D(x, v)$, for $x\in M$ and cotangent vector $v\in T_x^\ast M$. 
	We define the propagation speed of $D$ to be the  positive number 
	\begin{equation}\label{eq:propagationspeed}
	c_D = \sup\{ \|\sigma_D(x, v)\| \colon x\in M, v\in T^\ast_xM, \|v\| =1\}. 
	\end{equation} 
When $D$ is the Dirac operator on a spin manifold, we have $c_D = 1$.
\begin{definition}\label{def:speclarge}
With the above notation, let us define	\begin{equation}\label{eq:sigma_h}
	\sigma_{\langle h\rangle}\coloneqq \frac{2K_{\langle h\rangle}\cdot c_D}{\uptau_{\langle h\rangle}}.
	\end{equation}
We say the spectral gap of $\widetilde D$ is sufficiently large if  the spectral gap of $\widetilde D$ at zero is larger than $\sigma_{\langle h\rangle}$, i.e. $\textup{spectrum}(\widetilde D) \cap [-\sigma_{\langle h\rangle}, \sigma_{\langle h\rangle}]$ is either $\{0\}$ or empty.
\end{definition}

In the following, we shall show that the convergence of the formula \eqref{delocalizedeta}  holds if  $\widetilde D$ has a spectral gap at zero larger than $\sigma_{\langle h\rangle}$.  In fact, the convergence of the formula \eqref{delocalizedeta} for the case where $\widetilde  D$ has a spectral gap at zero can be deduced from the invertible case by replacing $\widetilde  D$ with its restriction to the orthogonal complement of the kernel of $\widetilde  D$. Without loss of generality, we will only give the details of the proof for the case where $\widetilde D$ is invertible and its spectral gap at zero is larger than $\sigma_{\langle h\rangle}$.

\begin{theorem}\label{convergence}
	 With the same notation as above, for any nontrivial conjugacy class $\left\langle h\right\rangle $, suppose that $\widetilde D$ is invertible and its  spectral gap at zero is larger than $\sigma_{\langle h\rangle}$. Then the delocalized eta invariant $\eta_{\left\langle h\right\rangle }(\widetilde D)$ given in line \eqref{delocalizedeta} converges absolutely. 
\end{theorem}
\begin{proof}
    The proof is divided into three steps. In the first step, we show that $\tr_{\left\langle h\right\rangle }(\widetilde De^{-t^2\widetilde D^2})$ is finite for any fixed $t>0$ (Proposition \ref{convergenceofTRforA}).  In the second step, we prove the convergence of the integral for small $t$  (Proposition \ref{convergenceofetanearzero}). In the last step, we show the convergence of the integral for large $t$ (Proposition \ref{convergencofetanearinf}). 
  In fact,  only the last step requires the assumption that  the spectral gap of $\widetilde  D$ is larger than $\sigma_{\langle h\rangle}$. 
\end{proof}
We say that $\langle h\rangle$ has sub-exponential growth if we can choose $K_{\langle h\rangle}$ in line \eqref{eq:growthofh} to be arbitrarily small. In this case, we have the following corollary.
\begin{corollary}
	If $\langle h\rangle$ has sub-exponential growth, then $\eta_{\left\langle h\right\rangle }(\widetilde D)$ given in line \eqref{delocalizedeta} converges absolutely,  as long as $\widetilde D$ has a spectral gap at zero. 
\end{corollary}

In Theorem $\ref{convergence}$, the condition that the spectral gap of $\widetilde D$ is larger than $\sigma_{\langle h\rangle}$  might first appear to be rather ad hoc. In the following, we shall show that such a condition in fact holds for an abundance of natural examples whose higher rho invariant is nontrivial. 
	
	Suppose that $N$ is a closed spin manifold  equipped with a positive scalar curvature metric $g_N$, whose fundamental group $F = \pi_1(N)$ is finite and its higher rho invariant $\rho(\widetilde D_N)$ is nontrivial. Here $\widetilde D_N$  is the Dirac operator on the universal covering $\widetilde N$ of $N$. For instance, let  $N$ to be a lens space, that is, the quotient of the 3-dimensional sphere by a free action of a finite cyclic group. In this case, the classical equivariant Atiyah-Patodi-Singer index theorem implies that the higher rho invariant of $N$ is nontrivial, cf. \cite{MR511246}.   
	
	Now let $V$ be an even dimensional closed spin manifold, whose Dirac operator $D_V$ has nontrivial higher index in $K_0(C^\ast_r(\Gamma))$, where $\Gamma = \pi_1(V)$. In particular, it follows that $D_V$ defines a nonzero element in  the equivariant \mbox{$K$-homology} $K_0(C^\ast_L(E\Gamma)^\Gamma)$ of the universal space $E\Gamma$ for free $\Gamma$ actions.  Consider the product  space $M = V\times N$  equipped with a metric $g_{M}=g_V+\varepsilon\cdot  g_N$, where $g_V$ is an arbitrary Riemannian metric on $V$ and the metric $g_N$ on $N$ is scaled by a positive number $\varepsilon$. Denote the Dirac operator on the universal covering $\widetilde M$ of $M$ by $\widetilde D_M$.  The product formula for secondary invariants (cf. \cite[Claim 2.19]{Xiepos}, \cite[Corollary 4.15]{MR3551834}) shows that the higher rho invariant $\rho(\widetilde D_M)$ is the product of the $K$-homology class of $D_V$ and the higher rho invariant $\rho(\widetilde D_N)$. By the assumptions above, we see that $\rho(\widetilde D_M)$ is nonzero in  $K_1(C_{L,0}^\ast(E(\Gamma\times F))^{\Gamma\times F})$.

Now let us return to the proof of Theorem $\ref{convergence}$. First, we need a few technical lemmas. Let $M$ be a closed Riemannian manifold. 	Suppose $T$ is an integral operator on the space $L^2(M)$ of $L^2$ functions on $M$, and assume the Schwartz kernel of $T$ is a  continuous function on $M\times M$. In particular, we have 
$$T(f)(x)=\int_M T(x,y)f(y)dy,$$
for all $f\in L^2(M)$. 
We denote  the operator norm of an operator $A$ on $L^2(M)$ by $\|A\|_{op}$.

\begin{lemma}\label{kernelupperbound}
	 Let $M$ be a closed Riemannian manifold and $D$ a first-order self-adjoint elliptic differential operator on $M$. Suppose $T$ is a bounded linear operator on $L^2(M)$ such that 
	 \[ \sup_{k+j \leqslant \frac{3}{2}\dim M+3}\|D^{k}TD^{j}\|_{op} <\infty. \] Then   $T$ is an integral operator with a continuous Schwartz kernel $K_T(x,y)$, and there exists a positive number $C$ $($independent of $T)$  such that 
	\begin{equation}
	\sup_{x,y\in M}|K_T(x,y)|\leqslant C\cdot \sup_{k+j \leqslant \frac{3}{2}\dim M+3}\|D^{k}TD^{j}\|_{op}.
	\end{equation}

\end{lemma}
\begin{proof}
	Let $p$ be the smallest even integer that is greater than $\frac{1}{2}\dim M$. Then $(1 + D^{p})^{-1}$ is a Hilbert-Schmidt operator on $L^2(M)$. Denote the Hilbert-Schmidt norm of $(1 + D^{p})^{-1}$ by $\|(1 + D^{p})^{-1}\|_{HS}$, 
	
	By assumption, $(1 + D^{p})T$ is bounded. It  follows that $T = (1 + D^{p})^{-1}\circ (1 + D^{p})T$ is also a Hilbert-Schmidt operator, and furthermore 
	\[  \|T\|_{HS} \leq \|(1 + D^{p})^{-1}\|_{HS} \cdot \|(1 + D^{p})T\|_{op}.\] 	
	It follows that $T$ is an integral operator whose Schwartz kernel $K_T(x, y)$ is an \mbox{$L^2$-function} on $M\times M$.  
	
	To see that $K_T$ is continuous on $M\times M$, we consider the elliptic differential operator 
	\[ \mathscr D = D\otimes 1 + 1\otimes D \]
	on $M\times M$. Observe that 
\begin{align*}
\|\mathscr D^n K_T\|_{L^2(M\times M)} & = \Big\|\sum_{r=0}^{n} {n \choose r} D^r TD^{n-r}\Big\|_{HS} \\
& \leq \sum_{r=0}^{n} {n \choose r} \|(1 + D^{p})^{-1}\|_{HS} \cdot \|(1 + D^{p})D^r TD^{n-r}\|_{op}.
\end{align*} 
Therefore, our assumption 
\[ \sup_{k+j \leqslant \frac{3}{2}\dim M+3}\|D^{k}TD^{j}\|_{op} <\infty \] 
implies that $\|\mathscr D^n K_T\|_{L^2(M\times M)}$ is finite for all $n\leq \dim M + 1$. 
It follows from the Sobolev embedding theorem that there exists $C_1 >0$ such that 
\begin{align*}
	\sup_{x,y\in M}|K_T(x,y)|& \leqslant C_1 \cdot \sup_{n\leq \dim M + 1}  \|\mathscr D^n K_T\|_{L^2(M\times M)}  
\end{align*}
where the right hand side is dominated by 
\[  C\cdot \sup_{k+j \leqslant \frac{3}{2}\dim M+3}\|D^{k}TD^{j}\|_{op}  \]
for some constant $C >0$. This finishes the proof. \end{proof}

\begin{remark}
	Let $E$ be a Hermitian vector bundle over $M$ and $\widetilde M$ a regular \mbox{$G$-covering} space of $M$. Denote the lift of $E$ to $\widetilde M$ by $\widetilde E$. The above lemma admits an obvious analogue for $G$-equivariant operators $T$ acting on $L^2(\widetilde M, \widetilde E)$.
\end{remark}

\begin{remark}
As an immediate consequence of the above lemma, we see that for a closed Riemannian manifold $M$ and any bounded linear operator $T$ on $L^2(M)$, if
	\[ \sup_{k+j\leqslant \frac{3}{2}\dim M+3}\|D^{k}TD^{k'}\|_{op} < \infty,\] then $T$ is of trace class. In fact, in this case, we have 
	$$\tr(T)=\int_M T(x,x)dx,$$
	cf. \cite[Chapter V, Proposition 3.1.1]{MR0444836}.
	
\end{remark}

	Suppose $D$ is a first-order self-adjoint elliptic differential operator acting on a vector bundle $E$ over $M$, and $\widetilde E$ (resp. $\widetilde D$) is the lift of $E$ (resp. $D$) to $\widetilde M$.
	If $f$ is a function on $\mathbb R$  such that  
	\begin{equation}\label{3.8}
	\|x^{m}f(x)\|_{L^\infty}<\infty
	\end{equation} for all $m\leqslant \frac{3}{2}\dim M+3$, then the corresponding Schwartz kernel of the operator  $f(\widetilde D)$ is  continuous. Denote the Schwartz kernel of  $f(\widetilde D)$ by $K_f$.
	\begin{lemma}\label{functionalcal}
	  With the notations above, for  any $\mu>1$ and $r>0$,  there exists a constant $C>0$ such that 
	\begin{equation*}
	\|K_f(x,y)\|\leqslant C\cdot F_f\left(\frac{\dist(x, y)}{\mu \cdot c_D}\right), 
	\end{equation*}
for $\forall x,y\in \widetilde M$ with  $\dist(x,y)>r$ and any $f$ satisfying line \eqref{3.8}. Here $\dist(x, y)$ stands for the distance between $x$ and $y$, the notation $\|\cdot\|$ denotes a matrix norm of elements in $\mathrm{End}(\widetilde E_y,\widetilde E_x)$, and the function $F_f$ is defined by 
	\begin{equation*}
	F_f(s)\coloneqq \sup_{n\leqslant \frac{3}{2}\dim M+3}\int_{|\xi|>s}\left|\frac{d^n}{d\xi^n}\widehat f(\xi)\right|d\xi, 
	\end{equation*}
	where $\widehat f$ is the Fourier transform of $f$. 	
\end{lemma}	
\begin{proof}
	The condition on $f$ implies that  $F_f(s)<\infty$ for any $s\in \mathbb R$ and $f(\widetilde  D)$ is an integral operator with continuous Schwartz kernel (cf. Lemma \ref{kernelupperbound}).
	
	By the Fourier inverse transform formula, we have
	\begin{equation*}
	f(\widetilde D)=\frac{1}{2\pi}\int_{-\infty}^{+\infty} \widehat f(\xi)e^{i\xi\widetilde D}d\xi. 
	\end{equation*}
	Fix $r>0$. Let $x_0, y_0\in \widetilde M$ such that  $\lambda \coloneqq \dist(x_0, y_0) \geqslant r$.    Choose a smooth function $\varphi$ over $\mathbb R$ such that 
	$\varphi(\xi)=1$ for $|\xi|\geqslant 1 $ and $\varphi(\xi)=0$ for $|\xi|\leqslant 1/\mu$.  
	Let $\varphi_\lambda(\xi)=\varphi(\xi\cdot c_D/\lambda)$.
	Let $g$ be the function with Fourier transform
	  \[ \widehat{g}(\xi)=\varphi_\lambda(\xi)\widehat f(\xi),  \]
	and  $L\in C(\widetilde M \times \widetilde M)$ be the Schwartz kernel corresponding to the operator $g(\widetilde  D)$. It follows from standard finite propagation estimates of wave operators that  $L(x,y)=K_f(x,y)$ for all $x, y \in \widetilde M$ with $\dist(x,y)\geqslant \lambda$. In particular, we have  $L(x_0, y_0) = K_f(x_0, y_0)$.
	
	By Lemma \ref{kernelupperbound}, it suffice to estimate the operator norm of 
	\[  \widetilde D^{k}g(\widetilde D)\widetilde D^{j}=\widetilde D^{k+j}g(\widetilde D) \]
	for all $k+j \leqslant \frac{3}{2}\dim M+3$. 
Now for a given $n\leqslant \frac{3}{2}\dim M+3$,  define $\psi_n(x)=x^{n}g(x)$. We have
	\begin{equation*}
	\widehat{\psi}_n(\xi)=\left(\frac{1}{i}\frac{d}{ d\xi}\right)^n(\varphi_\lambda\,\widehat f)(\xi).
	\end{equation*}
	Since $\varphi$ is supported on $|\xi|\geqslant \lambda/(\mu c_D)$, there exist positive numbers $C_1$ and $C$ such that 
	\begin{align*}
	\|\psi_n(\widetilde  D)\|_{op}&\leqslant 
	\frac{1}{2\pi}\int_{ |\xi|\geqslant \lambda/(\mu c_D)}|\widehat\psi_n(\xi)|d\xi\\
	&\leqslant C_1 \sum_{j=0}^n\left( \frac{c_D}{r}\right) ^j\int_{ |\xi|\geqslant \lambda/(\mu c_D)}|\widehat f^{(n-j)}(\xi)|d\xi\\
	&\leqslant C\cdot F_f\left(\frac{\lambda}{\mu\cdot c_D}\right)
	\end{align*}
	By Lemma $\ref{kernelupperbound}$, we see that
	$$\|K_f(x_0,y_0)\|=\|L(x_0,y_0)\|\leqslant C\cdot F_f\left(\frac{\dist(x_0, y_0)}{\mu \cdot  c_D}\right).$$
	This finishes the proof. 
\end{proof}

To streamline our estimates later, we consider the following class of functions in $C_0(\mathbb R)$ whose Fourier transform has exponential decay. 

\begin{definition}\label{ALambdaN}
	Let $\mathcal A_{\Lambda,N}$ be the subspace of $C_0(\mathbb R)$ consisting of functions $f$ satisfying the following conditions:
	\begin{enumerate}
		\item $f$ admits an analytic continuation $\widetilde f$ on the strip $\{|\mathrm{Im}(z)|<\Lambda\}$,
		\item for any $n\leqslant N$, $|z^n\widetilde f(z)|$ is uniformly bounded on the strip.
	\end{enumerate}
	
	Equip $\mathcal A_{\Lambda,N}$ with the norm $\|\cdot\|_{\mathcal A}$ defined by
	$$\|f\|_\mathcal A=\sup_{0\leqslant n\leqslant N}\sup_{|\mathrm{Im}(z)|<\Lambda}|z^n\widetilde f(z)|$$
\end{definition}

For any fixed $\Lambda$ and $N$, it is easy to verify that $\mathcal A_{\Lambda,N}$ is closed under multiplication and conjugation. In fact,  $\mathcal A_{\Lambda,N}$ is a Banach $*$-subalgebra of $C_0(\mathbb R)$ under the norm $\|\cdot\|_\mathcal A$. Clearly,  $e^{-x^2}$ and $xe^{-x^2}$ lie in $\mathcal A_{\Lambda,N}$. 

\begin{lemma}\label{estimateA}
	Suppose $f\in \mathcal A_{\Lambda,N}$ for some $N\geqslant 2$.
	Let $\widehat f$ be its Fourier transform. Then 
	for any $0<\lambda<\Lambda$ and $0\leqslant n\leqslant N-2$ there exists some constant $C = C_{\lambda, n}$ such that 
	\begin{equation}
	\int_{|\xi|>s}|\frac{d^n}{d\xi^n}\hat f(\xi)|d\xi \leqslant C\cdot \|f\|_{\mathcal A}\cdot e^{-\lambda s}.
	\end{equation}	
\end{lemma}
\begin{proof}
	For notational simplicity, let us denote the analytic continuation of $f$ on the strip $\{|\mathrm{Im}(z)|<\Lambda\}$ still by $f$. For any $|y|<\Lambda$, $ f(x-iy)$ is a smooth $L^1$-integrable function since $|z^2  f(z)|$ is uniformly bounded. Denote the Fourier transform of $f(x-iy)$ with respect to $x$ by $\widehat{f_y}$, namely
	\begin{equation*}
	\widehat{f_y}(\xi)=\int_{-\infty}^{+\infty}  f(x-iy)e^{-i\xi x}dx
	\end{equation*}
	
	Since the right-hand side is differentiable in $y$, and uniformly differentiable in $x$ and $y$ by Cauchy inequality, the left-hand side is also differentiable in $y$ with
	\begin{equation*}
	\frac{\partial}{\partial y}\widehat{f_y}(\xi)=\int_{-\infty}^{+\infty} \frac{\partial}{\partial y}[f(x-iy)]e^{-i\xi x}dx
	\end{equation*}
	Since $f$ is holomorphic, we have that $\frac{\partial}{\partial y}f=i\frac{\partial}{\partial x}f$ by the Cauchy-Riemann equation. Thus
	\begin{equation*}
	\frac{\partial}{\partial y}\widehat{f_y}(\xi)=\int_{-\infty}^{+\infty} i\frac{\partial}{\partial x}[  f(x-iy)]e^{-i\xi x}dx
	=\xi \widehat{f_y}(\xi)
	\end{equation*}
	It follows that
	\begin{equation*}
	\widehat{f_y}(\xi)=\widehat{f}(\xi)e^{y\xi}
	\end{equation*}
	
	Therefore by our assumption, for any $n\leqslant N-2$ and $0<\lambda<\Lambda$ there exists some constant $C_1>0$ such that 
	\begin{equation*}
	\begin{split}
	|\frac{d^n}{d\xi^n}(\widehat f(\xi)e^{\lambda\xi})|&=|\frac{d^n}{d\xi^n}\widehat{f}_y(\xi)|
	=\int_{-\infty}^{+\infty}e^{ix\xi} (ix)^n f(x-i\lambda)dx \\
	&\leqslant \int_{-\infty}^{\infty}|x^n f(x-i\lambda)|dx\leqslant C_1\|f\|_{\mathcal A}
	\end{split}
	\end{equation*}
	Thus by induction on $n$, for any $n\leqslant N-2$ and $0<\lambda<\Lambda$,  there exists a constant $C_2>0$ such that 
	\begin{equation*}
	|\frac{d^n}{d\xi^n}(\widehat f(\xi))|\leqslant C_2\|f\|_{\mathcal A} e^{-\lambda\xi}.
	\end{equation*}
	Hence for $\xi > s$, we have
	\begin{equation*}
	\int_{\xi>s}|\frac{d^n}{d\xi^n}\widehat f(\xi)|d\xi \leqslant Ce^{-\lambda s}.
	\end{equation*}
	The estimates for the part $\int_{\xi <  -s}|\frac{d^n}{d\xi^n}\widehat f(\xi)|d\xi$ are completely similar. This finishes the proof. 
\end{proof}

If we fix a fundamental domain $\mathcal F\subset \widetilde M$ for the action of $G$, then one naturally identifies  $L^2(\widetilde M)$ with $L^2(\mathcal F) \otimes \ell^2(G)$ through the mapping $\tilde h \mapsto h$ by $h(x, \gamma) = \tilde h(\gamma x)$ for $x\in \mathcal F$ and $\gamma\in G$. In particular, every $G$-equivariant Schwartz kernel $A$ on $\widetilde M\times \widetilde M$ becomes a formal sum 
\[  A = \sum_{g\in G}A_g R_g \]
where  $A_g(x, y) = A(x, gy)$ for $x, y \in \mathcal F$ and $R_g$ denotes the right translation of $g$ on $\ell^2(G)$ corresponding to the right regular representation of $G$ on $\ell^2(G)$. Now suppose $A = f(\widetilde D)$ as in Lemma $\ref{functionalcal}$ above. In this case, each $A_g$ is a trace class operator and 
\[ \tr(A_g) = \int_{\mathcal F} A_g(x, x) dx = \int_{\mathcal F} A(x, gx) dx. \]

Now we are ready to proceed with the first step of the proof for Theorem $\ref{convergence}$.

\begin{proposition}\label{convergenceofTRforA}
	Suppose $\langle h\rangle$ is a nontrivial conjugacy class of $G$. If $f\in\mathcal A_{\Lambda,N}$ with $N\geqslant \frac{3}{2}\dim M+5$ and $\Lambda$ sufficiently large, then $\tr_{\langle h\rangle}(f(\widetilde{D}))$ is finite.
\end{proposition}
\begin{proof}
	Fix a symmetric generating set $S$ of $G$. Let $\ell$ be the length function on $G$ determined by $S$. Since $G$ acts freely and cocompactly on $\widetilde M$, there exists $\varepsilon>0$ such that $\dist(x,gx)>\varepsilon$ for any $x\in \widetilde M$ and $g\in\langle h\rangle$. 
	Let $k(x,y)$ be the Schwartz kernel of $f(\widetilde D)$. Since the action of $G$ on $\widetilde M$ is free and cocompact,  there exist $C_1,C_2 >0$ such that 
	\[ \dist(x, gx) > C_1 \cdot \ell(g) - C_2. \]
	It follows from Lemma \ref{functionalcal} that there exists $C_3>0$ such that 
	\begin{align*}
	|\tr_{\langle h\rangle}(f(\widetilde{D}))|&\leqslant\sum_{g\in\langle h\rangle}\int_{\mathcal F}|\tr(k(x,gx))|dx
	\\
	&\leqslant C_3\sum_{g\in\langle h\rangle}F_f\left(\max\{\varepsilon,C_1\cdot\ell(g)-C_2\}\right)\\
	&\leqslant C_3\sum_{n=1}^\infty |S|^nF_f\left(\max\{\varepsilon,C_1\cdot n-C_2\}\right),
	\end{align*}
	where $|S|$ is the cardinality of the generating set $S$. By Lemma \ref{estimateA}, when $N\geqslant \frac{3}{2}\dim M+5$, for any $\lambda<\Lambda$, there exists $C>0$ such that
	$$F_f(x)\leqslant Ce^{-\lambda x}.$$
	Observe that the summation
	$$\sum_{n=1}^\infty |S|^ne^{-\lambda\cdot \max\{\varepsilon, \, C_2\cdot n-C_3\}}$$
	 converges absolutely, as long as $\lambda$ is sufficiently large.  This finishes the proof. 
\end{proof}

Now let us prove the convergence of the integral \eqref{delocalizedeta} for small $t$. 
\begin{proposition}\label{convergenceofetanearzero}
Suppose $\langle h\rangle$ is a nontrivial conjugacy class of $G$. If $f\in\mathcal A_{\Lambda,N}$ with $N\geqslant \frac{3}{2}\dim M+5$ and $\Lambda$ sufficiently large, then the following integral
	$$\int_0^1 t^k\cdot \tr_{\langle h\rangle}(f(t\widetilde D)) dt$$
	is absolutely convergent for any $k\in \mathbb R$. 
\end{proposition}
\begin{proof}
	Let us denote $f_t(x)=f(tx)$. Clearly,  $f_t\in\mathcal A_{\Lambda,N}$ since $f\in\mathcal A_{\Lambda,N}$.
	Similar to the proof of Proposition $\ref{convergenceofTRforA}$, we have 
	\[ |t^k \tr_{\langle h\rangle}(f_t(\widetilde{D}))|\leqslant C_3\sum_{n=1}^\infty |S|^n t^k F_{f_t}\left(\max\{\varepsilon,C_1\cdot n-C_2\}\right). \]
	Recall that  $\widehat{f_t}(\xi)=t^{-1}\cdot \widehat{f}(\xi/t)$. In particular,  we have 
	$$\frac{d ^n}{d\xi^n}\widehat{f_t}(\xi)=\frac{1}{t^{n+1}}\widehat{f}^{(n)}(\xi/t),$$
	where $\widehat f^{(n)}$ is the $n$-th derivative of $\widehat f$. 
	It follows from Lemma $\ref{estimateA}$ that there exists $C_\Lambda >0$ such that 
	\begin{equation}\label{eq:F_fwitht}
	\begin{split}
	F_{f_t}(s)=&\sup_{n\leqslant \frac{3}{2}\dim M+3}\frac{1}{t^{n+1}}\int_{|\xi|>s}\Big|
	\widehat f^{(n)}(\xi/t)\Big|d\xi\\
	\leqslant& \frac{F_f(s\cdot t^{-1})}{t^{\frac{3}{2}\dim M+3}}
	\leqslant \frac{C_\Lambda}{t^{\frac{3}{2}\dim M+3}} \cdot  \|f\|_{\mathcal A}\cdot \exp(\frac{-\Lambda s}{2t})
	\end{split}
	\end{equation}
	for all $t\in (0, 1]$. The following summation
	$$\sum_{n=1}^\infty |S|^n {t^{k-(\frac{3}{2}\dim M+3)}} \exp\Big(\frac{-\Lambda\cdot \max\{\varepsilon,C_1\cdot n-C_2\} }{2t}\Big)$$
	is integrable on $(0, 1]$, as long as $\Lambda$ is sufficiently large. This finishes the proof.
\end{proof}

\begin{proposition}\label{convergencofetanearinf}
	    Let $\sigma^2$ be the infimum of the spectrum of $\widetilde D^2$.
	    If $\sigma>\sigma_{\langle h\rangle}$ defined in line \eqref{eq:sigma_h}, then the following integral
	     $$\int_1^{+\infty} \tr_{\langle h\rangle}(\widetilde  De^{-t^2\widetilde  D^2}) dt$$
	     is absolutely convergent.
\end{proposition}
\begin{proof}
	
	View $\widetilde  De^{-t^2\widetilde  D^2}$ as an element in $\mathcal K\otimes C_r^*(G)$ and write
	$$\widetilde  De^{-t^2\widetilde  D^2}=\sum_{g\in G} A_{g,t} g.$$
Note that $A_{g,t}$ are compact operators for all $g\in G$ and $t\geqslant 1$. 
	By Lemma \ref{kernelupperbound}, $A_{g,t}$ is of trace class and there exist $C_1,N_1>0$ such that
	\begin{equation*}
	|A_{g,t}|_1\leqslant C_1t^{N_1}\cdot e^{-t^2\sigma^2},
	\end{equation*}
	where $|\cdot|_1$ stands for the trace norm.
	
By the definition of $\uptau_{\langle h\rangle}$ in line \eqref{eq:equivalentwithlength}, for any $\varepsilon>0$ there exists $L>0$ such that 
	\[ \dist(x,gx)\geqslant (\uptau_{\langle h\rangle}-\varepsilon)\ell(g),  \]
	for all $x\in\mathcal F$ and $g\in \langle h \rangle$ satisfying $\ell(g)>L$. 	By Lemma \ref{functionalcal} and Proposition \ref{convergenceofTRforA}, if $\ell(g)>L$, then there exist $C_2,N_2>0$ such that 
	\begin{equation}
	|A_{g,t}|_1\leqslant C_2\left(t\cdot \ell(g)\right)^{N_2}\cdot\exp\left(-\frac{(\uptau_{\langle h\rangle}-\varepsilon)^2 \ell^2(g)}{4t^2(\mu\cdot  c_D)^2}\right).
	\end{equation}
	Write $N=\max\{N_1,N_2\}$ and $C=\max\{C_1,C_2\}$. Then for $\delta>0$ sufficiently small, we have 
	\begin{align*}
		|A_{g,t}|_1\leqslant &C\left(t\cdot \ell(g)\right)^{N}\cdot\exp\left(-\frac{1}{2}\frac{(\uptau_{\langle h\rangle}-\varepsilon)^2 \ell^2(g)}{4t^2(\mu\cdot  c_D)^2}-\frac{t^2\sigma^2}{2}\right)\\
		\leqslant &	C\left(t\cdot \ell(g)\right)^{N}\cdot \exp\left(-\frac{(\uptau_{\langle h\rangle}-\varepsilon)(\sqrt{\sigma^2-\delta^2})\ell(g)}{2\mu\cdot  c_D}\right)\cdot e^{-t^2\delta^2}, 
	\end{align*}
	 as long as $\ell(g) >L$.
	By the assumption that $\sigma>\sigma_{\langle h\rangle}$, we may find suitable $\mu,\varepsilon,\delta$ such that $$(\uptau_{\langle h\rangle}-\varepsilon)(\sqrt{\sigma^2-\delta^2})>2\mu\cdot  c_D \cdot K_{\langle h\rangle}.$$
	Therefore, by line \eqref{eq:growthofh}, there exist $C_3,C_4$ such that
	\begin{align*}
	|\tr_{\langle h\rangle}(\widetilde  De^{-t^2\widetilde  D^2})| \leqslant&
	\sum_{g\in\langle h\rangle,\ell(g)\leqslant L}|A_{g,t}|_1+
	\sum_{g\in\langle h\rangle,\ell(g)>L}|A_{g,t}|_1\\
	\leqslant& C_3e^{K_{\langle h\rangle}L}\cdot C(tL)^Ne^{-t^2\sigma^2}\\&+\sum_{n>L}C_3e^{K_{\langle h\rangle}n}\cdot C\left(nt\right)^{N}\exp\left(-\frac{(\uptau_{\langle h\rangle}-\varepsilon)(\sqrt{\sigma^2-\delta^2})n}{2\mu c_D}\right)e^{-t^2\delta^2}\\
	\leqslant & C_4t^N e^{-\delta^2 t^2}.
	\end{align*}
	Hence the integral
	$$\int_0^\infty\tr_{\langle h\rangle}(\widetilde  De^{-t^2\widetilde  D^2})  dt $$
	converges absolutely. This finishes the proof.
\end{proof}
\subsection{Delocalized higher eta invariants}\label{HigherEtaFormula}
In this subsection, we shall generalize the results of the previous subsection to higher degree cyclic cocyles. For simplicity, we only give details for the odd case;  the even case is similar.

Let us first recall the definition of cyclic cocycles. 
\begin{definition}\label{cycliccohomology}
	Let $C^n(\mathbb CG)$ be the space spanned by all $(n+1)$-linear functionals $\varphi$ on $\mathbb CG$ such that
	$$\varphi(g_1,g_2,\cdots,g_n,g_0)=(-1)^n\varphi(g_0,g_1,\cdots,g_n).$$
	The coboundary map $b:C^n(\mathbb CG)\to C^{n+1}(\mathbb CG)$ is defined to be
	\begin{align*}
	b\varphi(g_0,g_1,\cdots,g_{n+1}) & =\sum_{j=0}^n(-1)^j\varphi(g_0,g_1,\cdots,g_jg_{j+1},\cdots,g_{n+1})\\
	&\hspace{1.5cm} +
	(-1)^{n+1}\varphi(g_{n+1}g_0,g_1,\cdots,g_n). 
	\end{align*}
\end{definition}
	The cohomology of this cochain complex $(C^n(\mathbb CG),b)$ is the cyclic cohomology of $\mathbb CG$, denoted by $HC^*(\mathbb CG)$.
\begin{definition}
Fix any conjugacy class $\left\langle h\right\rangle $ of $G$. Let $C^n(\mathbb CG,\left\langle h\right\rangle )$ be the space spanned by all elements $\varphi\in C^n(\mathbb CG)$ satisfying the condition:
$$g_0g_1\cdots g_n\notin\left\langle h\right\rangle\implies \varphi(g_0,g_1,\cdots,g_n)=0.$$
\end{definition}
It is easy to verify that that $(C^n(\mathbb CG,\left\langle h\right\rangle ),b)$ is a subcomplex of $(C^n(\mathbb CG),b)$.
We denote the cohomology of $(C^n(\mathbb CG,\left\langle h\right\rangle ),b)$ by $HC^*(\mathbb CG,\left\langle h\right\rangle)$.

\begin{definition}\label{def:delocalcocyc}
	If the conjugacy class $\langle h\rangle$ is nontrivial, then a cyclic cocycle in $C^n(\mathbb CG,\left\langle h\right\rangle )$ is called a \emph{delocalized cyclic cocycle} at $\langle h\rangle$, and  $HC^*(\mathbb CG,\left\langle h\right\rangle)$ is called the \emph{delocalized cyclic cohomology} of $\mathbb CG$ at $\langle h\rangle$. 
\end{definition}
Recall that  (cf. \cite{Nistor})
$$HC^*(\mathbb CG)\cong\prod_{\left\langle h\right\rangle }HC^*(\mathbb CG,\left\langle h\right\rangle).$$
Moreover, it is easy to verify that $HC^0(\mathbb CG,\langle h\rangle)$ is a one dimensional vector space generated by $tr_{\langle h\rangle}$.

Let us write $w =\sum_g w^g g$ for an element $w\in  C^*(\widetilde  M)^G$.  Given $\varphi\in C^n(\mathbb CG,\left\langle h\right\rangle )$ and $w=w_0\otimes w_1\otimes\cdots\otimes w_n$ in $(C^*(\widetilde  M)^G)^{\otimes n+1}$, we define the following map 
	\begin{equation}\label{phiTR}
	(\varphi\# \tr)(w)=\sum_{g_0,\cdots,g_n\in G} \tr(w_0^{g_0}\cdots w_n^{g_{n}})\varphi(g_0,\cdots,g_n)
	\end{equation}
	whenever the above formula converges.
	Here $\tr(w_0^{g_0}\cdots w_n^{g_{n}})$ stands for the trace of  the  operator $w_0^{g_0}\cdots w_n^{g_{n}}$. 

Following the definition of higher rho invariant in line \eqref{higherrho}, we define 
    \begin{equation}\label{eq:smallu}
        u_t(x)=U_{1/t}(x)=\exp\Big(2\pi i\frac{1}{\sqrt{\pi}}\int_{-\infty}^{xt} e^{-s^2}ds\Big).
    \end{equation}
Note that the functions $u_t(x)-1$, $u_t^{-1}(x)-1$ and $\dot{u}_t(x)u^{-1}_t(x)=2\sqrt \pi ix e^{-t^2x^2}$ are Schwartz functions. 
We define the delocalized higher eta invariant as follows.
\begin{definition}\label{higherdelocalizedeta}
    For any $\varphi\in C^{2m}(\mathbb CG,\left\langle h\right\rangle )$ with $\left\langle h\right\rangle $ nontrivial,  we define the delocalized higher eta invariant of $\widetilde D$ with respect to $\varphi$  to be
    \begin{equation}\label{etahigher}
    \eta_{\varphi}(\widetilde D):=\frac{m!}{\pi i}\int_0^\infty \eta_{\varphi}(\widetilde D,t)dt,
    \end{equation}
    where 
	\begin{equation}\label{etaoftensor}
	\eta_\varphi(\widetilde D,t)=\varphi\#\tr(
	\dot{u}_t(\widetilde D) u_t^{-1}(\widetilde D)\otimes ((u_t(\widetilde D)-1)\otimes (u_t^{-1}(\widetilde D)-1))^{\otimes m}).
	\end{equation}
	More precisely, we have 
	\begin{equation}\label{eta(t)}
\begin{split}
\eta_{\varphi}(\widetilde D,t)& \coloneqq  \sum_{g_i\in G}\Big[\varphi(g_0,g_1,\cdots,g_{2m})\int_{\mathcal F^{2m+1}} \tr\big(k_{0,t}(x_0,g_0x_1) k_{1,t}(x_1,g_1x_2) \\
&\hspace{4.5cm}  \cdots k_{2m,t}(x_{2m},g_{2m}x_0)\big)dx_0\cdots dx_{2m}\Big].
\end{split}
\end{equation}
where  $k_{i,t}(x,y)$ is the corresponding Schwartz kernel of $\dot{u}_t(\widetilde D)u^{-1}_t(\widetilde D)$, $u_t(\widetilde D)-1$ and $u_t^{-1}(\widetilde D)-1$ respectively.

\end{definition}

\begin{remark}
	Clearly, if $m=0$, then $\eta_{\tr_{\langle h\rangle}}(\widetilde D)=\eta_{\langle h\rangle}(\widetilde D)$.  Hence the delocalized higher eta invariant is indeed a natural generalization of Lott's delocalized eta invariant.
\end{remark}

\begin{remark}
	In fact, the delocalized higher eta invariant can be defined for a more general class of representatives besides $u_t(\widetilde D)$ above.  We shall deal with the general case in Section \ref{DelocalizedPairing}.
\end{remark}
\begin{remark}
	The map $\varphi\#\tr$ naturally extends to  the unitization  $(C^\ast(\widetilde M)^G)^+$ of $C^\ast(\widetilde M)^G$ by defining
	 $\widetilde{\varphi\#\tr}$ to be the map that vanishes on the identity element in $(C^\ast(\widetilde M)^G)^+$. With this notation, the formula of $\eta_{\varphi}(\widetilde D)$ becomes 
	\begin{equation*}
	\eta_{\varphi}(\widetilde D)=\frac{m!}{\pi i}\int_0^\infty 
	\widetilde{\varphi\#\tr}(
	\dot{u}_t(\widetilde D) u_t^{-1}(\widetilde D)\otimes (u_t(\widetilde D)\otimes u_t^{-1}(\widetilde D)^{\otimes m})
	dt.
	\end{equation*}
\end{remark}
\begin{remark}
	We have only discussed the odd case so far. The even case is completely analogous. In this case, in the construction of the higher rho invariant $\rho(\widetilde D)$,  the path of invertibles $\{u_t(\widetilde D)\}_{0\leqslant t < \infty}$ is replaced by  a path of projections $\{p_t(\widetilde D)\}_{0\leqslant t< \infty }$ (see for example \cite{Xiepos}).  Given  $\varphi\in C^{2m+1}(\mathbb CG,\langle h\rangle)$, the delocalized higher eta invariant of $\widetilde D$ with respect to $\varphi$ is defined to be
	    \begin{equation*}
	\eta_{\varphi}(\widetilde D)=\frac{1}{\pi i}\frac{(2m)!}{m!}\int_0^\infty \eta_{\varphi}(\widetilde D,t)dt,
	\end{equation*}
	where 
	\begin{equation*}
	\eta_\varphi(\widetilde D,t)=\widetilde{\varphi\#\tr}([\dot p_t(\widetilde D), p_t(\widetilde D)]\otimes p_t(\widetilde D)^{\otimes 2m+1}).
	\end{equation*}
\end{remark}

The integral formula in line $\eqref{etahigher}$ does not converge in general. In the following,  we shall show that the convergence holds  whenever $\widetilde D$ has a sufficiently large spectral gap at zero.
\begin{definition}\label{def:exgrow}
Let $\langle h\rangle$ be a nontrivial conjugacy class of $G$.	An element $\varphi\in C^n(\mathbb CG,\left\langle h\right\rangle )$ is said to have  exponential growth if there exist $C$ and $K_{\varphi}>0$ such that
	\begin{equation}\label{eq:growthofphi}
	|\varphi(g_0,g_1,\cdots,g_n)|\leqslant C e^{K_{\varphi}\cdot(\ell (g_0)+\ell (g_1)+\cdots+\ell (g_n))}
	\end{equation}
	for $\forall (g_0,g_1,\cdots,g_n)\in G^{n+1}$.
\end{definition}

 Similar to the definition of $\uptau_{\langle h\rangle}$ in line \eqref{eq:equivalentwithlength},
we define $\uptau$ to be the following positive number 
\begin{equation}\label{eq:equivalentwithlengthG}
\uptau=\liminf_{\ell(g)\to\infty}\Big(\inf_{x\in\mathcal F}\frac{\dist(x,gx)}{\ell(g)}\Big).
\end{equation}
Since $G$ is finitely generated, there exist $C$ and $K_G>0$ such that
\begin{equation}\label{eq:growthofgroup}
\#\{g\in G:\ell(g)=n\}\leqslant Ce^{K_G\cdot n}.
\end{equation}
\begin{definition}\label{def:specbd}
We define \begin{equation}\label{eq:higherspectralLB}
	\sigma_\varphi\eqqcolon \frac{2(K_G+K_{\varphi})\cdot c_D}{\uptau}.
	\end{equation}
	where  $c_D$ is the propagation speed of $D$ as defined in line \eqref{eq:propagationspeed}.
\end{definition}

\begin{theorem}\label{higherordereta}
	Suppose that $\varphi\in C^n(\mathbb CG,\left\langle h\right\rangle )$ has  exponential growth. With the same notation from above, if the spectral gap of $\widetilde D$ at zero   is larger than $\sigma_\varphi$ given in line \eqref{eq:higherspectralLB} above, then the delocalized higher eta invariant  $\eta_{\varphi}(\widetilde D)$ given in line $\eqref{etahigher}$  converges absolutely.
\end{theorem}
\begin{proof}
    The proof is divided into three steps. First we show that $\eta_{\varphi}(\widetilde D,t)$ is well-defined for any fixed $t>0$ (Proposition \ref{convergenceofTRhigher}). Next we show that its integral for small $t$ converges absolutely (Proposition \ref{convergenceatzerohigher}). The last step is to show the convergence of the integral for large $t$ (Proposition \ref{convergenceatinfhigher}). In fact, only the last step requires the spectral gap of  $\widetilde  D$ to be larger than $\sigma_{\varphi}$.
\end{proof}

We say $G$ (resp. $\varphi$) has \mbox{sub-exponential} growth if we may choose $K_G$ in line \eqref{eq:growthofgroup} (resp. $K_{\varphi}$ in line \eqref{eq:growthofphi}) to be arbitrarily small. The following corollary is an immediate consequence of Theorem $\ref{higherordereta}$ above. 
\begin{corollary}
	With the same notation from above, if both $G$ and $\varphi$ have \mbox{sub-exponential} growth, then the delocalized higher eta invariant  $\eta_{\varphi}(\widetilde D)$ given in line $\eqref{etahigher}$  converges absolutely,  as long as $\widetilde D$ has a spectral gap at zero.
\end{corollary}
	
\begin{proposition}\label{convergenceofTRhigher}
	Suppose that $\langle h\rangle$ be a nontrivial conjugacy class of $G$ and $\varphi\in C^n(\mathbb CG,\left\langle h\right\rangle )$ has exponential growth. If $N\geqslant \frac{3}{2}\dim M+5$ and $\Lambda$ is sufficiently large, then there exists $C>0$ such that for any $f_i\in\mathcal A_{\Lambda,N}$ ($i=0,1,\cdots,n$),
	$$\varphi\#\tr(f_0(\widetilde D)\otimes \cdots \otimes f_n(\widetilde D))\leqslant C\|f_0\|_\mathcal A\cdots \|f_n\|_\mathcal A.$$
\end{proposition}
\begin{proof}
	Fix a symmetric generating set $S$ of $G$. Let $\ell$ be the length function on $G$ determined by $S$. Denote the cardinality of $S$ by $|S|$. 
	
	For each $0\leqslant i \leqslant n$, let  $w_i=f_i(\widetilde D)$ for some function $f_i\in\mathcal A_{\Lambda,N}$. Let us write $w_i$ as a formal sum $w_i=\sum_g w_i^g g$. 
	If we denote by $k_i(x,y)$ the Schwartz kernel of $w_i$, then
	\begin{equation}\label{3.30}
	\tr(w_0^{g_0}\cdots w_n^{g_{n}})
	=\int_{\mathcal F^{n+1}}
	\tr(k_0(x_0,g_0x_1)\cdots k_n(x_{n},g_{n}x_0))dx_0\cdots dx_{n}.   
	\end{equation}
	It follows from Lemma \ref{kernelupperbound} and the definition of $\|\cdot \|_{\mathcal A}$ (cf. Definition $\ref{ALambdaN}$)  that there exists a constant $C_1$ such that
	$$\|k_i(x,y)\|\leqslant C_1\|f_i\|_{\mathcal A} \textup{ for } 0\leqslant i \leqslant n.$$
	For any $(g_0,g_1,\cdots,g_n)\in G^{n+1}$, we divide $\mathcal F^{n+1}$ into $(n+1)$ disjoint (possibly empty) Borel sets $\mathcal F^{n+1}_{(j),(g_0,g_1,\cdots,g_n)}$ such that $\dist(x_j,g_jx_{j+1})$ is the maximum of the set  $\{\dist(x_i,g_ix_{i+1})\}_{0\leq i \leq n}$, where $x_{n+1} = x_0$. In other words, 
	$$\dist(x_j,g_jx_{j+1})\geqslant \dist(x_i,g_ix_{i+1}) \textup{ on } \mathcal F^{n+1}_{(j),(g_0,g_1,\cdots,g_n)} \textup{ for all } 0\leqslant i \leqslant n. $$
	If no confusion is likely to arise, we shall write $\mathcal F^{n+1}_{(j)}$ in place of $\mathcal F^{n+1}_{(j), (g_0,g_1,\cdots,g_n)}$.
	
 Since the action of $G$ on $\widetilde M$ is free and cocompact, there exist $C_1,C_2$ such that
	$$\dist(x,gy)\geqslant C_1\ell(g)-C_2, \textup{ for } \forall x,y\in\mathcal F.$$
It follows that  
	$$\dist(x_j,g_jx_{j+1})\geqslant C_1\frac{\sum_{i=0}^n\ell(g_i)}{n+1}-C_2$$
	on $\mathcal F^{n+1}_{(j)}$. 
	
Again, as the action of $G$ on $\widetilde M$ is free and cocompact, there exists $\varepsilon>0$ such that 
$\dist(x,gx)>\varepsilon$ for all $x\in \widetilde M$ and all $g\ne e$. Note that, since the metric on $\widetilde M$ is $G$-equivariant, we have
$$\sum_{i=0}^n\dist(x_i,g_ix_{i+1})\geqslant \dist(x_0,g_0g_1\cdots g_n x_0).$$
To summarize,  we have 
	$$\dist(x_j,g_j x_{j+1})\geqslant\max\big\{\frac{\varepsilon}{n+1},C_1\frac{\sum_{i=0}^n\ell(g_i)}{n+1}-C_2\big\},$$
on $\mathcal F^{n+1}_{(j)}$.	It follows from  Lemma \ref{functionalcal} that, for any $g_0g_1\cdots g_n\in\langle h\rangle$, there exist $C_3>0$ and $C_4>0$ such that
	\begin{align*}
	&
	\left|\int_{\mathcal F^{n+1}_{(j)}}\tr(k_0(x_0,g_0x_1)\cdots k_n(x_{n},g_{n}x_0))dx_0\cdots dx_{n}\right|\\
	\leqslant& C_3\cdot F_{f_j}\Big(C_4\cdot\max\big\{\frac{\varepsilon}{n+1},C_1\frac{\sum_{i=0}^n\ell(g_i)}{n+1}-C_2\big\}\Big)
	\cdot\prod_{i\ne j}\|f_i\|_{\mathcal A}.
	\end{align*}
	Since $\varphi$ has exponential growth, there exist $C,C',K_\varphi>0$ such that 
	\begin{align*}
	&|\varphi\#\tr(w_0\otimes w_1\otimes\cdots\otimes w_n)|\\
	\leqslant & C\sum_{g_0g_1\cdots g_n\in\langle h\rangle} e^{K_\varphi\sum_{i=0}^n\ell(g_i)}\cdot |\tr(w_0^{g_0}\cdots w_n^{g_{n}})|\\
	\leqslant & C\sum_{g_0g_1\cdots g_n\in\langle h\rangle}e^{K_\varphi\sum_{i=0}^n\ell(g_i)}
	\cdot \Big( \sum_{j=0}^{n}\int_{\mathcal F^{n+1}_{(j)}}|\tr(k_0(x_0,g_0x_1)\cdots k_n(x_{n},g_{n}x_0))|dx_0\cdots dx_{n}\Big)\\
	\leqslant & CC_3\sum_{m=1}^\infty e^{K_\varphi\cdot m}\cdot |S|^{(n+1)m}\left(\sum_{j=0}^n
	\Big[F_{f_j}\Big(C_4\cdot\max\big\{\frac{\varepsilon}{n+1},\frac{mC_1}{n+1}-C_2\big\}\Big)
	\cdot\prod_{i\ne j}\|f_i\|_{\mathcal A}\Big]\right)\\
	\leqslant & C'C_3\sum_{m=1}^\infty e^{K_\varphi\cdot m}\cdot |S|^{(n+1)m}
	\exp\Big(-\frac{\Lambda}{2}C_4\cdot\max\big\{\frac{\varepsilon}{n+1},\frac{mC_1}{n+1}-C_2\big\}\Big)
	\prod_{i=0}^n\|f_i\|_{\mathcal A},
	\end{align*}
	where the last summation converges for sufficiently large $\Lambda$ by Lemma $\ref{estimateA}$. This finishes the proof. 

\end{proof}

Let us now prove the convergence of the integral in line $\eqref{etahigher}$ for small $t$. 
\begin{proposition}\label{convergenceatzerohigher}
	Suppose that $\langle h\rangle$ be a nontrivial conjugacy class of $G$ and $\varphi\in C^n(\mathbb CG,\langle h\rangle)$ has exponential growth. If $N\geqslant \frac{3}{2}\dim M+5$ and $\Lambda$ is sufficiently large, then  the following integral
	$$\int_0^1 t^k\varphi\#\tr(f_0(t\widetilde  D)\otimes\cdots\otimes f_n(t\widetilde  D))dt$$
	is absolutely convergent for all $f_i\in\mathcal A_{\Lambda,N}$ and any $k\in\mathbb R$. 
\end{proposition}
\begin{proof}
	Define $f_{j,t}(x)=f_j(tx)$ for $j=0,1,\cdots,n$.
	From the definition of $\|\cdot\|_\mathcal A$ (cf.  Definition \ref{ALambdaN}), we have
	\[ \|f_{j,t}\|_{\mathcal A}\leqslant t^{-N}\|f_j\|_{\mathcal A} \textup{ for all } t\in (0, 1] \textup{ and } f_j \in \mathcal A_{\Lambda, N}. \]	
	Recall the inequality in line $\eqref{eq:F_fwitht}$:
	$$F_{f_{j,t}}(s)\leqslant \frac{F_{f_j}(s\cdot t^{-1})}{t^{\frac{3}{2}\dim M+3}}\leqslant \frac{C_\Lambda}{t^{\frac{3}{2}\dim M+3}} \cdot  \|f_j\|_{\mathcal A}\cdot \exp\Big(\frac{-\Lambda s}{2t}\Big).$$
	Similar to the proof of Proposition $\ref{convergenceofTRhigher}$, there exist positive constants $K_\varphi, C_1, C_2,$ $C_3,$ $C_4, C_5, C_6$ and  $C_7$ such that 	
	\begin{align*}
	&|t^k\varphi\#\tr(f_0(t\widetilde  D)\otimes\cdots\otimes f_n(t\widetilde  D))|\\
	\leqslant & t^k C_3\sum_{m=1}^\infty e^{K_\varphi\cdot m}\cdot |S|^{(n+1)m}\left(\sum_{j=0}^n
	\Big[F_{f_{j, t}}\Big(C_4\cdot\max\big\{\frac{\varepsilon}{n+1},\frac{mC_1}{n+1}-C_2\big\}\Big)
	\cdot\prod_{i\ne j}\|f_{i, t}\|_{\mathcal A}\Big]\right) \\
	\leqslant & t^{k-(\frac{3}{2}\dim M+3)n-nN} C_5\prod_{i=0}^n\|f_i\|_{\mathcal A}\sum_{m=1}^\infty
	e^{m C_6}  \exp\left(-mC_7\cdot \frac{\Lambda}{2t}\right), 
	\end{align*}
	which proves the proposition as long as $\Lambda$ is sufficiently large. 
\end{proof}

To prove the convergence of the integral in line $\eqref{etahigher}$ for large $t$, we need to fix some notation.
    \begin{definition}\label{def:A_Phi}
    	For any $K>0$, let $\mathcal L_K$ be the subspace of $\mathcal K\otimes C_r^*(G)$ consisting of  operators $A=\sum A_g g$ such that
    	\begin{enumerate}
    		\item for any $g\in G$, $A_g$ is of  trace class; 
    		\item  and we have $$\sum_{g\in G} e^{K\cdot \ell(g)} |A_g|_1<\infty,$$
    		where $|\cdot|_1$ is the trace norm.
    	\end{enumerate}
    	Equip $\mathcal L_K$ with the following norm
    	$$\|A\|_{\mathcal L}=\sum_{g\in G}e^{K\cdot \ell(g)}|A_g|_1.$$
    \end{definition}
    \begin{lemma}\label{lemma:A_Phi}
    	The space $\mathcal L_K$ is a Banach algebra with the norm $\|\cdot\|_{\mathcal L}$.
    \end{lemma}
    \begin{proof}
    	It is not difficult to see that $\mathcal L_K$ is a Banach space under the norm $\|\cdot\|_{\mathcal L}$. It remains to show that $\|\cdot\|_{\mathcal L}$ is sub-multiplicative. Indeed, given $A_1 = \sum_{g\in G}A_{1,g} g$ and $A_2 = \sum_{g\in G}A_{2,g} g$ in $\mathcal L_K$, we have 
    	\begin{align*}
    	\|A_1A_2\|_\mathcal L\leqslant& \sum_{g\in G}e^{K\cdot \ell(g)}\sum_{g_1\in G}|A_{1,gg_1}|_1\cdot |A_{2,g_1^{-1}}|_1\\
    = & \sum_{g_1\in G}\sum_{g_2\in G}e^{K\cdot  (\ell(g_1g_2)-\ell(g_1)-\ell(g_2))} \left(e^{K\cdot \ell(g_1)} |A_{1,g_1}|_1 \right) 
   \left(e^{K\cdot \ell(g_2)} |A_{2,g_2}|_1\right)\\
    	\leqslant & \|A_1\|_{\mathcal L}\cdot \|A_2\|_{\mathcal L}.
    	\end{align*}
    \end{proof}

Now we prove the convergence of the integral in $\eqref{etahigher}$ for large $t$.
    \begin{proposition}\label{convergenceatinfhigher}
       Suppose that\footnote{Here $\varphi$ is not necessarily delocalized.} $\varphi\in C^{2m}(\mathbb CG)$ has exponential growth. If the spectral gap of $\widetilde D$  at zero is larger than $\sigma_\varphi$ given in line \eqref{eq:higherspectralLB}, then the integral
        $$\int_1^{\infty} \eta_\varphi(\widetilde  D,t)dt$$
        is absolutely convergent.
    \end{proposition}
    \begin{proof}
   The case where $\widetilde  D$ has a spectral gap at zero can be deduced from the invertible case by replacing $\widetilde  D$ with its restriction to the orthogonal complement of the kernel of $\widetilde  D$. Without loss of generality, let us assume $\widetilde D$ is invertible. 
   
   	Since $\varphi$ has exponential growth, there exist $C$ and $K_{\varphi}$ such that
$$
|\varphi(g_0,g_1,\cdots,g_n)|\leqslant C e^{K_{\varphi}(\ell (g_0)+\ell (g_1)+\cdots+\ell (g_n))}
$$
for $\forall (g_0,g_1,\cdots,g_n)\in G^{n+1}$. 
We first show that $\dot{u}_t(\widetilde D)u^{-1}_t(\widetilde D)$, $u_t(\widetilde D)-1$ and $u_t^{-1}(\widetilde D)-1$ all lie in $\mathcal L_{K_{\varphi}}$, where  $\mathcal L_{K_{\varphi}}$ is given in Definition \ref{def:A_Phi} above.

Let $\sigma$ be the spectral gap of  $\widetilde D$ at zero.  Note that 
$$\dot{u}_t(\widetilde D)u^{-1}_t(\widetilde D)=2\sqrt \pi i\widetilde D e^{-t^2\widetilde D^2}.$$
Since $\sigma>\sigma_\varphi = \frac{2(K_G+K_{\varphi})\cdot c_D}{\uptau}$, the same proof of Proposition \ref{convergencofetanearinf} shows that $\dot{u}_t(\widetilde D)u^{-1}_t(\widetilde D)$ lies in $\mathcal L_{K_{\varphi}}$ and furthermore there exists sufficiently small $\omega>0$ such that
\begin{equation}\label{Lbound}
\|\dot{u}_t(\widetilde D)u^{-1}_t(\widetilde D)\|_\mathcal L\leqslant 2\sqrt{\pi} e^{-\omega t^2}.
\end{equation}
Observe that 
$$u_t(\widetilde D)-1=
\exp\left(-\int_t^{+\infty}\dot{u}_s(\widetilde D)u^{-1}_s(\widetilde D)ds\right)-1.$$
 By the inequality in line $\eqref{Lbound}$ above, we have
$$\int_t^{+\infty}\|\dot{u_s}(\widetilde D)u^{-1}_s(\widetilde D)\|_\mathcal L ds\leqslant \frac{2\pi e^{-\omega t^2}}{\sqrt{\omega}}.$$
Thus for any $t>1$, we have
$$\int_t^{+\infty}\dot{u}_s(\widetilde D)u^{-1}_s(\widetilde D)ds\in\mathcal L_{K_\varphi}$$
and 
$$\left\|\int_t^{+\infty}\dot{u}_s(\widetilde D)u^{-1}_s(\widetilde D)ds\right\|_\mathcal L\leqslant \frac{2\pi e^{-\omega t^2}}{\sqrt{\omega}}.$$
By Lemma $\ref{lemma:A_Phi}$, $\mathcal L_{K_{\varphi}}$ is a Banach algebra under the norm $\|\cdot \|_{\mathcal L}$. It follows that 
$$u_t(\widetilde D)-1=
\exp\left(-\int_t^{+\infty}\dot{u}_s(\widetilde D)u^{-1}_s(\widetilde D)ds\right)-1\in\mathcal L_{K_\varphi}.$$
Moreover, we have 
$$\|u_t(\widetilde D)-1\|_\mathcal L\leqslant 
\sum_{n=1}^\infty\frac{1}{n!}\left(\left\|\int_t^{+\infty}\dot{u_s}(\widetilde D)u^{-1}_s(\widetilde D)ds\right\|_\mathcal L\right)^n\leqslant e^{C_1\cdot e^{-\omega t^2}}-1, $$
where $C_1 = \frac{2\pi}{\sqrt{\omega}}$.  A similar estimate holds for $u_t^{-1}(\widetilde D)-1$ as well. A straightforward calculation shows that 
$$|\eta_\varphi(\widetilde D,t)|\leqslant C\cdot \|\dot{u}_t(\widetilde D)u^{-1}_t(\widetilde D)\|_{\mathcal L} \cdot 
\|u_t(\widetilde D)-1\|_{\mathcal L}^m \cdot \|u_t^{-1}(\widetilde D)-1\|_\mathcal L^m.$$
It follows that 
$$|\eta_\varphi(\widetilde D,t)|\leqslant 2\sqrt{\pi} e^{-\omega t^2}(e^{C_1\cdot e^{-\omega t^2}}-1)^{2m},$$ where the right hand side is clearly integrable over $[1, \infty)$. 
This finishes the proof.
    \end{proof}

The following proposition shows that the delocalized higher eta invariant in Definition $\ref{higherdelocalizedeta}$ is independent of the choice of representative within the same cyclic cohomology class.  
\begin{proposition}\label{welldefinedofeta}
    Let $\langle h\rangle$ be a nontrivial conjugacy class of $G$. Suppose that $\varphi_1$ and $\varphi_2$ are two cyclic cocycles in $C^{2m}(\mathbb CG,\langle h\rangle)$ with exponential growth and $\varphi_1$ and $\varphi_2$ are cohomologous via a cochain with exponential growth. If $\widetilde D$ has a sufficiently large spectral gap at zero, then $\eta_{\varphi_1}(\widetilde  D)=\eta_{\varphi_2}(\widetilde  D)$.
\end{proposition}

\begin{proof}
    It suffices to show that if $\varphi\in C^{2m-1}(\mathbb CG,\langle h\rangle)$ has exponential growth, then $\eta_{b\varphi}(\widetilde  D)=0$.  By the definition of $\widetilde{\varphi\#\tr}$, if $w_j =1$ for some $j$, then 
	\begin{equation*}
	\widetilde{\varphi\#\tr}(w_0\otimes w_1\otimes w_2\otimes\cdots\otimes w_n)=0,
	\end{equation*}
	For notational simplicity, let us write $u_t$ in place of $u_t(\widetilde D)$ in the following. 
	
Observe that  
	\begin{align*}
	\frac{d}{dt}\widetilde{\varphi\#\tr}((u_t\otimes u_t^{-1})^{\otimes m}) =  &m\, \widetilde{\varphi\#\tr}(\dot{u}_t\otimes u_t^{-1}\otimes (u_t\otimes u_t^{-1})^{\otimes m-1}) \\
	& - m\, \widetilde{\varphi\#\tr}(u_t \otimes u_t^{-1}\dot{u}_t u_t^{-1}\otimes(u_t\otimes u_t^{-1})^{\otimes m-1}),
	\end{align*}
	and 
	\begin{align*}
	(\widetilde{b\varphi\#\tr})(\dot u_t u_t^{-1}\otimes(u_t\otimes u_t^{-1})^{\otimes m}) =&\widetilde{\varphi\#\tr}(\dot{u}_t\otimes u_t^{-1}\otimes (u_t\otimes u_t^{-1})^{\otimes m-1})\\
	&-\widetilde{\varphi\#\tr}(u_t \otimes u_t^{-1}\dot{u}_t u_t^{-1}\otimes(u_t\otimes u_t^{-1})^{\otimes m-1}).
	\end{align*}
It follows that 
	\begin{equation}\label{1sttransgression}
	\frac{d}{dt}\widetilde{\varphi\#\tr}((u_t\otimes u_t^{-1})^{\otimes m})
	=m\cdot \eta_{b\varphi}(\widetilde D,t).
	\end{equation}
In particular, we have 
\begin{align*}
\int_{t_0}^{t_1}\eta_{b\varphi}(\widetilde D,t)dt
=&\widetilde{\varphi\#\tr}((u_t\otimes u_t^{-1})^{\otimes m})\big| _{t=t_1}-
\widetilde{\varphi\#\tr}((u_t\otimes u_t^{-1})^{\otimes m})\big| _{t=t_0}\\
=&\varphi\#\tr(((u_t-1)\otimes (u_t^{-1}-1))^{\otimes m})\big| _{t=t_1}\\
&-\varphi\#\tr(((u_t-1)\otimes (u_t^{-1}-1))^{\otimes m})\big| _{t=t_0}.
\end{align*}
By the proof of Proposition \ref{convergenceatzerohigher}, we have 
$$\lim_{t\to 0}\varphi\#\tr(((u_t-1)\otimes (u_t^{-1}-1))^{\otimes m}) = 0.$$ Furthermore, it follows from the proof of Proposition \ref{convergenceatinfhigher} that 
\[\lim_{t\to \infty}\varphi\#\tr(((u_t-1)\otimes (u_t^{-1}-1))^{\otimes m}) = 0, \] 
 as long as the spectral gap of $\widetilde D$ is sufficiently large.  This finishes the proof.

\end{proof}

In the remaining part of this section, we show that the delocalized higher eta invariant is stable under Connes' periodicity map. First, let us recall the definition of Connes' periodicity map (cf. \cite[Page 121]{Connes}):
\[ S\colon HC^{2m}(\mathbb CG)\to HC^{2m+2}(\mathbb CG)\] 
First, let us fix some notation. For $0\leqslant i\leqslant n$, we define $b_i\colon L^n(\mathbb CG)\to L^{n+1}(\mathbb CG)$ by 
	\begin{align*}
	b_if(g_0,g_1,\cdots,g_{n+1}) &=      f(g_0,g_1,\cdots,g_jg_{j+1},\cdots,g_{n+1}),
	\end{align*}
	and for $i=n+1$, we define $b_{n+1}\colon L^n(\mathbb CG)\to L^{n+1}(\mathbb CG)$ by 
	\[b_{n+1}f(g_0, g_1, \cdots, g_{n+1}) = (-1)^{n+1}f(g_{n+1}g_0,g_1,\cdots,g_n); \]

\begin{definition}\label{def:periodic}
Connes' periodicity map $S\colon  C^n(\mathbb CG)\to C^{n+2}(\mathbb CG)$ is defined to be 
 \[ S\coloneqq  \frac{-1}{(n+2)(n +1)}\sum_{0\leqslant i < j \leqslant n+2}  b_j b_i.\]
\end{definition}
We point out that the above explicit formula of $S$ differs by a constant from Connes' formula in \cite{Connes} (compare with the formula for $S$ in the proof of \cite[Part II, Corollary]{Connes}). This is because the constants appearing in the Connes-Chern character of the current paper are different from those in Connes' original paper  \cite{Connes}.

The following proposition shows that the delocalized higher eta invariant is stable under Connes' periodicity map.
\begin{proposition}\label{suspension}
 Let $\langle h\rangle$ be a nontrivial conjugacy class of $G$. If a cyclic cocycle $\varphi\in C^{2m}(\mathbb CG,\langle h\rangle)$ has exponential growth and the spectral gap of $\widetilde D$ at zero is sufficiently large, then 
	$\eta_{S\varphi}(\widetilde D)=\eta_{\varphi}(\widetilde D)$.
\end{proposition}
\begin{proof}
	We will prove the proposition by a direct computation. It is much easier to follow Connes and carry out the computation in  the context of universal graded differential algebras \cite[Part II]{Connes}. Let us recall the construction of the universal graded differential algebra $\Omega(\mathscr A)$ associated to an algebra $\mathscr A$. Denote by $\ntilde{\mathscr A}$  the algebra obtained from $\mathscr A$ by adjoining a unit: $\ntilde{\mathscr A} = \{ a+ \lambda I \mid a\in  \mathscr A, \lambda \in \mathbb C\}$. Let $\Omega^0(\mathscr A) = \mathscr A$ and 
	\[ \Omega^n(\mathscr A) = \ntilde{\mathscr A}\otimes (\mathscr A)^{\otimes n} \]
	for $n \geq 1$.  The differential $\dbar\colon \Omega^n(\mathscr A) \to \Omega^{n+1}(\mathscr A)$
	is given by 
	\[ \dbar\big((a_0+\lambda I)\otimes a_1\otimes \cdots \otimes a_n) = I\otimes a_0\otimes a_1\otimes \cdots \otimes a_n.  \] 
Clearly, one has $d^2 =0$. 
	The product structure (cf. \cite[Part II, Proposition 1]{Connes})  on $\Omega^\ast(\mathscr A)$ is defined so that the following are satisfied:	\begin{enumerate}
		\item $\dbar(\omega_1 \omega_2) = (\dbar\omega_1) \omega_2 + (-1)^{|\omega_1|} \omega_1 \dbar\omega_2  $ for $\omega_1 \in \Omega^i$ and $\omega_2 \in \Omega^j$, where $|\omega_1| = i$ is the degree of $\omega_1$;
		\item $\tilde a_0 \dbar a_1 \dbar a_2 \cdots \dbar a_n = \tilde a_0 \otimes a_1\otimes a_2 \otimes \cdots \otimes a_n$ in $\Omega^n(\mathscr A)$. 
	\end{enumerate}
An $(n+1)$-linear functional $\varphi$ on $\mathscr A$ induces a linear functional $\widehat \varphi$ on $\Omega^n(\mathscr A)$ by setting
\[ \widehat \varphi\big((a_0 + \lambda I)da_1\cdots da_n\big) = \varphi(a_0, a_1, \cdots, a_n).  \]
By \cite[Part II, Proposition 1]{Connes}, since $\varphi$ is a cyclic cocycle, $\widehat \varphi$ is a closed graded trace on  $\Omega^n(\mathscr A)$. In particular, we have 
\[ \widehat \varphi( \omega_1 \omega_2) = (-1)^{|\omega_1|\cdot|\omega_2|} \widehat\varphi(\omega_2 \omega_1).  \]
By using the equality $\dbar (\omega_1 \omega_2) = (\dbar \omega_1) \omega_2 + (-1)^{|\omega_1|} \omega_1  \dbar \omega_2$ above, a straightforward calculation shows that the formula for the periodicity operator $S$  becomes the following (compare with \cite[Part II, Corollary 10]{Connes}):
\begin{equation}\label{eq:periodic2}
\widehat{S\varphi}(a_0  \dbar a_1\cdots  \dbar a_{n+2}) =  c_n\sum_{j=1}^{n+1}\widehat \varphi(a_0 \dbar a_1\cdots  \dbar a_{j-1} (a_{j}a_{j+1}) \dbar a_{j+2}\cdots \dbar  a_{n+2}),
\end{equation}
where $\displaystyle c_n = \frac{1}{(n+2)(n +1)}.$

Now we shall prove the proposition. Let  $\varphi\in C^{2m}(\mathbb CG,\langle h\rangle)$ be a cyclic cocycle with exponential growth,  then $S\varphi$ also has exponential growth. Recall that  the delocalized higher eta invariant of $\widetilde D$ with respect to $\varphi$  to be\footnote{Here we us the notation $u_t^{-1}-I$ for the corresponding term $u_t^{-1}-1$ in Definition $\ref{higherdelocalizedeta}$. This is to emphasis that $I$ is the identity operator, which is the unit adjoined.}
	\begin{equation*}
	\eta_{\varphi}(\widetilde D)=\frac{m!}{\pi i}\int_0^\infty \varphi\#\tr(
	\dot{u}_t u_t^{-1}\otimes ((u_t-I)\otimes (u_t^{-1}-I))^{\otimes m}) dt,
	\end{equation*}
	where $u_t= u_t(\widetilde D)$ as in line $\eqref{eq:smallu}$.  
For notational simplicity, let us write 	
\[  \varphi(
\dot{u}_t u_t^{-1}\otimes ((u_t-I)\otimes (u_t^{-1}-I))^{\otimes m}) \]
in place of $\varphi\#\tr(
\dot{u}_t u_t^{-1}\otimes ((u_t-I)\otimes (u_t^{-1}-I))^{\otimes m})$, and furthermore write  
\[ a_0 = \dot{u}_tu_t^{-1} \textup{ and for $j\geq 1$,  }    a_j = \begin{cases}
u_t -I & \textup{ if $j$ is odd,} \\[5pt]
u_t^{-1} - I & \textup{ if $j$ is even.}
\end{cases} \]
By the above discussion, we have 
\begin{align*}
\widehat \varphi(
a_0da_1\cdots da_{2m})  = \varphi(
\dot{u}_t u_t^{-1}\otimes ((u_t-1)\otimes (u_t^{-1}-1))^{\otimes m}). 
\end{align*}
The following observations will be useful in the computation below. 
\begin{enumerate}
		\item $(u_t - I) (u_t^{-1} - I) = 2I - u_t - u_t^{-1}$;
	\item $\dbar (u_t - I) = \dbar u_t$ and $\dbar (u^{-1}_t - I) = \dbar u_t^{-1}$ ;
	\item $ (\dbar u_t) u_t^{-1} = -u_t (\dbar u_t^{-1})$ and $ (\dbar u^{-1}_t) u_t = -u_t^{-1} (\dbar u_t).$
\end{enumerate}
Observation $(1)$ immediately implies that 
\begin{align*}
& \widehat \varphi(a_0\dbar a_1\cdots \dbar a_{j-1} (a_{j}a_{j+1})\dbar a_{j+2}\cdots \dbar a_{2m+2})  \\
= & 2 \widehat \varphi(a_0\dbar a_1\cdots \dbar a_{j-1} \dbar a_{j+2}\cdots \dbar a_{2m+2})   \\
& -\widehat \varphi(a_0\dbar a_1\cdots \dbar a_{j-1} (u_t) \dbar a_{j+2}\cdots \dbar a_{2m+2})  -\widehat \varphi(a_0\dbar a_1\cdots \dbar a_{j-1} (u_t^{-1}) \dbar a_{j+2}\cdots \dbar a_{2m+2}).
\end{align*}
Now observations $(2)$ and $(3)$ imply the following: 
\begin{enumerate}[label=(\roman*)]
	\item if $j$ is odd, then 
	\begin{align*}
	& \widehat \varphi(a_0\dbar a_1\cdots \dbar a_{j-1} (u_t) \dbar a_{j+2}\cdots \dbar a_{2m+2}) \\
	=&\widehat\varphi(\dot{u}_t \underbrace{\dbar u_t^{-1}\dbar u_t\cdots \dbar u_t^{-1} \dbar u_t}_{\textup{ $(j-1)$ terms}}\underbrace{\dbar u_t\dbar u_t^{-1}\cdots \dbar u_t \dbar u^{-1}_t}_{\textup{ $(2m-j+1)$ terms}}),
	\end{align*}
and 
\begin{align*}
	& \widehat \varphi(a_0\dbar a_1\cdots \dbar a_{j-1} (u_t^{-1}) \dbar a_{j+2}\cdots \dbar a_{2m+2}) \\
	=&- \widehat\varphi(\dot{u}^{-1}_t \underbrace{\dbar u_t\dbar u_t^{-1}\cdots \dbar u_t \dbar u_t^{-1}}_{\textup{ $(j-1)$ terms}}\underbrace{\dbar u_t^{-1}\dbar u_t\cdots \dbar u^{-1}_t \dbar u_t}_{\textup{ $(2m-j+1)$ terms}}); 
	\end{align*}
	\item 
	if $j$ is even, then 
	\begin{align*}
	& \widehat \varphi(a_0\dbar a_1\cdots \dbar a_{j-1} (u_t) \dbar a_{j+2}\cdots \dbar a_{2m+2}) \\
	=&\widehat\varphi(\dot{u}^{-1}_t \underbrace{\dbar u_t\dbar u_t^{-1}\cdots \dbar u^{-1}_t \dbar u_t}_{\textup{ $(j-1)$ terms}}\underbrace{\dbar u_t\dbar u_t^{-1}\cdots \dbar u_t \dbar u^{-1}_t \dbar u_t}_{\textup{ $(2m-j+1)$ terms}}), 
	\end{align*}
	and 	\begin{align*}
	& \widehat \varphi(a_0\dbar a_1\cdots \dbar a_{j-1} (u_t^{-1}) \dbar a_{j+2}\cdots \dbar a_{2m+2}) \\
	=&- \widehat\varphi(\dot{u}_t \underbrace{\dbar u_t^{-1}\dbar u_t\cdots \dbar u_t \dbar u_t^{-1}}_{\textup{ $(j-1)$ terms}}\underbrace{\dbar u_t^{-1}\dbar u_t\cdots \dbar u_t^{-1}\dbar u_t \dbar u_t^{-1}}_{\textup{ $(2m-j+1)$ terms}}). 
	\end{align*}
\end{enumerate}
Since $\widehat \varphi$ is a closed graded trace, it follows that 
\begin{align*}
& \sum_{j=1}^{2m+1} \widehat \varphi(a_0\dbar a_1\cdots \dbar a_{j-1} (u_t) \dbar a_{j+2}\cdots \dbar a_{2m+2})  =  \frac{d}{dt}  \widehat \varphi\big((u_t-1) \underbrace{\dbar u_t\dbar u^{-1}_t\cdots \dbar u_t \dbar u_t^{-1}}_{\textup{ $2m$ terms}}\big) 
\end{align*} 
and 
\begin{align*}
& \sum_{j=1}^{2m+1} \widehat \varphi(a_0\dbar a_1\cdots \dbar a_{j-1} (u^{-1}_t) \dbar a_{j+2}\cdots \dbar a_{2m+2})  =  \frac{d}{dt}  \widehat \varphi\big((u^{-1}_t-1) \underbrace{\dbar u_t^{-1}\dbar u_t\cdots \dbar u_t^{-1} \dbar u_t}_{\textup{ $2m$ terms}}\big) 
\end{align*} 
On the other hand, by the proof of Proposition \ref{convergenceatzerohigher} and Proposition \ref{convergenceatinfhigher}, we have 
\begin{align*}
& \int_0^\infty \frac{d}{dt}  \widehat \varphi\big((u_t-1) \underbrace{\dbar u_t\dbar u^{-1}_t\cdots \dbar u_t \dbar u_t^{-1}}_{\textup{ $2m$ terms}}\big) dt \\
= & \lim_{t\to \infty} \varphi\big( (u_t-1) \otimes ((u_t-1)\otimes (u_t^{-1}-1))^{\otimes m}\big)  \\
& - \lim_{t\to 0} \varphi\big( (u_t-1) \otimes ((u_t-1)\otimes (u_t^{-1}-1))^{\otimes m}\big)  \\
=&  0
\end{align*}
Similarly, we have
\[ \int_0^\infty \frac{d}{dt}  \widehat \varphi\big((u^{-1}_t-1) \underbrace{\dbar u_t^{-1}\dbar u_t\cdots \dbar u_t^{-1} \dbar u_t}_{\textup{ $2m$ terms}}\big)dt = 0. \]   
To summarize, we have 
\[ \int_0^\infty\widehat{S\varphi}(\dot{u}_t u_t^{-1} \underbrace{\dbar u_t\dbar u^{-1}_t\cdots \dbar u_t \dbar u_t^{-1}}_{\textup{ $(2m+2)$ terms}}) dt = \frac{1}{m+1} \int_0^\infty\widehat{\varphi}(\dot{u}_t u_t^{-1}  \underbrace{\dbar u_t\dbar u^{-1}_t\cdots \dbar u_t \dbar u_t^{-1}}_{\textup{ $2m$ terms}}) dt, \]
which implies that 
\[ 	\eta_{S\varphi}(\widetilde D)  = 	\eta_{\varphi}(\widetilde D).\]
This finishes the proof. 
\end{proof}
\subsection{A higher Atiyah-Patodi-Singer index formula}
In this subsection,  for each cyclic cocycle of at most exponential growth, we prove a formal higher Atiyah-Patodi-Singer index theorem (abbr. higher APS index theorem) on manifolds with boundary, under the condition that the operator on the boundary has a sufficiently large spectral gap at zero. We point out that there is no condition on the fundamental group in this formal higher APS index theorem.

Leichtnam and Piazza proved a higher APS index theorem in terms of noncommutative differential forms on a certain smooth dense subalgebra of the reduced $C^\ast$-algebra of the fundamental group \cite[Theorem 4.1]{MR1708639}. Heuristically speaking, our version of higher APS theorem is the pairing between their version of higher APS theorem and cyclic cocycles of the fundamental group. However, this is not the approach we take in this section. In fact, in general it is rather difficult to make this heuristic argument rigorous. A main difficulty here is whether cyclic cocycles of a group algebra extends continuously to cyclic cocycles on a given smooth dense subalgebra of the reduced group $C^\ast$-algebra. In this section, our approach is based on the convergence results from the  previous sections, and avoids the subtle issue of continuous extension of cyclic cocylces. On the other hand, in order to apply the higher APS index theorem
 in this subsection to problems in geometry and topology (cf. \cite{Piazzarho, MR3590536, Weinberger:2016dq}), one actually needs to extend the pairing to  be defined at the level of (periodic) cyclic cohomology and $K$-theory of $C^\ast$-algebras. In later sections, we shall use Puschnigg's smooth dense subalgebra to define such a pairing at the level of (periodic) cyclic cohomology and $K$-theory of $C^\ast$-algebras for all hyperbolic groups. As a consequence, in the case of hyperbolic groups, we shall prove a higher APS index theorem without the assumption that the operator on the boundary has a sufficiently large spectral gap at zero (cf. Section $\ref{Application}$).

 For notational simplicity, we shall only discuss the case of even dimensional spin manifolds. The exact same strategy clearly works for the more general case of Dirac-type operators acting on Clifford modules over Riemannian manifolds of all dimensions. 
 
 Let $W$ be an even dimensional compact spin manifold with boundary $M$, and $D$ the Dirac operator on $W$. Suppose the metric of $W$ is a product metric when restricted to the boundary $M$. Let $G$ be a finitely presented discrete group and $\widetilde W$ a regular $G$-covering space of $W$.  Let $\widetilde D$ be the lift of $D$ to $\widetilde W$  and $\widetilde D_\partial$ the restriction $\widetilde D$ to the boundary of $\widetilde W$.

Let us briefly review Lott's noncommutative differential higher eta invariant. We shall follow  closely the notation in Lott's paper \cite{Lott1}. 
For each $q>0$, we define  $\mathfrak B_q^{\omega}$ to be the following dense subalgebra of $C_r^*(G)$: 
\[ \mathfrak B_q^{\omega}=\big\{f: G\to \mathbb C  \mid \sum_{g\in G}e^{q\cdot \ell(g)}|f(g)|<\infty\big\}, \]
where $\ell$ is a word-length function on $G$. Note that $\mathfrak B_q^{\omega}$ is generally \emph{not} closed under holomorphic functional calculus in $C_r^*(G)$.   The universal graded differential algebra of $\mathfrak B^{\omega}_q$ is 
\[ \Omega_\ast(\mathfrak B_q^{\omega}) = \bigoplus_{k=0}^\infty \Omega_k(\mathfrak B_q^{\omega}) \]
where as a vector space, $\Omega_k(\mathfrak B_q^{\omega}) = \mathfrak B_q^{\omega} \otimes (\mathfrak B_q^{\omega}/\mathbb C)^{\otimes k}$. As $\mathfrak B_q^{\omega}$ is a Banach algebra (cf. Lemma $\ref{lemma:A_Phi}$ above), we consider the Banach completion of $\Omega_\ast(\mathfrak B_q^{\omega})$, which will still be denoted by $\Omega_\ast(\mathfrak B_q^{\omega})$.
 

Let $S$ be the restriction of spinor bundle of $W$ on $M$. We denote the corresponding $\mathfrak B_q^{\omega}$-vector bundle by $\mathfrak S = (\widetilde M\times_G \mathfrak B_q^{\omega})\otimes S$ and  the space of smooth sections by $C^\infty (M; \mathfrak S)$. Now suppose $\psi$ is a smooth function on $\widetilde M$ with comapct support  such that 
\[ \sum_{g\in G}g^*\psi=1. \] Then we have a superconnection  $\nabla\colon C^\infty(M; \mathfrak S)\to C^\infty(M; \mathfrak S\otimes_{\mathfrak B_q^{\omega}} \Omega_1(\mathfrak B_q^{\omega}))$ given by 
$$\nabla(f)=\sum_{g\in G}(\psi \cdot  g^*f)\otimes_{\mathfrak B_q^{\omega}} dg.$$  See \cite{MR790678} for more details of the superconnection formalism. 
\begin{definition}[{\cite[Section 4.4 \& 4.6]{Lott1}}]
	Lott's higher eta invariant $\widetilde \eta(\widetilde D_\partial)$ is defined by the formula
	\[ \widetilde \eta(\widetilde D_\partial)=\int_0^\infty \STR(\widetilde De^{-(t\widetilde D_\partial +\nabla)^2})dt, \]
	where $\STR$ is the corresponding supertrace, cf. \cite[Proposition 22]{Lott1}. 
\end{definition}

\begin{remark}
	The above integral formula for $\widetilde \eta(\widetilde D_\partial)$ generally does not converge to define an element in $\Omega_*(\mathfrak B_q^\omega)$. On the other hand, the estimates in the previous subsections show that the above integral converges absolutely and defines an element in $\Omega_*(\mathfrak B_q^\omega)$,  provided the spectral gap of $\widetilde D_\partial$ at zero is sufficiently large.  
\end{remark}

Now let $W_\infty$ be the complete Riemannian manifold obtained by attaching an infinite cylinder $M\times [0, \infty)$ to $W$. We still denote the associated Dirac operator on $W_\infty$ by $D$ and its lift to $\widetilde W_\infty$ by $\widetilde D$. Combining Quillen's superconnection formalism with Melrose's $b$-calculus formalism \cite{RM93}, for each $t>0$, one defines the $b$-Connes-Chern character of $\widetilde D$ on $\widetilde W_\infty$  to be\footnote{It is not difficult to adapt the estimates from the previous subsections to the $b$-calculus setting and show that $b\text{-}\tr_s(e^{-(t\widetilde D + \nabla )^2})$ indeed defines an element in $\Omega_*(\mathfrak B_q^\omega)$.}
\[ b\textup{-}\Ch_t(\widetilde D) =  b\text{-}\STR(e^{-(t\widetilde D + \nabla )^2}) \in \Omega_*(\mathfrak B_q^\omega),  \] 
where $b\text{-}\STR$ is the corresponding $b$-supertrace in this $b$-calculus setting. See for example \cite{Leichtnam} for more details.

\begin{theorem}\label{b-aps}
	Assume that $\varphi\in C^{2m}(\mathbb CG)$ is a cyclic cocycle with exponential growth. Let $\sigma_\varphi$ be the positive number from Definition $\ref{def:specbd}$. If the spectral gap of $\widetilde D_\partial$ at zero is larger than  $\sigma_\varphi$, then the pairing $\big\langle\varphi,\widetilde\eta(\widetilde D_\partial)\big\rangle$ converges,  the limit $\lim_{t\to\infty}\big\langle \varphi, b\textup{-}\Ch_t(\widetilde D))\big\rangle$ exists, and furthermore
	\begin{equation}\label{eq:b-APS}
	\lim_{t\to\infty}\big\langle \varphi, b\textup{-}\Ch_t(\widetilde D))\big\rangle =\big\langle\varphi,\int_W \hat A\wedge\omega\big\rangle -\frac{1}{2}\big\langle\varphi,\widetilde\eta(\widetilde D_\partial)\big\rangle.
	\end{equation}
	where $\hat A$ is the associated $\hat A$-form on  $W$ and $\omega$ is an element in $\Omega^\ast(W)\otimes \Omega_\ast(\mathfrak B_q^\omega)$ for some $q>0$, \textup{(}cf. \cite[Theorem 13.6]{Leichtnam}\textup{)}.
	In particular, if both $G$ and $\varphi$ have sub-exponential growth, then the equality in line $\eqref{eq:b-APS}$ holds as long as $\widetilde D_\partial$ is invertible. 
\end{theorem}
\begin{proof}
By Proposition \ref{convergenceofTRhigher} and Proposition \ref{convergenceatzerohigher}, we have that, 
for any $\varepsilon>0$, $\mu>1$ and $k\in\mathbb N$,  there exist $C,N_1,N_2$ such that
\begin{equation}\label{3.27}
\|(\widetilde D^k e^{-t^2\widetilde D^2})(x,y)\|\leqslant C\frac{d(x,y)^{N_1}}{t^{N_2}}\exp\Big(-\frac{d(x,y)^2}{4\mu^2 c_{D}^2t^2}\Big),
\end{equation}
for all  $x,y\in\widetilde W_\infty$ with $d(x,y)>\varepsilon$,
and 
\begin{equation}
\|(\widetilde D_\partial^k e^{-t^2\widetilde D_\partial^2})(x,y)\|\leqslant C\frac{d(x,y)^{N_1}}{t^{N_2}}\exp\Big(-\frac{d(x,y)^2}{4\mu^2 c_{D_\partial}^2t^2}\Big).
\end{equation}
for  all  $x,y\in\widetilde M$ with $d(x,y)>\varepsilon$,
where $c_{D}$ (resp. $c_{D_\partial}$) is the propagation speed of $D$ (resp. $ D_\partial$), cf. line \eqref{eq:propagationspeed}.

Since the spectral gap of $\widetilde D_\partial$ is larger than $\displaystyle \sigma_\varphi = \frac{2(K_G+K_{\varphi})\cdot c_{D_\partial}}{\uptau}$, the proof of Proposition $\ref{convergencofetanearinf}$ shows that  for each $k\in \mathbb N$, the operator  $  \widetilde D_\partial^k e^{-t^2\widetilde D_\partial^2}$ lies in $\mathcal L_{K_{\varphi}}$ (cf. Definition $\ref{def:A_Phi}$) and furthermore there exists sufficiently small $\omega>0$ such that
\begin{equation*}
\|\widetilde D_\partial^k e^{-t^2\widetilde D_\partial^2}\|_\mathcal L\leqslant 2\sqrt{\pi} e^{-\omega t^2}.
\end{equation*} 

Now apply the commutator formula for $b$-trace (cf. \cite[(In.22) on Page 8]{RM93}), and a straightforward calculation shows that  for any $0<t_0<t_1$, the equality 
	\begin{align}
	&b\text{-}\Ch_{t_1}(\widetilde D)-b\text{-}\Ch_{t_0}(\widetilde D)  \notag\\
	= &  - \frac{1}{2}\int_{t_0}^{t_1}\STR(\widetilde D_\partial e^{-(s\widetilde D_\partial +\nabla)^2})ds + d\int_{t_0}^{t_1} b\text{-}\STR(\widetilde D e^{-(\nabla +s\widetilde D)^2})ds \label{eq:bvar}
	\end{align}
holds   in $\Omega_\ast(\mathfrak B_q^\omega)$ with $q = K_\varphi$,  where $d\colon \Omega_\ast(\mathfrak B_q^\omega) \to \Omega_{\ast+1}(\mathfrak B_q^\omega)$ is the differential on $\Omega_\ast(\mathfrak B_q^\omega)$,  cf. \cite[Section 6]{EG93b}\cite[Proposition 14.2]{Leichtnam}. In particular, by pairing both sides of  $\eqref{eq:bvar}$ with $\varphi$, we have 
\[ \big\langle \varphi, b\text{-}\Ch_{t_1}(\widetilde D) \big\rangle - \big\langle \varphi,, b\text{-}\Ch_{t_0}(\widetilde D) \big\rangle
=  - \frac{1}{2}\int_{t_0}^{t_1} \big\langle \varphi,  \STR(\widetilde D_\partial e^{-(s\widetilde D_\partial +\nabla)^2})\big \rangle ds. \]
Let us write $\varphi=\varphi_e+\varphi_d$, where $\varphi_d$ is the delocalized part of $\varphi$, i.e.,  
$$\varphi_d(g_0,g_1,\cdots,g_{2m})=\begin{cases}
\varphi(g_0,g_1,\cdots,g_{2m})& \textup{ if } g_0g_1\cdots g_{2m}\ne e,\\
0& \text{ otherwise.}
\end{cases}$$
A similar argument as in the proof of Proposition $\ref{convergenceofTRhigher}$,  combined with Getzler's symbol calculus (cf. \cite{MR836727}),  shows that 
	$$\lim_{t\to 0}\big\langle \varphi_e,b\text{-}\Ch_t(\widetilde D)\big\rangle=\int_W \hat A\wedge\omega,$$
where $\omega$ is an element in $\Omega^\ast(W)\otimes \Omega_\ast(\mathfrak B_q^\omega)$, cf. \cite[Theorem 13.6]{Leichtnam}.	
Moreover,  it follows from the inequality in line \eqref{3.27} that 
		$$\lim_{t\to 0}\big\langle \varphi_d,b\text{-}\Ch_t(\widetilde D)\big\rangle=0.$$
Therefore, as $t_0\to 0$, we obtain the following formula:
		\begin{equation*}
	\big\langle\varphi,b\text{-}\Ch_{t_1}(\widetilde D)\big\rangle - \big\langle \varphi,\int_W \hat A\wedge\omega\big\rangle = - \frac{1}{2}\int_{0}^{t_1} \big\langle \varphi,  \STR(\widetilde D_\partial e^{-(s\widetilde D_\partial +\nabla)^2})\big \rangle ds.
	\end{equation*}
Now it follows from the discussion above that the integral on the right hand side converges absolutely as $t_1\to \infty$, under the condition that the spectral gap of $\widetilde D_\partial$ is larger than $\sigma_\varphi$. This finishes the proof.
\end{proof}

\begin{remark}
	Formally speaking, the term $\lim_{t\to\infty}\big\langle \varphi, b\textup{-}\Ch_t(\widetilde D))\big\rangle$ represents the pairing between the higher index class\footnote{In general, there is no natural way to define the higher index class of a Dirac operator on a manifold with boundary. However, in our setup above, due to the invertibility of the operator  $\widetilde D_\partial$ on the boundary, there is a natural higher index class associated to $\widetilde D$.} of $\widetilde D$  and the cyclic cocycle $\varphi$. However, to make this formal assertion rigorous, one needs to extend the the pairing in $\eqref{eq:b-APS}$ from $\mathfrak B_q^\omega$ to a smooth dense subalgebra of $C_r^\ast(G)$, which is a rather subtle issue in general. In the remaining sections below, we will show the existence of such an extension of the pairing, in the special case where $G$ is hyperbolic.
\end{remark}

\section{Puschnigg smooth dense subalgebra for hyperbolic Groups}\label{SmoothDenseSubalgebras}
In this section, we review the construction of Puschnigg's smooth dense algebra of $C_r^*(G)$ for hyperbolic groups \cite{Puschnigg}. One particular feature is that  every trace on $\mathbb CG$ admits a continuous extension to of this Puschnigg smooth dense subalgebra,   cf. \cite[Theorem 5.2]{Puschnigg}. We shall generalize this extension result to  cyclic cocycles of all degrees.

\subsection{Unconditional seminorms and tensor products}
In this subsection, we review the construction of Puschnigg's smooth dense subalgebra of  $ C^*_r(G)$ for hyperbolic groups $G$ \cite{Puschnigg}. 

Let $X$ be a set and $R$ a normed algebra equipped with a sub-multiplicative norm $|\cdot|$. We denote by $RX$ the algebra consisting of all finitely supported functions on $X$ with values in $R$. For each element
$
A=\sum A_x x\in RX
$, we define its absolute value to be
\begin{equation*}
|A|=\sum |A_x|x\in\mathbb CX.
\end{equation*}
Define a partial order on elements in $RX$ by
\begin{equation*}
A\leqslant A'\iff |A_x|\leqslant |A'_x| \textup{ for }\forall x\in X.
\end{equation*}
Recall the following notion of unconditional seminorm due to Bost and Lafforgue (cf. \cite{L1}).
\begin{definition}
	A seminorm $\|\cdot\|$ on $RX$ is called unconditional if
	\begin{equation*}
	|A|\leqslant |A'| \implies \|A\|\leqslant \|A'\|,\textup{ for }  \forall A,A'\in RX.
	\end{equation*}
\end{definition}

Any seminorm $\|\cdot \|$ on $RX$ naturally determines an  unconditional seminorm $\|\cdot \|^+$ by 
\begin{equation}\label{eq:ucnormwrt}
\|A\|^+:=\inf_{|A'|\geqslant |A|} \left\|  |A'| \right\|.
\end{equation}

\begin{lemma}[{\cite[Lemma 2.3]{Puschnigg}}]\label{unconditionaloperatornorm}
	Let $X,Y$ be two sets and $\|\cdot\|_X,\|\cdot\|_Y$ be seminorms on $RX$ and $RY$ respectively.
	Let $\varphi:(RX,\|\cdot\|_X)\to(RY,\|\cdot\|_Y)$ be a bounded linear map. Assume that $\varphi$ is expressed by a positive integral kernel, that is
	\begin{equation}
	\varphi(\sum A_y y)(x)=\sum_{y\in Y} \varphi_{x,y}A_y y,
	\end{equation}
	where $\varphi_{x,y}\in\mathbb R_{\geqslant 0}$,  for $\forall x\in X,\forall y\in Y$.
	Then $\varphi$ is also bounded with respect to the corresponding unconditional seminorms 
	$\|\cdot\|_X^+,\|\cdot\|_Y^+$, and
	\begin{equation}
	\|\varphi\|^+\leqslant\|\varphi\|.
	\end{equation}
\end{lemma}

Now we recall the notion of unconditional tensor product seminorm.
\begin{definition}\label{unconditionaltensor}
	Let $X,Y$ be sets and $\|\cdot\|_X,\|\cdot\|_Y$ be unconditional seminorms on $RX$ and $R'Y$ respectively. Let $R\otimes R'$ be the algebraic tensor product of $R$ and $R'$ equipped with the projective seminorm.
	The unconditional tensor product seminorm $\|\cdot\|_{uc}$ on 
	$RX\otimes R'Y\cong (R\otimes R')(X\times Y)$ is defined to be the unconditional norm determined by this projective seminorm.  More precisely, $\|\cdot\|_{uc}$ is given by
	\begin{equation*}
	\|A\|_{uc}:=\inf_{|A|\leqslant\sum|A_i'|\otimes|A_i''|}
	\sum_i\|A'_i\|_X\|A''_i\|_Y,\textup{ for } \forall A\in RX\otimes R'Y.
	\end{equation*}
\end{definition}

This norm is less than or equal to projective seminorm over $RX\otimes R'Y$. An example where these two are not equal is given in \cite[Example 2.4]{Puschnigg}.

Given a finitely generated group $G$, we fix a symmetric generating set $S$ of $G$. Let $\ell$ be the corresponding word metric on $G$. 
In the following, let $\mathcal S$ be the collection of all trace class operators equipped with the trace norm $|\cdot|_1$. Let $\mathcal SG$ be the subalgebra of  $\mathcal K\otimes C_r^*(G)$ consisting of all finite sums $\sum A_g g$ with $A_g\in\mathcal S$.
\begin{definition}
	For any fixed $p\geq 1$, we define an unconditional norm $\|\cdot\|_{RD, p}$ on $\mathcal SG$ by 
	\begin{equation}\label{rdnorms}
	\|A\|_{RD,p}^2=\sum_g |A_g|_1^2(1+\ell (g))^{2p},
	\end{equation}
	for $ A=\sum_g A_g g\in \mathcal SG$. 
\end{definition} 
We denote the completion of $\mathcal SG$ with respect to $\|\cdot\|_{RD, p}$ by $RD_p(\mathcal SG)$. Similarly, the same formula also defines an unconditional norm $\|\cdot\|_{RD,p}$ on $\mathbb CG$. We denote the completion of $\mathbb CG$ under this norm by $RD_p(G)$. In the following, if no confusion is likely to arise, we shall omit $p$ from the notation.  
 
Let us assume $G$ is hyperbolic for the rest of this section. In this case, it is known that $RD(G)$ is a smooth dense subalgebra of $C_r^*(G)$, cf.  \cite{Jolissaint,Harpe,Laff}. Similarly, $RD(\mathcal SG)$ is a smooth dense algebra of $\mathcal K\otimes C_r^*(G)$.

Recall the following quasiderivation map defined by Puschnigg:
\begin{equation*}
\Delta:\mathcal SG\to	\mathcal SG\otimes\mathbb CG\cong \mathcal S(G\times G),\ 
A_g g\mapsto \sum_{\substack{g_1g_2=g\\\ell (g_1)+\ell (g_2)=\ell (g)}}A_g g_1\otimes g_2.
\end{equation*}

\begin{definition}\label{def:Bpnorm}
Let $\|\cdot \|_{B, p}$ be the norm on $\mathcal SG$ given by 
\begin{equation}\label{norminB(RG)}
	\|A\|_{B, p}:=\|A\|_{RD, p}+\|\Delta A\|_{uc},\ \forall A\in\mathcal SG.
	\end{equation}
	Here $\|\cdot\|_{uc}$ is the unconditional tensor product norm on $\mathcal SG\otimes \mathbb CG\cong\mathcal S(G\times G)$ determined by the unconditional norm $\|\cdot\|_{RD, p}$ on both $\mathcal SG$ and $\mathbb CG$.
\end{definition}

Let $B_p(\mathcal SG)$ be the completion of $\mathcal SG$ with respect to $\|\cdot\|_{B,p}$. Apply the same construction to $\mathbb CG$ and we obtain $B_p(\mathbb CG)$. If no confusion is likely to arise, we shall omit $p$ from the notation. 

We define a more flexible quasiderivation as follows.
\begin{definition}\label{flexibleqd}
	For any $g\in G$ and $q\geqslant0$, let $C(q,g)$ be the collection of all pairs $(g_1,g_2)\in G\times G$ satisfying the following conditions:
	\begin{enumerate}[label=(\arabic*)]
		\item $g_1g_2=g$,
		\item there exists a geodesic $[e,g]$ connecting the identity $e$ and $g$ in the Cayley graph of $G$ such that $g_1$ lies in the $q$-neighborhood of $[e,g]$.
	\end{enumerate}  
We define
	\begin{equation*}
	\Delta_q\colon \mathbb CG\to	\mathbb CG\otimes\mathbb CG\cong\mathbb C(G\times G),\ 
	g\mapsto \sum_{(g_1,g_1)\in C(q,g)} g_1\otimes g_2.
	\end{equation*}
\end{definition}

If $q=0$, then $\Delta_0$ agrees with $\Delta$. By definition, for any $(g_1,g_2)\in C(q,g)$, there exists a group element $v\in G$ with $\ell(v) \leqslant  q$ such that $\ell (g_1v^{-1})+\ell (vg_2)=\ell (g)$. For each $v\in G$, define a map
\begin{equation*}
s_v:\mathbb CG\otimes \mathbb CG\to\mathbb CG\otimes \mathbb CG,\ 
g_1\otimes g_2\mapsto g_1v\otimes v^{-1}g_2.
\end{equation*}
Then we have  for $A\in \mathbb CG$
\begin{equation*}
\Delta_q|A|\leqslant \sum_{\ell (v)\leqslant q}s_v\Delta|A|.
\end{equation*}
By Lemma \ref{unconditionaloperatornorm}, the operator norm of $s_v$ (with respect to the unconditional norm $\|\cdot\|_{uc}$) does not exceed $(1+\ell (v))^2$. Since the number of elements in $G$ of length  $\leqslant q$ is finite, we see that there exists a constant $K_q$ such that $\|(\Delta_q|A|)\|_{uc}\leqslant K_q\cdot \|(\Delta|A|)\|_{uc}$ for all $\in\mathbb CG$.

\begin{proposition}{\cite[Proposition 3.5]{Puschnigg}}\label{estimationondelta}
	If $G$ is a hyperbolic group whose Cayley graph is $\delta$-hyperbolic, then there exists $C>0$ such that 
	\begin{equation*}
	\|\Delta(AA')\|_{uc}\leqslant C(\|\Delta(A)\|_{uc}\|A'\|_{RD}+
	\|A\|_{RD}\|\Delta(A')\|_{uc}),
	\end{equation*}
	for  $A,A'\in\mathcal SG$
\end{proposition}
\begin{proof}
	By the discussion above, it suffices to  prove the following pointwise inequality 
	\begin{equation}\label{derivationineq}
	\Delta |AA'|\leqslant \Delta_\delta (|A|)( 1\otimes |A'|)+(|A|\otimes 1)\Delta_\delta (|A'|),
	\end{equation}
	for $A,A'\in\mathcal SG$. Without loss of generality,  it suffices to consider the case where $A=g$ and $A'=g'$, for $g,g'\in G$.
	
	Let $k=gg'$. If a term $k_1\otimes k_2$ appears in the summation expression of $\Delta |AA'|$, then $k_1$ is a point on the geodesic $[e,k]$. Since the Cayley graph of $G$ is $\delta$-hyperbolic, $k_1$ lies in the $\delta$-neighborhood of the union of $[e,g]$ and $[g,k]$. Either there is a group element $g_1\in[e,g]$ such that $\dist(k_1,g_1)<\delta$, or there is a group element $g_2\in[g,k]$ such that $\dist(k_1,g_2)<\delta$. We prove the former case; the latter case is similar. In the former case,   we see that the term $k_1\otimes k_1^{-1}g$ appears in the summation expression of $\Delta_\delta(g)$. This implies that the term $k_1\otimes k_2=(k_1\otimes k_1^{-1}g)(1\otimes g')$ appears in $\Delta_\delta (|A|)( 1\otimes |A'|)$. This finishes the proof.
\end{proof}
\begin{remark}
	The above proof in fact shows that the pointwise inequality in line $\eqref{derivationineq}$ is equivalent to the hyperbolicity of the group.
\end{remark}
\begin{proposition}\label{smoothdensesubalgebra}
	$B(\mathcal SG)$ and $B(\mathbb CG)$ are smooth dense algebras of $\mathcal K\otimes C^*_r(G)$ and $C_r^*(G)$ respectively.
\end{proposition}
\begin{proof}
	We prove the case of  $B(\mathbb CG)$; the other case is similar. 
	
	Since $RD(G)$ is a smooth dense subalgebra of $C_r^\ast(G)$,  it suffices to show that if an element $T\in B(\mathbb CG)^+$ is invertible in $RD(G)^+$, then $T$ is invertible in $B(\mathbb CG)^+$. In fact, it suffices to show that there exists a constant $\varepsilon>0$ such that if an element $T\in B(\mathbb CG)^+$ satisfies $\|T -1\|_{RD} < \varepsilon$, then $T$ is invertible in $B(\mathbb CG)^+$. Indeed, let $S$ be an element in $B(\mathbb CG)^+$ such that $S$ is invertible in $RD(G)^+$ with inverse $R$. Since  $B(\mathbb CG)^+$ is dense in $RD(G)^+$, there exists an invertible element $U$ of $RD(G)^+$ such that $\|U - R\|_{RD} < \varepsilon \cdot \|S\|_{RD}$. It follows that $SU$ is invertible in $RD(G)^+$ and $\|SU-1\|_{RD} < \varepsilon$. Then by our assumption, $SU$ is invertible in $B(\mathbb CG)^+$, which implies $S$ is invertible in $B(\mathbb CG)^+$.
	
	 Now suppose $A\in B(\mathbb CG)^+$ such that $\|A-1\|_{RD}<\min\{1/C,1\}$, where $C$ is the same constant as in Proposition \ref{estimationondelta}. It follows from Proposition \ref{estimationondelta} that
	$$\|\Delta((A-1)^n)\|_{uc}\leqslant n(C\|A-1\|_{RD})^{n-1}\|\Delta(A-1)\|_{uc}.$$
	This immediately implies that 
	$A^{-1}= (1- (1-A))^{-1} = \sum_{n=0}^\infty (1-A)^n$ lies in $B(\mathbb CG)^+.$
	Therefore $B(\mathbb CG)$ is a smooth dense subalgebra of $RD(G)$. 
\end{proof}
\subsection{Continuous extension of traces}
In this subsection, we review Puschnigg's result on continuous extension of traces from $\mathbb CG$ to $B(\mathbb CG)$ for hyperbolic groups  \cite[Theorem 5.2]{Puschnigg}. In fact, for the purposes of this paper, we only need a weaker version of Pushnigg's theorem, to which we give a slightly different proof. 

\begin{lemma}[{\cite[Lemms 4.1]{Puschnigg}}]\label{estonconjugacy}
	Let $G$ be a group whose Cayley graph is  $\delta$-hyperbolic. Given $h\in G$, if $g$ lies in the conjugacy class $\left\langle h\right\rangle$, then there exist $g_1,g_2\in G$ such that $g_1g_2=g$, $\ell (g_1)+\ell (g_2)=\ell (g)$ and $\ell (g_2g_1)\leqslant 6\delta+6+ 3\ell (h)$.
\end{lemma}
\begin{proof}
	Suppose $h=ugu^{-1}$ for some $u\in G$. In the following, we denote  by $[a, b]$ a geodesics connecting  $a, b\in G$. By hyperbolicity, there exist vertices $w\in[e,ug]$, $u_1\in [e,u]$ and $ug_1\in[u,ug]$ such that $\dist(w,u_1) <\delta+1$ and $d(w,ug_1)<\delta+1$. Moreover, there exists $v_1\in G$ such that  $ugv_1\in [ug, h]$  and $\dist(w, ugv_1)< \delta +1 +\ell (h)$, since the $[e, ug]$ lies entirely in the $\delta + \ell(h)$ neighborhood of $[ug, h]$. Let $u_2, g_2, v_2$ be elements in $G$ such that $u = u_1u_2$, $g=g_1g_2$ and $u^{-1}=v_1v_2$. Let us write  $$g_2g_1=(g_2v_1)(v_1^{-1}u_2^{-1})(u_2g_1).$$
	Clearly, we have that 
	\begin{align*}
	\ell (g_2v_1)&=\dist(ug_1,ug_1g_2v_1) = \dist(ug_1, ugv_1)< 2\delta+ 2+ \ell (h),\\
	\textup{ and } \ell (u_2g_1)&=\dist(u_1,u_1u_2g_1) = \dist(u_1,ug_1)< 2\delta + 2.
	\end{align*}
Furthermore, observe that $v_2^{-1}$ is a vertex on the geodesic $[e,u]$. It follows that  $$\ell (v_1^{-1}u_2^{-1})=\ell (v_2u_1)=\dist(v_2^{-1},u_1) = |\ell (v_2^{-1})-\ell (u_1)|.$$
Therefore, we have 
\begin{align*}
\ell (v_1^{-1}u_2^{-1})& =|\ell (v_2^{-1})-\ell (u_1)|  = |\dist(h, hv_2^{-1}) -\dist (e, u_1)|\\
& \leqslant | \dist(e, h) + \dist(e, hv_2^{-1}) -\dist (e, u_1)| \\
& \leqslant | \dist(e, h)| + |\dist(u_1, hv_2^{-1})| \\
& < 2\delta + 2+ 2\ell(h).
\end{align*}
 This finishes the proof. 
\end{proof}

\begin{theorem}[{\cite[Theorem 5.2]{Puschnigg}}]\label{extensionoftrh}
	Let $G$ be a hyperbolic group and $B(\mathbb CG)$  the Puschnigg  smooth dense subalgebra of $C_r^*(G)$. For any conjugacy class $\left\langle h\right\rangle$ of $G$, the map 
	$$tr_{\left\langle h\right\rangle }:\mathbb CG\to\mathbb C,\ \sum_{g\in G}a_g g\mapsto \sum_{g\in\left\langle h\right\rangle }a_g.$$
	admits a continuous extension to $B(\mathbb CG)$.
\end{theorem}
\begin{proof}
	Define a map $\mu:\mathbb CG\to\mathbb C(G\times G)$ as follows: if $g\in \left\langle h\right\rangle$, then 
	\[  \mu(g)\coloneqq g_1\otimes g_2\] where $(g_1,g_2)\in G\times G$ is a pair of elements as given in Lemma \ref{estonconjugacy}; if $g\notin \langle h\rangle$, define $\mu(g)=0$. Clearly,  $\mu(|A|)\leqslant \Delta |A|$. Thus $\mu$ admits a continuous extension from $B(\mathbb CG)$ to $RD(G)\otimes_{uc}RD(G)$, which we will still denote by $\mu$.
	
	By Lemma \ref{unconditionaloperatornorm}, the maps 
	\begin{align*}
	&\sw\colon RD(G)\otimes_{uc} RD(G)\to RD(G)\otimes_{uc} RD(G),\ g_1\otimes g_2\mapsto g_2\otimes g_1,\\
	&\m \colon RD(G)\otimes_{uc} RD(G)\to RD(G),\ g_1\otimes g_2\mapsto g_1g_2,
	\end{align*}
	are continuous. 
	
	We define an evaluation map $\mathrm{E}\colon RD(G)\to\mathbb C$ as follows:
	\[ \mathrm{E}(g) = \begin{cases*}
	   1 &  \textup{ if $\ell (g)\leqslant 6\delta+ 6 + 3\ell (h)$}, \\
	  0 & \textup{ otherwise.}
	\end{cases*} \]  
	Clearly,  $\mathrm E$ is also well-defined and continuous. It follows that the composition 
	$$B(\mathbb CG)\xrightarrow{\mu}RD(G)\otimes_{uc} RD(G)
	\xrightarrow{\ \sw\ }RD(G)\otimes_{uc} RD(G)\xrightarrow{\ \m\ } RD(G)\xrightarrow{\mathrm{E}}\mathbb C$$
	is a continuous extension of $\tr_{\left\langle h\right\rangle}$. This finishes the proof. 
\end{proof}

The previous theorem has the following obvious analogue where the coefficient $\mathbb C$ is replaced by the algebra of trace class operators $\mathcal S$. 
\begin{proposition}\label{extendTR}
	Let $B(\mathcal SG)$ be the smooth dense subalgebra of $\mathcal K\otimes C_r^*(G)$ defined above. For any conjugacy class $\left\langle h\right\rangle $ in $G$,  let $\tr_{\left\langle h\right\rangle}\colon \mathcal SG \to \mathbb C$ be the trace map defined by
	\begin{equation*}
	\tr_{\left\langle h\right\rangle}(A)=\sum_{g\in\left\langle h\right\rangle} \tr(A_g),\ \textup{ for } A=\sum A_g g\in\mathcal SG.
	\end{equation*}
	Then $\tr_{\langle h\rangle}$ extends to a continuous trace map on $B(\mathcal SG)$.
\end{proposition}

\subsection{Continuous extension of higher degree cyclic cochains}\label{sec:higherext}
In this subsection, we  generalize the continuous extension result for traces to higher degree cyclic cochains.

\begin{definition}
	Fix a length function $\ell$ on $G$. For any $\varphi\in C^n(\mathbb CG,\langle h\rangle)$, we say $\varphi$ has polynomial growth if there exist constants $C$ and $k$ such that
	$$|\varphi(g_0,g_1,\cdots,g_n)|\leqslant C\prod_{i=0}^n(1+\ell (g_i))^k.$$
\end{definition}
In Section $\ref{CyclicCohomologyGrowth}$ below, we will show that, when $G$ is hyperbolic,  every element in $HC^n(\mathbb CG,\langle h\rangle)$ have a representative with polynomial growth. 

 Denote by $(\mathcal SG)^{\otimes n+1}$ the algebraic tensor product of $(n+1)$ copies of $\mathcal SG$. Recall the unconditional tensor product defined in Definition \ref{unconditionaltensor}.
 We construct the unconditional tensor product norms $\|\cdot\|_{RD}$ and $\|\cdot\|_{B}$ over $(\mathcal SG)^{n+1}$, and denote their completions by $RD(\mathcal SG)^{\otimes n+1}$ and $B(\mathcal SG)^{\otimes n+1}$ respectively.

\begin{proposition}\label{highertrace}
Let $G$ be a hyperbolic group whose Cayley graph is $\delta$-hyperbolic.	If $\varphi\in C^n(\mathbb CG,\langle h\rangle)$ has polynomial growth, then the map $\varphi\#\tr\colon (\mathcal SG)^{\otimes n+1} \to \mathbb C$ given by 
	\begin{equation}
	\begin{split}
	A_{g_0} g_0\otimes A_{g_1} g_1 \otimes\cdots\otimes A_{g_n}g_n&\mapsto 
	\tr(A_{g_0}A_{g_1}\cdots A_{g_n})\varphi(g_0,g_1,\cdots,g_n)
	\end{split}
	\end{equation}
	extends continuously to  $B(\mathcal SG)^{\otimes n+1}$.
\end{proposition}
\begin{remark}
	Before we prove the proposition, let us point out that the construction of $B_p(\mathcal SG)$  involves a choice of some sufficiently large $p$. In order to extend $\varphi\#\tr$ continuously to $B_p(\mathcal SG)^{\otimes n+1}$,   we assume that $p$ is sufficiently large so that it ``dominates"  the growth rate of $\varphi$. Hence, strictly speaking, the algebra $B(\mathcal SG)^{\otimes n+1}$ may vary for different cyclic cochains.
\end{remark}
\begin{proof}
	Suppose that 
	$$|\varphi(g_0,g_1,\cdots,g_n)|\leqslant C\prod_{i=0}^n(1+\ell (g_i))^k$$ 
	for $(g_0,g_1,\cdots,g_n)\in G^{n+1}$. 
	
	Define the following maps:
	\begin{enumerate}
		\item $\pi_\varphi:(\mathcal SG)^{\otimes n+1} \to (\mathcal SG)^{\otimes n+1} $ by 
		\begin{align*}
		A_{g_0} g_0\otimes A_{g_1} g_1\otimes\cdots\otimes A_{g_n}g_n&\mapsto 
		\varphi(g_0,g_1,\cdots,g_n)A_{g_0} g_0\otimes A_{g_1} g_1\otimes\cdots\otimes A_{g_n}g_n;
		\end{align*}
		\item $\m:\mathcal (\mathcal SG)^{\otimes n+1}  \to \mathcal SG$ by
		\begin{align*}
		A_{g_0} g_0\otimes A_{g_1} g_1\otimes\cdots\otimes A_{g_n}g_n&\mapsto 
		A_{g_0}A_{g_1}\cdots A_{g_n}g_0g_1\cdots g_n.
		\end{align*}
	\end{enumerate} 
	Clearly, the composition
	$$(\mathcal SG)^{\otimes n+1}\xrightarrow{\ \pi_\varphi\ }(\mathcal SG)^{\otimes n+1}\xrightarrow{\ \m\ }\mathcal SG\xrightarrow{\tr_{\langle h\rangle}}\mathbb C$$
	is exactly the map $\varphi\#\tr$.
Therefore, it suffices to show that $\m\circ\pi_\varphi$ extends to a continuous map 
	$B_p(\mathcal SG)^{\otimes n+1}$ to $B_{p-k}(\mathcal SG)$. It follows from Lemma \ref{unconditionaloperatornorm} that the map $\m\circ \pi_\varphi$ extends to a continuous map from $RD_p(\mathcal SG)^{\otimes n+1}$ to $ RD_{p-k}(\mathcal SG)$. 
	
In the following, let us prove the case where $k=0$, that is, $\varphi$ is uniformly bounded over $G^{n+1}$. The general case is similar. If $\varphi$ is uniformly bounded, then for any $A\in(\mathcal SG)^{\otimes n+1}$, we observe that
	$$|\Delta (\m\circ \pi_\varphi A)|=\Delta(|\m\circ \pi_\varphi A|)\leqslant \Delta (\m|\pi_\varphi A|)\leqslant
	C\cdot \Delta (\m|A|).$$
	
	Define another multiplication map $	\m'
	\colon \mathbb CG^{2(n+1)}\to\mathbb CG^2$ by 
	\begin{equation*}
	g_0\otimes g_0'\otimes g_1\otimes g_1'\otimes\cdots\otimes g_n\otimes g_n'\mapsto
	g_0g_1\cdots g_n\otimes g_0'g_1'\cdots g_n'.
	\end{equation*}
	This is also bounded with respect to the unconditional norm. We claim that
	\begin{equation}\label{4.22}
	\Delta (\m |A|)\leqslant \m'\Delta_{n\delta}^{\otimes (n+1)}(|A|),\ \forall A\in\mathbb CG^{n+1}.
	\end{equation}
	Here $\Delta_{n\delta}$ is  the quasiderivation from  Definition \ref{flexibleqd}, and $\Delta_{n\delta}^{\otimes (n+1)}$ stands for the tensor product of $(n+1)$-copies of $\Delta_{n\delta}$ from $\mathbb CG^{n+1}$ to $\mathbb C(G^2)^{(n+1)}\cong\mathbb CG^{2(n+1)}$.
	
	Assume the claim holds for the moment. Clearly, the  map $\m'$ is bounded with respect to the unconditional norm. Moreover, by the discussion before Proposition $\ref{estimationondelta}$, there exists a constant $K$ such that   \[ \|\Delta_{n\delta}^{\otimes (n+1)}(|A|)\|_{uc} \leqslant K \|\Delta^{\otimes (n+1)}(|A|)\|_{uc}\]
	for all $A\in \mathbb CG^{n+1}$. This proves the proposition. 
	
Now let us prove the claim.  It suffices to prove the inequality $\eqref{4.22}$ when $|A|=g_0\otimes\cdots\otimes g_n$. Denote $g_0g_1\cdots g_n$ by $g$. Suppose  $g'\otimes g''$ appears on the left-hand side of the inequality $\eqref{4.22}$, where  by definition $g'$ is a point on the geodesic $[e,g]$. We will show that $g'\otimes g''$ also appears on the right-hand side of the inequality \eqref{4.22}. Indeed, by hyperbolicity, there exists a point $x$ on the path $[e,g_0],[g_0,g_0g_1],\cdots, [g_0g_1\cdots g_{n-1},g]$ such that the distance from $x$ to $g'$ is less than $ n\delta$. More precisely, there exist $j\geqslant 0$ and $v, v'\in G$ such that $vv'=g_j$, $\ell (v)+\ell (v')=\ell (g_j)$ and $d(g',g_0g_1\cdots g_{j-1}v)=d(v,(g_0g_1\cdots g_{j-1})^{-1}g')< n\delta$.
	Thus the following element 
	\begin{align*}
	&(g_0\otimes 1)\otimes\cdots\otimes (g_{j-1}\otimes 1)
	\otimes\big( (g_0g_1\cdots g_{j-1})^{-1}g'\otimes g''(g_{j+1}g_{j+2}\cdots g_n)^{-1}\big)\\
	&\otimes (1\otimes g_{j+1})\otimes\cdots\otimes (1\otimes g_n)
	\end{align*}
	appears in the summation expression of $\Delta_{n\delta}^{\otimes n+1}(g_0\otimes g_1\otimes\cdots\otimes g_n)$. After applying the map $\m'$, we see that  $g'\otimes g''$ indeed appears on the right-hand side of the inequality \eqref{4.22}. This proves the claim, hence finishes the proof of the proposition. 
\end{proof}

\section{Cyclic gohomology of hyperbolic Groups}\label{CyclicCohomologyGrowth}

In this section, we show that every cyclic cohomology class of a hyperbolic group  has a uniformly bounded representative if its degree is $\geq 2$. Since for any group, the equivalence class of a cyclic cocycle of degree $\leq 1$ always has a representative of polynomial growth, it follows that all cyclic cohomology classes of  a hyperbolic group can be represented by cyclic cocycles of polynomial growth.

We will need the following results on the geometry of hyperbolic groups \cite{Gromov}. Suppose $G$ is a word hyperbolic group. For each $h\in G$, let $G_h$ be the centralizer of $h$ in $G$, and $N_h$ the quotient of $G_h$ by the cyclic group generated by $h$. 
\begin{enumerate}[label=(\textbf{G}\arabic*)]
	\item \label{finiteorder} If $h\in G$ has infinite order, then $N_h$  is finite.
	\item \label{quasiconvex}
	For any $h\in G$, the centralizer $G_h$ is a quasi-convex subspace of $G$, that is, there exists some $K>0$ such that any geodesic in $G$ connecting a pair of points in $G_h$ lies in a $K$-neighborhood of $G_h$.
	\item \label{centralizerishyperbolic}
	For any $h\in G$, its centralizer $G_h$ is also word hyperbolic, and the inclusion $G_h \hookrightarrow G$ is a quasi-isometry.
\end{enumerate}
Moreover, we will use the following result of Mineyev \cite[Theorem 11]{Mineyev} in an essential way. 
\begin{enumerate}[label=(\textbf{M}\arabic*)]
	\item  \label{polcohomology} Suppose $G$ is a word hyperbolic group with a given length function $\ell$. If $n\geqslant 2$, then every element in  $H^n(G;\mathbb C)$ --- the group cohomology of $G$ --- admits a uniformly bounded representative. Here a cocycle element $\varphi$ is said to be uniformly bounded if there exists $C>0$ such that $
	|\varphi(g_0,g_1,\cdots,g_n)|\leqslant C$ for $g_i\in G$. 
\end{enumerate}

We will also need the following lemma, the proof of which is communicated to us by Denis Osin. 

\begin{lemma}\label{lm:conj}
	Let $G$ be a $\delta$ hyperbolic group with a word length function $\ell$. For each element $h\in G$, there exists a constant $K_h >0$  such that 
	\[  \min \{\ell(\gamma) \mid \gamma^{-1} h \gamma = g\} \leqslant \ell(g) + K_h. \]
\end{lemma}
\begin{proof}
	Let $\beta$ be a group element of minimal length such that $\beta^{-1} h \beta = g$. Consider the geodesic quadrilateral $[e, \beta^{-1}], [\beta^{-1}, \beta^{-1} h], [\beta^{-1} h, g]$ and $[g, e]$ in the Cayley graph of $G$ (see Figure $\ref{fig:geo}$ below). 
	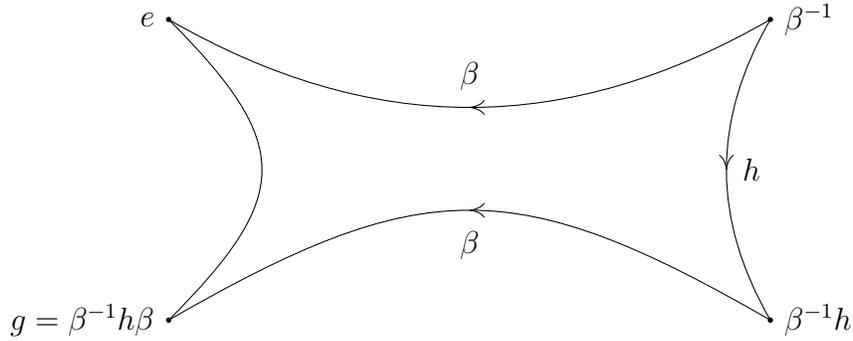
\begin{figure}[h]
		\centering
		\begin{tikzpicture}[decoration={
			markings,
			mark=at position 0.5 with {\arrow{To[length=2mm, width=2mm]}}}]
		\node [fill, circle, inner sep=0pt,minimum size=2pt, label={left: $e$}] (0) at (-4, 2) {};
		\node [fill, circle, inner sep=0pt,minimum size=2pt, label={left: $g =\beta^{-1}h\beta$}] (1) at (-4, -2) {};
		\node [fill, circle, inner sep=0pt,minimum size=2pt, label={right: $\beta^{-1}$}] (2) at (4, 2) {};
		\node [fill, circle, inner sep=0pt,minimum size=2pt, label={right: $\beta^{-1}h$}] (3) at (4, -2) {};
		\node [] (4) at (0, 1.25) {$\beta$};
		\node [] (5) at (0, -1) {$\beta$};
		\node [] (6) at (3.75, 0) {$h$};

		\draw [postaction={decorate}, bend right, looseness=1.25]  (3.center) to (1.center);
		\draw [postaction={decorate}, bend left]  (2.center) to (0.center);
		\draw [bend left=45, looseness=1.50] (0.center) to (1.center);
		\draw [postaction={decorate}, bend right] (2.center) to (3.center);
		\end{tikzpicture}
		\caption{Geodesic quadrilateral.} \label{fig:geo}
	\end{figure}

	By continuity, there exists a point $x$ on $[e, g]$ such that $x$ is equidistant from the two sides $[e, \beta^{-1}]$ and $[\beta^{-1} h, g]$. Let $y$ (resp. $z$) be a closest point on $[e, \beta^{-1}]$ (resp. $[\beta^{-1}h, g]$) to $x$, that is,  $d(x, y) = d(x, z)$ equals the distance between $x$ and the geodesic $[e, \beta^{-1}]$.    By hyperbolicity, it is not difficult to see that \[ d(x, y) = d(x, z) \leqslant 2\delta + \ell(h), \] cf. Figure $\ref{fig:hyper}$ below. It follows immediately that  both $d(e, y)$ and $d(z, g)$ are less than $2\delta + \ell(h) + \ell(g)$. 
	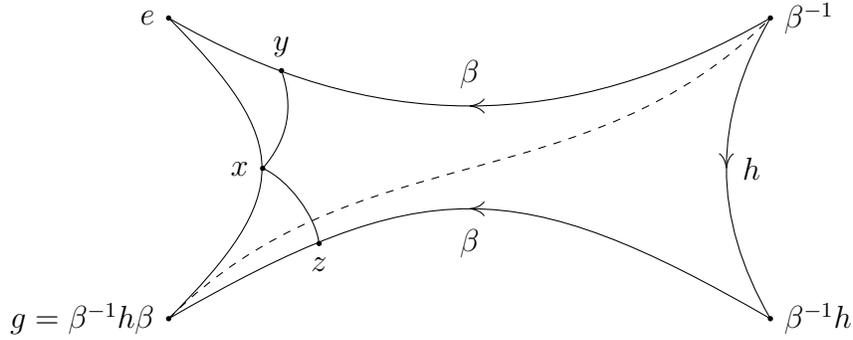
\begin{figure}[h]
		\centering 
		\begin{tikzpicture}[decoration={
			markings,
			mark=at position 0.5 with {\arrow{To[length=2mm, width=2mm]}}}]
		\node [fill, circle, inner sep=0pt,minimum size=2pt, label={left: $e$}] (0) at (-4, 2) {};
		\node [fill, circle, inner sep=0pt,minimum size=2pt, label={left: $g =\beta^{-1}h\beta$}] (1) at (-4, -2) {};
		\node [fill, circle, inner sep=0pt,minimum size=2pt, label={right: $\beta^{-1}$}] (2) at (4, 2) {};
		\node [fill, circle, inner sep=0pt,minimum size=2pt, label={right: $\beta^{-1}h$}] (3) at (4, -2) {};
		\node [] (4) at (0, 1.25) {$\beta$};
		\node [] (5) at (0, -1) {$\beta$};
		\node [] (6) at (3.75, 0) {$h$};
		\node [fill, circle, inner sep=0pt,minimum size=2pt, label={left: $x$}] (7) at (-2.75, 0) {};
		\node [fill, circle, inner sep=0pt,minimum size=2pt, label={above: $y$}] (8) at (-2.5, 1.3) {};
		\node [fill, circle, inner sep=0pt,minimum size=2pt, label={below: $z$}] (9) at (-2, -1) {};

		\draw [postaction={decorate}, bend right, looseness=1.25]  (3.center) to (1.center);
		\draw [postaction={decorate}, bend left]  (2.center) to (0.center);
		\draw [bend left=45, looseness=1.50] (0.center) to (1.center);
		\draw [postaction={decorate}, bend right] (2.center) to (3.center);
		\draw [bend left] (8.center) to (7.center);
		\draw [bend left, looseness=0.75] (7.center) to (9.center);
		\draw [dashed, in=-135, out=45] (1.center) to (2.center);
		
		\end{tikzpicture}
		\caption{$d(x, y) = d(x, z) \leqslant 2\delta + \ell(h)$.} \label{fig:hyper}
	\end{figure}

	It remains to estimate the length of either $[y, \beta^{-1}]$ or  $[z, \beta^{-1}h]$. Suppose  $\beta_1 \beta_2 \cdots \beta_m$ is a word of minimal length that represents $\beta$. Let $t_i$ be a geodesic connecting the points labeled by $\beta_i$ on $[e, \beta^{-1}]$ and $[\beta^{-1}h, g]$ (see Figure $\ref{fig:est}$ below). 
		\begin{figure}[h]
	\centering 
	\begin{tikzpicture}[decoration={
		markings,
		mark=at position 0.5 with {\arrow{To[length=2mm, width=2mm]}}}]
	\node [fill, circle, inner sep=0pt,minimum size=2pt, label={left: $e$}] (0) at (-4, 2) {};
	\node [fill, circle, inner sep=0pt,minimum size=2pt, label={left: $g =\beta^{-1}h\beta$}] (1) at (-4, -2) {};
	\node [fill, circle, inner sep=0pt,minimum size=2pt, label={right: $\beta^{-1}$}] (2) at (4, 2) {};
	\node [fill, circle, inner sep=0pt,minimum size=2pt, label={right: $\beta^{-1}h$}] (3) at (4, -2) {};
	\node [] (4) at (0, 1.25) {};
	\node [] (5) at (0, -1) {};
	\node [] (6) at (3.75, 0) {$h$};
	\node [fill, circle, inner sep=0pt,minimum size=2pt, label={left: $x$}] (7) at (-2.75, 0) {};
	\node [fill, circle, inner sep=0pt,minimum size=2pt, label={above: $y$}] (8) at (-2.5, 1.3) {};
	\node [fill, circle, inner sep=0pt,minimum size=2pt, label={below: $z$}] (9) at (-2, -1) {};
	\node [fill, circle, inner sep=0pt,minimum size=2pt, label={above: $\beta_1$}] (10) at (3, 1.5) {};
	\node [fill, circle, inner sep=0pt,minimum size=2pt, label={below: $\beta_1$}] (11) at (3, -1.45) {};
	\node [fill, circle, inner sep=0pt,minimum size=2pt, label={above: $\beta_2$}] (12) at (2, 1.12) {};
	\node [fill, circle, inner sep=0pt,minimum size=2pt, label={below: $\beta_2$}] (13) at (2, -0.98) {};
	\node [fill, circle, inner sep=0pt,minimum size=2pt, label={above: $\beta_{k}$}] (14) at (-0.75, 0.86) {};
	\node [fill, circle, inner sep=0pt,minimum size=2pt, label={below: $\beta_{k}$}] (15) at (-0.4, -0.56) {};
	\node [fill, circle, inner sep=0pt,minimum size=2pt, label={above: $\beta_3$}] (18) at (1.25, 0.94) {};
	\node [fill, circle, inner sep=0pt,minimum size=2pt, label={below: $\beta_3$}] (19) at (1.25, -0.73) {};
	\node [] (20) at (0.5, 0.25) {$\cdots$};
	\node  (21) at (-1.4, 0.25) {$\cdots $};
	\node  (22) at (-0.42, 0.15) {$t_{k}$};
	\node  (23) at (1.35, 0) {$t_3$};
	\node  (24) at (2.05, 0) {$t_2$};
	\node  (25) at (2.9, 0) {$t_1$};

	\draw [postaction={decorate}, bend right, looseness=1.25]  (3.center) to (1.center);
	\draw [postaction={decorate}, bend left]  (2.center) to (0.center);
	\draw [bend left=45, looseness=1.50] (0.center) to (1.center);
	\draw [postaction={decorate}, bend right] (2.center) to (3.center);
	\draw [bend left] (8.center) to (7.center);
	\draw [bend left, looseness=0.75] (7.center) to (9.center);
	\draw [dashed, bend left=15] (8.center) to (9.center);
	\draw [bend right=15] (14.center) to (15.center);
	\draw [bend right=20] (12.center) to (13.center);
	\draw [bend right=25] (10.center) to (11.center);
	\draw [bend right=15] (18.center) to (19.center);
	\end{tikzpicture}
	\caption{Estimates for the length of $[y, \beta^{-1}]$. } \label{fig:est}
\end{figure}
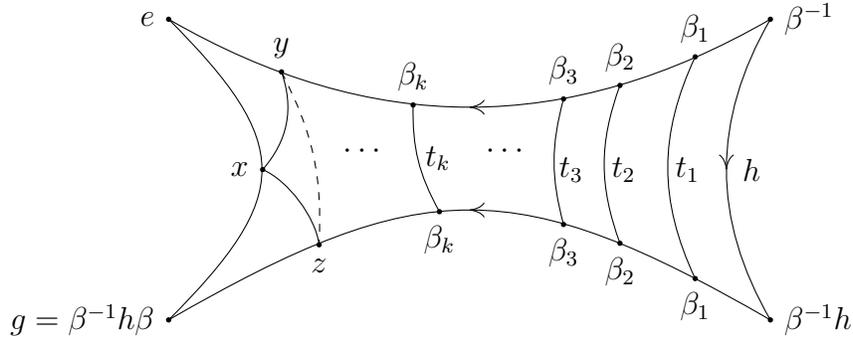	
	If $t_i = t_j$ for some $i< j$, we can cut the shaded region in Figure $\ref{fig:short}$ and  it follows that the element $\alpha = \beta_1\cdots \beta_i \beta_{j+1} \cdots \beta_m$  satisfies that  $\alpha^{-1} h \alpha = g$ and $\ell(\alpha) < \ell(\beta)$.    But  this contradicts the assumption that $\beta$  is a group element of minimal length such that $\beta^{-1} h \beta = g$.  Thus all $t_i$'s are pairwise distinct.

	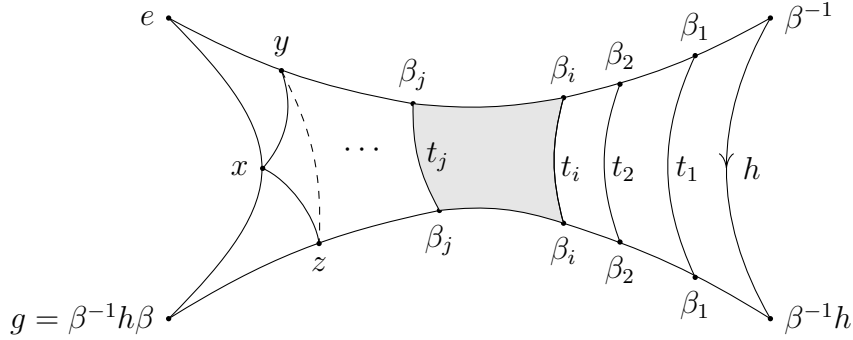
\begin{figure}[h]
		\centering 
		\begin{tikzpicture}[decoration={
			markings,
			mark=at position 0.5 with {\arrow{To[length=2mm, width=2mm]}}}]
		\node [fill, circle, inner sep=0pt,minimum size=2pt, label={left: $e$}] (0) at (-4, 2) {};
		\node [fill, circle, inner sep=0pt,minimum size=2pt, label={left: $g =\beta^{-1}h\beta$}] (1) at (-4, -2) {};
		\node [fill, circle, inner sep=0pt,minimum size=2pt, label={right: $\beta^{-1}$}] (2) at (4, 2) {};
		\node [fill, circle, inner sep=0pt,minimum size=2pt, label={right: $\beta^{-1}h$}] (3) at (4, -2) {};
		\node [] (4) at (0, 1.25) {};
		\node [] (5) at (0, -1) {};
		\node [] (6) at (3.75, 0) {$h$};
		\node [fill, circle, inner sep=0pt,minimum size=2pt, label={left: $x$}] (7) at (-2.75, 0) {};
		\node [fill, circle, inner sep=0pt,minimum size=2pt, label={above: $y$}] (8) at (-2.5, 1.3) {};
		\node [fill, circle, inner sep=0pt,minimum size=2pt, label={below: $z$}] (9) at (-2, -1) {};
		\node [fill, circle, inner sep=0pt,minimum size=2pt, label={above: $\beta_1$}] (10) at (3, 1.5) {};
		\node [fill, circle, inner sep=0pt,minimum size=2pt, label={below: $\beta_1$}] (11) at (3, -1.45) {};
		\node [fill, circle, inner sep=0pt,minimum size=2pt, label={above: $\beta_2$}] (12) at (2, 1.12) {};
		\node [fill, circle, inner sep=0pt,minimum size=2pt, label={below: $\beta_2$}] (13) at (2, -0.98) {};
		\node [fill, circle, inner sep=0pt,minimum size=2pt, label={above: $\beta_{j}$}] (14) at (-0.75, 0.86) {};
		\node [fill, circle, inner sep=0pt,minimum size=2pt, label={below: $\beta_{j}$}] (15) at (-0.4, -0.56) {};
		\node [fill, circle, inner sep=0pt,minimum size=2pt, label={above: $\beta_i$}] (18) at (1.25, 0.94) {};
		\node [fill, circle, inner sep=0pt,minimum size=2pt, label={below: $\beta_i$}] (19) at (1.25, -0.73) {};
		\node  (21) at (-1.4, 0.25) {$\cdots $};
		
		\node  (23) at (1.35, 0) {$t_i$};
		\node  (24) at (2.05, 0) {$t_2$};
		\node  (25) at (2.9, 0) {$t_1$};

		
		\draw [ bend left= 9]  (14.center) to (0.center);
		\draw [bend right = 8]  (18.center) to (2.center);

		\draw [ bend left= 10]  (1.center) to (15.center);
		\draw [bend left = 7]  (19.center) to (3.center);

		\draw [bend left=45, looseness=1.50] (0.center) to (1.center);
		\draw [postaction={decorate}, bend right] (2.center) to (3.center);
		\draw [bend left] (8.center) to (7.center);
		\draw [bend left, looseness=0.75] (7.center) to (9.center);
		\draw [dashed, bend left=15] (8.center) to (9.center);
		\filldraw [color=black, fill=gray!20] (14.center) to [bend right=15] (15.center) to [bend left=12]  (19.center) to [bend left = 15]  (18.center) to [bend left = 8] (14.center);
		\draw [bend right=20] (12.center) to (13.center);
		\draw [bend right=25] (10.center) to (11.center);
		\draw [bend right=15] (18.center) to (19.center);

		\node  (22) at (-0.42, 0.15) {$t_{j}$};
		\end{tikzpicture}
		\caption{If $t_i = t_j$ for some $i<j$, we can shorten $\beta$.}\label{fig:short}
		
	\end{figure}
	
	 Using hyperbolicity on the quadrilateral with vertices $\{y, z, \beta^{-1}h, \beta^{-1}\}$ (cf. \cite[Chapter III.H, Lemma 1.15]{MR1744486}), it is not difficult to see that there exists a constant $C_1 >0$ such that the following is satisfied:  if  $\beta_i$ is to the right of both $y$ and $z$ as shown in Figure $\ref{fig:est}$, then  $\ell(t_i)\leq C_{1}$.  Here  $C_1$ only depends on $\ell(h)$ and $\delta$, and in particular is independent of $g$.   It follows immediately the length of either $[y, \beta^{-1}]$ or $[z, \beta^{-1}h]$ is $\leqslant C_2 + 1$, where $C_2$ is the number of elements of $G$ of length at most $C_1$. 
	
	Combining the above estimates together, we see that there exists a constant $K_h$ (only dependent on $h$ and $\delta$) such that  
	\[ \ell(\beta) \leq \ell(g) + K_h.\]
	This finishes the proof.

\end{proof}

\begin{theorem}\label{growthofcocycle}
    	Suppose $G$ is a word hyperbolic group. Fix a conjugacy class $\langle h\rangle$ of $G$. Then every element in $HC^n(\mathbb CG, \langle h\rangle )$ has a representative $\varphi:G^{n+1}\to\mathbb C$ such that $\varphi$ is of polynomial growth. Furthermore, when $n\geqslant2$, such representative can be chosen to be uniformly bounded over $G^{n+1}$.
\end{theorem}

\begin{proof}
    
	 Elements of $HC^n(\mathbb CG,\left\langle h\right\rangle )$ have the following description,  cf. \cite[Section 4.1]{Lott1}.
	Let $C^n(G,G_h,h)$ be the space spanned by all $(n+1)$-linear maps on $\mathbb CG$ satisfying the following conditions:
	\begin{align}
	&\phi(g_{\sigma(0)},g_{\sigma(1)},\cdots,g_{\sigma(n)})=(-1)^\sigma\phi(g_0,g_1,\cdots,g_n)
	,\ \forall \sigma\in S_n; \label{chainofGh/h1} \\
	&\phi(zg_0,zg_1,\cdots,zg_n)=\phi(g_0,g_1,\cdots,g_n),\ \forall z\in G_h; \label{chainofGh/h2}\\
	&\phi(hg_0,g_1,\cdots,g_n)=\phi(g_0,g_1,\cdots,g_n). \label{chainofGh/h3}
	\end{align}
	Define a coboundary map $\partial \colon C^n(G,G_h,h)\to C^{n+1}(G,G_h,h)$ by
	$$\partial \phi(g_0,g_1,\cdots,g_{n+1})=\sum_{j=0}^{n+1}(-1)^{j}\phi(g_0,g_1,\cdots,g_{j-1},g_{j+1},\cdots,g_{n+1}.)$$
	Denote the resulting cohomology groups by $H^*(G,G_h,h)$. For each cocycle $\phi$ in $C^n(G,G_h,h)$, there is a cyclic cocycle $\mathcal T_\phi\in C^n(\mathbb CG,\left\langle h\right\rangle )$ given by
	\begin{equation}
	\mathcal T_\phi(g_0,g_1,\cdots,g_n)=
	\begin{cases}
	0& \text{if}\ g_0g_1\cdots g_n\notin \left\langle h\right\rangle,\\
	\phi(\gamma,\gamma g_0,\cdots,\gamma g_0\cdots g_{n-1})& \text{if}\ g_0g_1\cdots g_n=\gamma^{-1}h\gamma.
	\end{cases}
	\end{equation}

	Let us first consider the case of cyclic cocycles with degree $\geq 2$. 
  Observe that if $\phi$ is uniformly bounded, then $\mathcal T_\phi$ is also uniformly bounded. Therefore, it suffices to show that $\phi$ is uniformly bounded. It is not difficult to see that  $H^*(G,G_h,h)$ is isomorphic to $H^*(N_{h};\mathbb C)$ . If the order of $h$ is infinite, then $H^n(N_{h};\mathbb C)$ vanishes for $n>0$,  since $N_{h}$ is a finite group by item $\ref{finiteorder}$ above. Thus in the case, $H^{n}(G,G_h,h)$ with $n>0$ does not contribute to the cyclic cohomology of $\mathbb CG$.
	
Let us assume $h$ has finite order  for the rest of the proof. In fact, it is more convenient for us to work with the group cohomology $H^n(G_{h};\mathbb C)$  of $G_h$. By applying the transfer map, we immediately see that $H^n(G_h;\mathbb C)$ surjects onto $H^n(N_h;\mathbb C)$. More precisely,  consider the chain complex $(E^n(G,G_h),b)$ by removing the condition in line \eqref{chainofGh/h3}. There are two natural chain morphisms: the inclusion map $ \iota \colon (C^n(G,G_h,h),b)\to (E^n(G,G_h),b)$,  and the transfer map $\tau \colon (E^n(G,G_h),b)\to (C^n(G,G_h, h),b)$ defined by
	\begin{equation*}
		\tau(\psi)(g_0,g_1,\cdots,g_n)=\sum_{j=1}^{\textup{ord}(h)}
		\psi(h^{j}g_0,g_1,\cdots,g_n).
	\end{equation*}
	Since $\tau\circ \iota = \textup{ord}(h)\cdot  \textup{Id}$, it follows that $\tau$ induces a surjection on cohomology groups.  Clearly, if $\psi$ is uniformly bounded, then $\tau(\psi)$ is also uniformly bounded.  Therefore it suffices to show that for $n\geq 2$, every element  $H^n(E^\ast(G,G_h))$ admits a uniformly bounded representative.

	Let $Y$ (resp. $Y_h$) be the $\Delta$-complex consisting of all simplices of ordered $(n+1)$-tuple $\{g_0, g_1, \cdots, g_n\}$, where $g_i\in G$ (resp. $g_i\in G_h$). Observe that  $G_h$ acts freely on both $Y$ and $Y_h$. Moreover, we see that the cochain complex of $G_h$-equivariant simplicial cochain on $Y$ is essentially\footnote{To be precise, elements of $E^\ast(G,G_h)$ are assumed to be skew-symmetric (i.e. the condition in line $\eqref{chainofGh/h1}$). But this can be easily fixed by applying a standard anti-symmetrization map.} $E^\ast(G,G_h)$, and the cochain complex of $G_h$-equivariant simplicial cochain on $Y_h$ gives the standard resolution cochain complex for the group cohomology of $G_h$.
	
	Let $\pi\colon Y \to Y_h$ be any $G_h$-equivariant projection, that is, $\pi\circ i = \textup{Id}$ on $Y_h$, where $i\colon Y_h \to Y$ is the inclusion map; such map always exists (see the discussion below for a specific construction of such a map). Then $\pi$ induces a chain map from  the standard resolution cochain complex for the group cohomology of $G_h$ to $E^\ast(G, G_h)$, which is an isomorphism $\pi^\ast\colon H^n(G_h; \mathbb C) \xrightarrow{\, \cong \, } H^n(E^\ast(G, G_h))$ at the level of cohomology. In particular, any uniformly bounded group cocycle of $G_h$ pulls back to a uniformly bounded cocycle of the complex $E^\ast(G, G_h)$.  On the other hand, by the item $\ref{centralizerishyperbolic}$ above,  $G_h$ is hyperbolic. Therefore, by the item  \ref{polcohomology} above, every element of  $H^n(G_h;\mathbb C)$ has a uniformly bounded representative, when $n\geqslant 2$. This finishes the proof for cyclic cocycles of degree  $n\geqslant 2$.

	If $n=0$, $HC^0(\mathbb CG,\langle h\rangle)$ is a one dimensional linear space spanned by $\tr_{\langle h\rangle}$; and $\tr_{\langle h\rangle}$ is clearly uniformly bounded on $G$.

	The only remaining case is when $n=1$.  We divide the proof of this case as follows. First, we shall show that every element of $H^1(E^\ast(G, G_h))$ has a representative of polynomial growth. Second, we shall construct a  $G_h$-equivariant projection $\pi\colon Y \to Y_h$ such that $\pi$ is simplicial and furthermore Lipschitz with respect to the word length metric on $G$ and the corresponding  subspace metric on $G_h$. More precisely, we say $\pi$ is Lipschitz if  there exists a constant $C>0$ such that 
	\[ d(\pi(g_1), \pi(g_2))  \leqslant C d(g_1, g_2),  \]
	for all $g_1, g_2 \in G$, where $d$ is the given word length metric on $G$. 
	Now by the item  $\ref{centralizerishyperbolic}$, the metric on $G_h$ is quasi-isometric to the subspace metric inherited from $G$. It follows that if a degree $1$ group cocycle $\varphi$ of $G_h$ has polynomial growth, then the pullback $\pi^\ast(\varphi) $ of $\varphi$ to the complex $E^\ast(G, G_h)$ also has polynomial growth. In this case, by applying Lemma $\ref{lm:conj}$ above, it is not difficult to see that the corresponding  $\mathcal T_{\pi^\ast(\varphi)}\in C^n(\mathbb CG,\left\langle h\right\rangle )$ has polynomial growth. 
	
	Let us now show that every element of $H^1(E^\ast(G, G_h))$ has a representative of polynomial growth. By definition, a degree $1$ group cocycle $\varphi$ of $G_h$ is a $G_h$-equivariant function  $\varphi\colon G_h\times G_h\to\mathbb C$, such that $\varphi(g_2,g_3)-\varphi(g_1,g_3)+\varphi(g_1,g_2)=0$. In particular, it follows that $\varphi$ is determined by the function $\psi\colon G_h \to \mathbb C$, where $\psi(g)=\varphi(1,g)$. The cocycle condition implies that $\psi(g_1g_2)=\psi(g_1)+\psi(g_2)$, i.e., $\psi$ is a group homomorphism from $G_h$ to $\mathbb C$. As any homomorphism from $G_h$ to $\mathbb C$ factors through the abelianization of $G_h$, it follows that $\psi$ has polynomial growth, so does $\varphi$.

Now let us construct a  $G_h$-equivariant Lipschitz projection $\pi\colon Y \to Y_h$, which will finish the proof of the theorem by the above discussion.	Fix an element, say $\alpha_i\in G$, for each coset of $G_h$ in $G$. Let $\alpha'_i$ be an element of $G_h$ such that $d(\alpha_i, \alpha'_i) = d(\alpha_i, G_h)$. For simplicity,  if $\alpha_i= e\in G$, then we map $e$ to itself. Define a $G_h$-equivariant map $\pi\colon G\to G_h$ by mapping $\beta \cdot \alpha_i \mapsto \beta \cdot \alpha'_i, $ where $\beta\in G_h$. Clearly, $\pi$ extends by linear combination to a $G_h$-equivariant map from $Y$ to $Y_h$, which will still be denoted by $\pi$.

    Let us show that $\pi\colon Y \to Y_h$ is Lipschitz. For two distinct points $g_1,g_2\in G$, choose a geodesic $[\pi(g_1),\pi(g_2)]$ in the Cayley graph of $G$.  By the quasi-convexity of $G_h$ from item \ref{quasiconvex} above, $[\pi(g_1),\pi(g_2)]$ lies in a $K$-neighborhood of $G_h$.
	By assumption, $G$ is hyperbolic. More specifically, let us assume the Cayley graph of $G$ is $\delta$-hyperbolic.  Then the geodesic $[\pi(g_1),\pi(g_2)]$ lies in the $2\delta$-neighborhood of the union of geodesics $[\pi(g_1),g_1]$, $[g_1,g_2]$ and $[g_2,\pi(g_2)]$. Choose $\gamma\in G$ to be a ``midpoint" of $[\pi(g_1),\pi(g_2)]$, that is,  \[ |d(\pi(g_1),\gamma)-d(\pi(g_2),\gamma)|\leqslant 1. \] 
	Then there exists a point $\beta$ on one of the geodesics $[\pi(g_1),g_1]$, $[g_1,g_2]$ or $[g_2,\pi(g_2)]$ such that $d(\gamma,\beta)\leqslant 2\delta$. 
	\begin{enumerate}[label=$(\arabic*)$]
		\item we claim that, if there exists a point $\beta$ on $[\pi(g_1),g_1]$ such that $d(\gamma,\beta)\leqslant 2\delta$, then we have  $d(\beta,\pi(g_1))\leq 2\delta +K$. Indeed, otherwise, we could find an element $h_1\in G_h$ such that $d(g_1, h_1) < d(g_1,  \pi(g_1))$, which contradicts the fact $d(g_1,\pi(g_1)) = d(g_1, G_h)$. Therefore, in this case, we see that  \[ d(\pi(g_1),\pi(g_2))\leqslant 2(4\delta+K)+1. \] Things are similar for $p'$ lies on $[g_2,\pi(g_2)]$. 
		\item Similarly, if there exists $\beta$ on $[\pi(g_1),g_1]$ such that $d(\gamma,\beta)\leqslant 2\delta$, then we also have \[ d(\pi(g_1),\pi(g_2))\leqslant 2(4\delta+K)+1. \]
		\item If there exists $\beta$ on $[g_1,g_2]$ such that $d(\gamma,\beta)\leqslant 2\delta$, then both $d(g_1, G_h)$ and $d(g_2, G_h)$ are $\leqslant d(g_1, g_2)+2\delta+K$. It follows immediately that 
		\[ d(\pi(g_1),\pi(g_2))\leqslant 3 d(g_1, g_2) + 4\delta+2K. \] 
	\end{enumerate}

To summarize, we see that  there exists a constant $C >0$ such that  
\[ d(\pi(g_1),\pi(g_2))\leqslant Cd(g_1, g_2),\] for all $g_1, g_2\in G$.  This finishes the proof. 

\begin{remark}\label{rmk:bpg}
	 By a theorem of Meyer \cite[Theorem 5.2 \& Corollary 5.3]{Meyer}, the same strategy in the above proof can be used to show the following: given a cyclic cohomology class  $[\alpha]$ of a hyperbolic group $G$, 
	 if $\varphi_1$ and $\varphi_2$ are two representatives with polynomial growth of $[\alpha]$, then there exists a cyclic cocycle $\psi$ such that $\varphi_1 - \varphi_2 = b\psi$ and $\psi$ has polynomial growth. Here  $b:C^n(\mathbb CG)\to C^{n+1}(\mathbb CG)$ is the coboundary map of the cyclic cochain complex. 
\end{remark}

\end{proof}

\section{Delocalized Connes-Chern character of secondary invariants}\label{DelocalizedPairing}
In this section, we construct a delocalized Connes-Chern character map for $C^\ast$-algebraic secondary  invariants and prove the second main theorem (Theorem $\ref{thm:main}$) of the paper.  We will only give the details for the odd dimensional case; the even dimensional case is completely similar. 
\begin{theorem}\label{thm:main}
	Let $M$ be a closed manifold whose fundamental group $G$ is  hyperbolic. Suppose  $\left\langle h\right\rangle $  is non-trivial conjugacy class of $G$. Then every element $[\alpha]\in HC^{2k+1-i}(\mathbb CG,\left\langle h\right\rangle ) $ induces a natural map 
	$$\tau_{[\alpha]} \colon K_i(C^*_{L,0}(\widetilde M)^G)\to \mathbb C$$ such that the following are satisfied.
	\begin{enumerate}[label=$(\alph*)$]
		\item $\tau_{[S\alpha]} = \tau_{[\alpha]}$, where $S$ is Connes' periodicity map \[ S\colon HC^{\ast}(\mathbb CG,\left\langle h\right\rangle )\to HC^{\ast+2}(\mathbb CG,\left\langle h\right\rangle ).\]
		\item Suppose $D$ is an elliptic operator on $M$ such that the lift $\widetilde D$ of $D$ to the universal cover $\widetilde M$ of $M$ is invertible. Let  $\varphi$ be a representative of $[\alpha]$ with polynomial growth. Then the delocalized higher eta invariant $\eta_{\varphi}(\widetilde D)$  \textup{(}cf. Definition $\ref{higherdelocalizedeta}$\textup{)}  converges absolutely.  Moreover, we have 
		$$\tau_{[\alpha]}(\rho(\widetilde D))= - \eta_{\varphi}(\widetilde D),  $$
		where $\rho(\widetilde D)$ is the higher rho invariant of $\widetilde D$. 
	\end{enumerate}
\end{theorem}
 In more conceptual terms, the above theorem provides a formula to compute the Connes-Chern character of  elements of $K_i(C^*_{L,0}(\widetilde M)^G)$. Moreover,  the theorem establishes a precise connection between Higson-Roe's $K$-theoretic higher rho invariants and Lott's higher eta invariants.
 
 \begin{remark}
 Note that, in part $(b)$ of the theorem,  $\eta_{\varphi}(\widetilde D)$ converges absolutely  for all invertible $\widetilde D$. In particular, it is not necessary for the spectral gap of $\widetilde D$ to be sufficiently large. 
 \end{remark}
 
  This section is organized as follows. First, we show that each element in $K_i(C^*_{L,0}(\widetilde M)^G)$ has a particular type of nice representatives. Second, we construct an explicit formula for the map  $\tau_{[\alpha]}$  by using such nice representatives, and prove that the formula is well-defined. We shall only give the details for the case of $K_1(C^*_{L,0}(\widetilde M)^G)$; the even case is completely similar.

	\begin{definition}
		Let $B_L (\widetilde M)^G$ to be the subalgebra of $C_L^*(\widetilde M)^G$ consisting of elements $f(t)\in C_L^*(\widetilde M)^G$ such that $f(t)\in B(\mathcal SG)$ for all $t\in[0,\infty)$ and 	
		$f(t)$ is piecewise smooth with respect to $\|\cdot\|_B$.
	\end{definition}
	
	$B_L (\widetilde M)^G$ is a smooth dense subalgebra of $C_L^*(\widetilde M)^G$. Similarly, we define
	$B_{L,0} (\widetilde M)^G$ to be the kernel of the evaluation map
	$$ev\colon B_L (\widetilde M)^G \to B(\mathcal SG),\ f\mapsto f(0).$$
	Note that $B_{L,0} (\widetilde M)^G$ is a smooth dense subalgebra of $C_{L,0}^*(\widetilde M)^G$. In particular, it follows that \[ K_*(B_{L,0} (\widetilde M)^G)\cong K_*(C_{L,0}^*(\widetilde M)^G). \]
	
	\begin{definition}\label{localloop}
		Let $SC^\ast(\widetilde M)^G$ be the suspension of $C^\ast(\widetilde M)^G$, and  $\varphi\in$ be an invertible element in $ SC^\ast(\widetilde M)^G$, that is,  a loop $\varphi\colon S^1 = [0, 1]/\{0, 1\} \to (C^\ast(\widetilde M)^G)^+$ of invertible elements such that $\varphi(1) =1$, where  $(C^\ast(\widetilde M)^G)^+$ is the unitization $C^\ast(\widetilde M)^G$. We say  $\varphi$ is \emph{local} if it is the image of an invertible element $\psi\in SC^\ast_{L}(\widetilde{M})^G$ under the evaluation map $SC^\ast_{L}(\widetilde{M})^G\to SC^\ast(\widetilde{M})^G$. Similarly, an invertible element $\varphi \in SB(\mathcal SG)$ is called local if it is the image of an invertible element of $\psi\in SB_{L} (\widetilde M)^G$ under the evaluation map. 
	\end{definition} 
	\begin{definition}[{\cite[Definition 3.3]{Keswani}}]\label{polycontrolled}
		A path $\zeta \in B_L (\widetilde M)^G$ is said to have \emph{polynomial B-norm control}  if  
		\begin{enumerate}
			\item the propagation of $\zeta(t)$ is finite and goes to zero as $t\to \infty$;
			\item there exists some polynomial $q$ such that  $\|\zeta(t)\|_B\leqslant q\left(\frac{1}{\prop \zeta(t)}\right)$ for sufficiently large $t\gg 0$. Here $\prop \zeta(t)$ stands for the propagation of $\zeta(t)$. 
		\end{enumerate}
	\end{definition}

In the following, we shall prove a sharpened version of \cite[Proposition 3.6]{Xie}. We show that every element of $K_1(B_{L,0}(\widetilde M)^G)$ has nice representatives that satisfy certain regularity properties, in particular, the polynomial control property above.

Let us first prove the following technical lemma. 
\begin{lemma}\label{lemma:f(D)inB}
	Suppose $D$ is a self-adjoint first order elliptic differential operator over $M$ and $\widetilde D$ is the lifting of $D$ to the universal cover $\widetilde M$ of $M$. If $G=\pi_1(M)$ is hyperbolic and $f\in\mathcal A_{\Lambda,N}$ \textup{(}cf. Definition $\ref{ALambdaN}$\textup{)} with $N\geqslant \frac{3}{2}\dim M+5$ and $\Lambda$ sufficiently large, then $f(\widetilde D)\in B(\mathcal SG)$.
\end{lemma}
\begin{proof}
	Fix a symmetric generating set $S$ of $G$. Let $\ell$ be the corresponding word length function of $G$ determined by $S$ and $|S|$ the cardinality of $S$. Suppose $f(\widetilde D)=\sum_{g\in G}A_g g$. By lemma $\ref{functionalcal}$ and Lemma \ref{estimateA}, there exist $C_1,C_2>0$ such that 
	\begin{equation*}
	|A_g|_1\leqslant C_1 \cdot e^{-C_2\Lambda\ell (g)}, 
	\end{equation*}
	for all but finitely many $g\in G$, where  $|A_g|_1 $ stands for the trace norm of $A_g$.  
	Let us denote 
	\[
	A^{(n)}=\sum_{\ell (g)\leqslant n} A_g g.\]
	It suffices to show that $\{A^{(n)}\}$ is a Cauchy sequence under the norm 
	$\|\cdot\|_{B, p}$ (cf. Definition  $\ref{def:Bpnorm}$). 
	Now for any $m<n$, we have 
	\begin{align*}
	\|A^{(n)}-A^{(m)}\|_{RD, p}^2
	=&\sum_{m<\ell (g)\leqslant n} |A_g|_1^2(1+\ell (g))^{2p}\\
	\leqslant &C_1\cdot  \sum_{j=m+1}^n e^{-C_2\Lambda j}(1+j)^{2p}|S|^j,
	\end{align*}
	and 
	\begin{align*}
	&\|\Delta (A^{(n)}-A^{(m)})\|_{uc}\\
	\leqslant&\sum_{m<\ell (g)\leqslant n} |A_g|_1\sum_{\substack{g_1g_2=g\\\ell (g_1)+\ell (g_2)=\ell (g)}}(1+\ell (g_1))^p(1+\ell (g_2))^p\\
	\leqslant&\sum_{m<\ell (g)\leqslant n}|A_g|_1(1+\ell (g))^{2p}\cdot \#\{(g_1,g_2):g_1g_2=g,\ell (g_1)+\ell (g_2)=\ell (g)\}\\
	\leqslant& C_1\cdot C \cdot\sum_{m<\ell (g)\leqslant n} e^{-C_2\Lambda \ell (g)}(1+\ell (g))^{2p+1}\\
	\leqslant& C_1\cdot C \cdot\sum_{j=m+1}^n e^{-C_2\Lambda j}(1+j)^{2p+1}|S|^j,
	\end{align*}
	where we have used the fact that there exists $C>0$ such that 
	$$\#\{(g_1,g_2):g_1g_2=g,\ell(g_1)|+\ell(g_2)=\ell (g)\}\leqslant C\cdot \ell (g),$$
	since $G$ is hyperbolic. It follows that, as long as $\Lambda$ is sufficiently large,  both  $\|A^{(n)}-A^{(m)}\|_{RD,p}^2$ and \mbox{$\|\Delta (A^{(n)}-A^{(m)})\|_{uc}$ } go to zero,  as $m,n \to \infty$. 
\end{proof}

Now let us show that every element of $K_1(B_{L,0}(\widetilde M)^G)$ has nice representatives that satisfy certain regularity properties, in particular, the polynomial control property above. The main motivation for choosing such nice representatives is to justify  the explicit construction of  $\tau_{[\alpha]}\colon K_1(B_{L,0}(\widetilde M)^G) \to \mathbb C$ below (Definition $\ref{etaforhyperbolic}$).   Moreover, we show that  for a given element of $K_1(\mathscr B_{L,0}(\widetilde M)^\Gamma)$, two different such regularized representatives can be connected by a family of representatives of the same kind. This allows us to show that the integral in line  $\eqref{etapair}$ in Definition $\ref{etaforhyperbolic}$ is independent of the choice of such representatives.  

	\begin{proposition}\label{regularrep}
		Every element $[u]\in K_1(B_{L,0}(\widetilde M)^G)$ has a representative $w\in (B_{L,0}(\widetilde M)^G)^+$
		such that 
		$$
		w(t)=\begin{cases}
		u(t)& \text{if  } 0\leqslant t\leqslant 1,\\
		h(t)& \text{if  } 1\leqslant t\leqslant 2,\\
		e^{2\pi i\frac{F(t-1)+1}{2}}& \text{if  } t\geqslant 2,
		\end{cases}$$
		where $h$ is a path of invertible elements connecting $u(1)$ and 
		$\exp(2\pi i\frac{F(t-1)+1}{2})$, and $F$ is a piecewise smooth map $F\colon [1,\infty)\to D^*(\widetilde M)^G$ satisfying
		\begin{enumerate}[label =$(\arabic*)$]
			\item $F(t)^2-1\in B(\mathcal SG)$ and $F(t)^*=F(t)$,
			\item its derivative $F'(t)\in B(\mathcal SG)$,
			\item $\|F(t)\|_{op}\leqslant 1$, where $\|\cdot\|_{op}$ stands for operator norm, 
			\item propagation of $F(t)$ is finite, and goes to zero as $t\to\infty$.
			\item both $F(t)^2-1$ and $F'(t)$ have polynomial B-norm control in the sense of  Definition  $\ref{polycontrolled}$ above. 
		\end{enumerate}
		Moreover, if $v$ is another such representative, then there exists a piecewise smooth path of invertibles $u_s\in (B_{L,0}(\widetilde M)^G)^+$ and piecewise smooth maps $F_s:[0,\infty)\to D^*(\widetilde M)^G$ satisfying conditions above, with $s\in[0,1]$, such that
		\begin{enumerate}[label =\textup{(\Roman*)}]
			\item $u_0=w$, $u_1(t)=v(t)$ for $t\notin (1,2)$,
			\item $u_s(t)=\exp(2\pi i\frac{F_s(t-1)+1}{2})$ for all $t\geqslant 2$,
			\item $u_1v^{-1}:[1,2]\to (B(\mathcal SG))^+$ is a local loop of invertible elements,
			\item $\partial_s(F_s)$ has polynomial $B$-norm control, 
			\item the operator norm of $F_s(t)$ is uniformly bounded, and the degrees of polynomials used for the polynomial $B$-norm control of $F_s$ and $\partial_s F_s$ are uniformly bounded, and the propagation of $F_s(t)$ goes to zero uniformly in $s$, as $t\to \infty$.
		\end{enumerate}
	\end{proposition}
	
	\begin{remark}
	We shall call a representative appearing in the proposition above a \emph{regularized} representative from now on.  
	\end{remark}
	\begin{proof}
	    View the invertible element $u\in (B_{L,0}(\widetilde M)^G)^+$ as an invertible element in $(B_L (\widetilde M)^G)^+$. Consider the element $\hat u=u:[1,\infty)\to (B(\widetilde M)^G)^+$ in $K_1(B_L (\widetilde M)^G)$. Since the $K$-theory of $B_L (\widetilde M)^G$ is the $K$-homology of $M$, it follows from the Baum-Douglas geometric description of $K$-homology \cite{MR679698} that $\hat u$ can be represented by a twisted Dirac operator over a $spin^c$ manifold. More precisely, let $X$ be a $spin^c$ manifold   together with a vector bundle $E$ over $X$ and  a continuous map $\psi\colon X \to M$.  Suppose $D_E$ is the associated twisted Dirac operator on $X$. Let $\widetilde X$ be the $G$-covering space of $X$ induced by $\psi$, and $\widetilde D$ be the lift of $D_E$ to $\widetilde X$. Choose an odd continuous function $\chi\colon \mathbb R \to [-1, 1]$ such that $\chi(x) \to \pm 1$ as $x\to \pm\infty$ and
	    its distributional Fourier transform $\widehat \chi$ has compact support. We define $F(t) = \psi_\ast(\chi(\widetilde D/t))$, where $\psi_\ast\colon D^\ast(\widetilde X)^G \to D^\ast(\widetilde M)^G$ is the natural map induced by $\psi$. It is not difficult to see that $F$ satisfies the properties (1)-(5) listed above\footnote{ Since by construction the propagation of $F(t)$ is uniformly bounded (in particular finite) for all $t\in [1, \infty)$, the polynomial $B$-norm control in property $(5)$ follows from the work of \cite[Section 4]{Keswani}. Roughly speaking, the polynomial $B$-norm control is a consequence of the existence of partition of unity $\{\psi_{n, j}\}$ on the manifold $X$  for each $n\in \mathbb N$ such that  the diameter of $\psi_{n, j}$ is $\leqslant 1/n$ and  the norm of $d\psi_{n,j}$ is bounded by $q(1/n)$ for some polynomial $q$.  }.  Moreover, we have $[\hat u]=[e^{2\pi i\frac{F(t)+1}{2}}]\in K_1(B_L (\widetilde M)^G).$ In particular,  $u$ is homotopic to the invertible element $w$ defined by
	$$w(t)=\begin{cases}
	u(t)& \textup{ if } 0\leqslant t\leqslant 1,\\[5pt]
	h(t)&\textup{ if } 1\leqslant t\leqslant 2, \\[5pt]
	e^{2\pi i\frac{F(t-1)+1}{2}} & \textup{ if } t\geqslant 2,
	\end{cases}$$
	where $h$ is a path of invertible elements connecting $u(1)$ and $e^{2\pi i\frac{F(1)+1}{2}}$. 
	
	Now suppose $v$ is another representative of $[u]$ such that
		$$v(t)=\begin{cases}
	u(t)& \textup{ if } 0\leqslant t\leqslant 1,\\[5pt]
	g(t)&\textup{ if } 1\leqslant t\leqslant 2,\\[5pt]
	e^{2\pi i\frac{G(t-1)+1}{2}} & \textup{ if } t\geqslant 2,
	\end{cases}$$
	where $g$ is a path of invertible elements connecting $u(1)$ and $e^{2\pi i\frac{G(1)+1}{2}}$, and $G_t$ is a piecewise smooth map from $[1,\infty)$ to $D^*(\widetilde M)^G$ satisfying the properties (1)-(5) above. 
	
	By Theorem 3.8 in \cite{Keswani}, there exists a piecewise smooth family $F_s:[1,\infty)\to D^*(\widetilde M)^G$ with $s\in [0,1]$ such that $F_0=F$ and $F_1=G$; $F_s(t)^*=F_s(t)$; propagation of $F_s(t)$ goes to zero, as $t\to \infty$; and all $F_s(t)^2-1$, $\partial_t F_s(t)$, and $\partial_s F_s(t)$ lie in $B(\mathcal SG)$. Furthermore, since the propagation of $\partial_sF_s(t)$ (resp. $F_s(t)$) is finite and the propagation is bounded uniformly in $t$, it is not difficult to see that $\partial_s F_s$ has polynomial $B$-norm control and the degrees of polynomials used for the polynomial $B$-norm control of $F_s$ and $\partial_s F_s$ are uniformly bounded.
	
	Let $\varpi:[0,\infty)\to (B(\mathcal SG))^+$ be the path of invertibles defined as
			$$\varpi(t)=\begin{cases}
	u(t)& \textup{ if } 0\leqslant t\leqslant 1,\\[5pt]
	h(t)& \textup{ if } 1\leqslant t\leqslant 2,\\[5pt]
	e^{2\pi i\frac{F_s(1)+1}{2}} & \textup{ if } 2\leqslant t=s+2\leqslant 3,\\[5pt]
	e^{2\pi i\frac{G(t-2)+1}{2}} &t\geqslant 3.
	\end{cases}$$
	Clearly, $w$ is homotopic to $\varpi$. On the other hand, 
	after a re-parameterization, it is not difficult to see that $\varpi$ differs from $v$ by the loop $f\colon [0,1]\to(B(\mathcal SG))^+$ defined by
	$$f(t)=\begin{cases}
	g(t)^{-1}h(2t)& \textup{ if } 0\leqslant t\leqslant 1/2,\\[5pt]
	g(t)^{-1}e^{2\pi i\frac{F_{2t-1}(1)+1}{2}}& \textup{ if } 1/2\leqslant t\leqslant 1.
	\end{cases}
	$$
	Moreover, $f$ is a local loop in the sense of Definition \ref{localloop}. This finishes the proof. 
	\end{proof}

The following lemma will be useful in the proof of Theorem $\ref{thm:main}$. 

\begin{lemma}\label{exp,polycontrolled}
	Let $u = u_s(t)$ be the family of invertible elements from Proposition 
	$\ref{regularrep}$ above. Then for any delocalized cyclic cocycle $\varphi$ with polynomial growth,
	$$\lim_{t\to+\infty} \widetilde {\varphi\#\tr}(u_s(t)^{-1}\partial_s(u_s(t))\otimes(u_s(t)\otimes u_s(t)^{-1})^{\otimes m})=0$$
	uniformly in $s\in[0,1]$.
\end{lemma}
\begin{proof}
	By definition, 
	\[ u_s(t) = e^{2\pi i \frac{F_s(t-1)+1}{2}}, \]
	where  $F_s$ and $\partial_sF_s$ have polynomial $B$-norm control. Denote $P = P_s(t) = \frac{F_s(t) +1}{2}$. Then $P_s^2 - P_s = (F_s^2 - 1)/4$. 
	Let \[ f_n(x) = \sum_{k=0}^{n}\frac{(2\pi i x)^k}{k!}.  \] Note that we have 
	\begin{align*}
	f_n(P)  & =  \sum_{k=0}^{n}\frac{(2\pi i)^k}{k!}P^k \\
	&= 1 + \Big(\sum_{k=1}^{n}\frac{(2\pi i)^k}{k!}\Big) P +  \Big(\sum_{k=1}^{n} \frac{(2\pi i)^k}{k!}\sum_{j=0}^{k-2}P^j\Big)(P^2-P)
	\end{align*} 
	Define 
	\[ A_n = f_n(P) - \Big(\sum_{k=1}^{n}\frac{(2\pi i)^k}{k!}\Big) P^2. \]
	Clearly, $A_n\in B(\mathcal SG)^+$ for all $n\geqslant  1$, and 
	\[ u - A_n = \Big(\sum_{k=n+1}^{\infty}\frac{(2\pi i)^k}{k!}\Big) (P -P^2) +  \Big(\sum_{k=n+1}^{\infty} \frac{(2\pi i)^k}{k!} \sum_{j=0}^{k-2}P^j\Big)(P^2-P).   \]
	Recall that, by construction,  there exists a polynomial $q$ such that 
	\[  \|P^2(t) - P(t)\|_B \leq q\big(1/\prop F(t)\big).  \]
	Since the operator norm of $P(t)$ is uniformly bounded and the propagation of $P(t)$ goes to zero as $t\to \infty$, a routine calculation shows that   there exists $C>0$ such that
	\[ \| P^j(t) (P^2(t) -P(t))\|_B \leq C^j\cdot  q\big(1/\prop F(t)\big)   \]
	for all $j\in \mathbb N$ and sufficiently large $t\gg 0$.
	It follows that there exists $K>0$ such that 
	\[ \|u(t) - A_n(t)\|_B < K \frac{C^n}{(n+1)!} \cdot  q\big(1/\prop F(t)\big), \]
	and 
	\[ \|u(t) - 1 \|_B < K  e^C \cdot q\big(1/\prop F(t)\big) \]
	for sufficiently large $t\gg 0$.
	The same type of estimates also apply to $u^{-1}$ and $\partial_s u$. 
	
	Now fix  $\varepsilon >0$ so that, if $a_i\in B(\mathcal SG)^+$ has propagation $\leqslant \varepsilon$ for all $0\leqslant i \leqslant 2m$, then 
	\[ \widetilde {\varphi\#\tr}(a_0\otimes a_1\otimes  \cdots \otimes a_{2m}) = 0. \]
	By the proof of Proposition $\ref{convergenceofTRhigher}$, such an $\varepsilon$ exits. 
	For each $n\in \mathbb N$, there exists $t_n>0$ such that $\prop F_s(t) < \frac{\varepsilon}{n}$, 
	for all $t> t_n$. In particular, we have that 
	\begin{enumerate}[label=$(\arabic*)$]
		\item $\prop A_n(t) < \varepsilon$, 
		\item $\|u(t) - A_n(t)\|_B < K \frac{C^n}{(n+1)!} \cdot  q\big(n/\varepsilon)$,
		\item $\|u(t) - 1 \|_B < K  e^C \cdot q\big(n/\varepsilon)$
	\end{enumerate} 
	for all $t> t_n$. Similarly, the same type of estimates hold for $u^{-1}$ and $\partial_s u$. The lemma easily follows from these estimates. This finishes the proof. 
\end{proof}

Now for each class  $[\alpha]\in HC^{2m}(\mathbb CG,\left\langle h\right\rangle )$. we define a map 
\[ \tau_{[\alpha]}\colon K_1(B_{L,0}(\widetilde M)^G) \to \mathbb C\]
as follows.  
	\begin{definition}\label{etaforhyperbolic}
	Let $\varphi$ be a representative of $[\alpha]$ with polynomial growth. For each $[u]\in K_1(B_{L,0}(\widetilde M)^G)$, let $w$ be a regularized representative of $[u]$.  We define
			\begin{equation}\label{etapair}
		\tau_{[\alpha]}([u])\coloneqq \tau_{\varphi}(w) = \frac{m!}{\pi i} \int_{0}^{\infty} \widetilde{\varphi\#\tr}(\dot{w}(t)w(t)^{-1}\otimes (w(t)\otimes w(t)^{-1})^{\otimes m})dt.
		\end{equation}

	\end{definition}

The convergence of the integral in line $\eqref{etapair}$ follows from the following two observations.
 \begin{enumerate}[label=$(\arabic*)$]
 	\item By Proposition \ref{highertrace}, the integrand 
 	\[ \widetilde{\varphi\#\tr}(\dot{w}(t)w(t)^{-1}\otimes (w(t)\otimes w(t)^{-1})^{\otimes m}) \]  is a piecewise smooth function with respect to $t$ on $[0, \infty)$. In particular, this implies that the integral in line $\eqref{etapair}$ converges absolutely for small $t$. 
 	\item 	By the proof of Proposition \ref{regularrep}, when $t\geqslant 2$, we have that
 	$$w(t)=e^{2\pi i\frac{F(t)+1}{2}}=e^{2\pi i \frac{\chi(\widetilde D/t)+1}{2}}$$
 	Set $s = 1/t$, and then we have 
 	\begin{align*}
 	&\int_{2}^{\infty}\widetilde{\varphi\#\tr}(\dot{w}(t)w(t)^{-1}\otimes (w(t)\otimes w(t)^{-1})^{\otimes m})dt\\
 	=& \pi i \int_{1/2}^{0}\widetilde{\varphi\#\tr}\left( \dot\chi(s\widetilde D) \widetilde D
 	\otimes \big(e^{2\pi i\frac{\chi(s\widetilde D)+1}{2}}\otimes e^{-2\pi i\frac{\chi(s\widetilde D)+1}{2}}\big)^{\otimes m}\right)ds.
 	\end{align*}
 	Since the Fourier transform of $\chi$ has compact support,  	
 	it follows that $x\chi'(x)$ and 
 	$\exp(\pm 2\pi i\frac{\chi(x)+1}{2})-1$ lies in $\mathcal A_{\Lambda,N}$ for any $\Lambda,N$. Therefore by Proposition \ref{convergenceofTRhigher} and \ref{convergenceatzerohigher}, the integral with respect to $s$ converges absolutely for small $s$. Consequently the integral from line $\eqref{etapair}$ converges for large $t$.
 \end{enumerate}

 \begin{proof}[Proof of Theorem $\ref{thm:main}$]
 		
 	Let $\varphi$ be a representative of $[\alpha]$ with polynomial growth. 
 	     Let us first show that $\tau_{\varphi}([u])$ is independent of the choice of regularized representative of $[u]$. 
 	    Suppose $w$ and $v$ are two regularized representatives of $[u]$. By Proposition \ref{regularrep}, there exists a piecewise smooth family of invertibles $u_s\in B_{L,0}(\widetilde M)^G$ with the stated properties (I)-(V).

 	    Now a straightforward calculation shows that 
 	    \begin{equation}\label{transgression}
 	    \partial_s\left(\widetilde{\varphi\#\tr}(u^{-1}\partial_t u \otimes(u\otimes u^{-1})^{\otimes m})\right)=
 	    \partial_t\left(\widetilde{\varphi\#\tr}(u^{-1}\partial_s u\otimes(u\otimes u^{-1})^{\otimes m})\right).
 	    \end{equation}	
 	    It follows that 
 	    \begin{align*}
 	    & \int_0^T  \widetilde{\varphi\#\tr}(\dot{u}_1u_1^{-1}\otimes (u_1\otimes u_1^{-1})^{\otimes m})dt - \int_0^T  \widetilde{\varphi\#\tr}(\dot{u}_0u_0^{-1}\otimes (u_0\otimes u_0^{-1})^{\otimes m})dt \\
 	    & =
 	    \int_0^1 \widetilde{\varphi\#\tr}(u^{-1}\partial_s u\otimes(u\otimes u^{-1})^{\otimes m})\Big|_{t=0}^{t=T} ds
 	    \end{align*}
 	    By Lemma \ref{exp,polycontrolled} below, we have 
 	    $$\widetilde{\varphi\#\tr}(u^{-1}\partial_s u\otimes(u\otimes u^{-1})^{\otimes m})\to 0$$
 	    as $t\to \infty$. Also, note that $u_s(0)\equiv 1$ for all $s$. It follows that  
 	    \[  \tau_{\varphi}(u_1) = \tau_{\varphi}(u_0). \]
 	    On the other hand, $u_1$ differs from $v$ by a local loop $f \colon  S^1\to B(\mathcal SG)^+$. 
 	    By \cite[Lemma 3.4]{Xie}, for $\forall \varepsilon > 0$, there exists an idempotent $p\in B(\mathcal SG)^+$ such that the propagation of $p$ is $\leqslant \varepsilon$ and $f$ is homotopic, in the algebra  $SB(\mathcal SG)$,  to the element \[ \beta(t) = e^{2\pi i t} p + (1-p), \textup{ where } 0\leqslant t \leqslant 1. \]
 	    It follows  that 
 	    \begin{align*}
 	    & \int_{0}^{1}  \widetilde{\varphi\#\tr}(\dot{f}(t)f(t)^{-1}\otimes (f(t)\otimes f(t)^{-1})^{\otimes m})dt \\
 	    =& \int_{0}^{1}  \widetilde{\varphi\#\tr}(\dot{\beta}(t)\beta(t)^{-1}\otimes (\beta(t)\otimes \beta(t)^{-1})^{\otimes m})dt \\
 	    =&\int_{0}^{1}  \varphi\#\tr(2\pi ip\otimes ((e^{2\pi i t}-1)p\otimes (e^{-2\pi it}-1)p)^{\otimes m})dt
 	    \end{align*} 
 	    where the last integral is clearly zero,  as long as $\varepsilon$ is sufficiently small. Therefore, $\tau_{[\alpha]}([u])$ is independent of the choice of regularized representative of $[u]$. 
 	    
 	    For a given regularized representative $w$ of $[u]$, the same proof from Proposition \ref{welldefinedofeta}  show that 
 	      $\tau_{[\alpha]}(w)$ is independent of the choice of polynomial growth representative $\varphi$ of $[\alpha]$ (cf. Remark $
 	      \ref{rmk:bpg}$). This proves that the map \[ \tau_{[\alpha]}\colon K_1(B_{L,0}(\widetilde M)^G) \to \mathbb C\] is well-defined.  
 	    Furthermore, the same proof from Proposition \ref{suspension} shows that \[ \tau_{\varphi}(w) = \tau_{S\varphi}(w). \] 
 	   This proves part $(a)$ of the theorem. 
 	    
 	    We shall prove  part $(b)$ of  the theorem in three steps. 
 	    \begin{enumerate}[label=(\roman{*})]
 	    	\item Recall that the definition of $\eta_\varphi(\widetilde D)$ (cf. Definition $\ref{higherdelocalizedeta}$) uses the  representative $u_t = U_{1/t}(\widetilde D)$ of the higher rho invariant $\rho(\widetilde D)$, where 
 	    	\[ U_t=e^{2\pi iF_t(\widetilde D)} \textup{ with }
 	    	F_t(x)=\frac{1}{\sqrt \pi}\int_\infty^{x/t}e^{-s^2}ds \textup{ for } t>0 \textup{ and } U_0 = 1.\]
 	    	We first prove that the path $U_t$ is an element of $(B_{L,0}(\widetilde M)^G)^+$. 
 	    	\item Second, we prove that the integral 
 	    	\[ \eta_\varphi(\widetilde D)= (-1)\frac{m!}{\pi i }\int_0^\infty \widetilde{\varphi\#\tr}(\dot{U}_tU_t^{-1}\otimes(U_t\otimes U_t^{-1})^{\otimes m}) dt  \]
 	    	absolutely converges, where the minus sign is due to the change of variables $1/t\to t$. Note that here we do not require the spectral gap of $\widetilde D$ to be sufficiently large. In other words, the convergence of the integral holds as long as $\widetilde D$ is invertible. 
 	    	\item Recall that we defined  $\tau_{[\alpha]}(\rho(\widetilde D))$ by using a regularized representative of $\rho(\widetilde D)$. In the third step, we use a transgression formula as in line $\eqref{transgression}$ to prove that $\tau_{[\alpha]}(\rho(\widetilde D)) = - \eta_\varphi(\widetilde D)$. 	
 	    \end{enumerate}  
     The first and second steps are proved in Proposition $\ref{convergenceatinf-hypb}$ below. Let us now turn to the third step.  Let $\chi$ be a normalizing function from the proof of Proposition \ref{regularrep}, that is, an odd continuous function $\chi\colon \mathbb R \to [-1, 1]$ such that $\chi(x) \to \pm 1$ as $x\to \pm\infty$ and
     its distributional Fourier transform $\widehat \chi$ has compact support. Furthermore, without loss of generality, we can assume in addition $x\cdot \widehat\chi(x)$ is a smooth function. 
     	Denote $E_t(\widetilde D) = \frac{\chi(\widetilde D)+1}{2}$. It follows from  Lemma \ref{estimateA} and Lemma \ref{lemma:f(D)inB}  that 
     	 	$e^{2\pi iE_t(\widetilde D)}$ is a smooth path in $(B(\mathcal SG))^+$.
     	 Let us define 
     	 	$$V_t=\begin{cases}
     	 	U_t & \textup{ if } 0\leqslant t\leqslant 1,\\[5pt]
     	 	e^{2\pi i ((2-t)F_1(\widetilde D)+(t-1)E_1(\widetilde D))} 
     	 	& \textup{ if } 1\leqslant t\leqslant 2, \\[5pt]
     	 	 e^{2\pi i E_{t-1}(\widetilde D)} & \textup{ if } t\geqslant 2.
     	 	\end{cases}$$ 	    	 
 	    	 Then the path $V_t$ is a regularized representative of $\rho(\widetilde D)$ in $B_{L,0}(\mathcal SG)^+$.  Furthermore, $V_t$ and $U_t$ are homotopic in $B_{L,0}(\mathcal SG)^+$ by the following family of elements $H_s$, with $0\leq s \leq 1$:
 	    	 $$H_s(t)=\begin{cases}
 	    	 	U_t& \textup{ if } 0\leqslant t\leqslant 1,\\[5pt]
 	    	 		e^{2\pi i ((2-t)F_1+(t-1)(sE_1+ (1-s) F_1))} 
 	    	 	& \textup{ if }  1\leqslant t\leqslant 1+s,\\[5pt]
 	    	 e^{2\pi i((1-s)F_{t-1}+sE_{t-1})}
 	    	 	& \textup{ if }  t\geqslant 1+s.
 	    	 \end{cases}$$
 	    Now the same transgression formula in line $\eqref{transgression}$ can be applied to show that 
 	    $$\tau_{\varphi}(V)= \tau_{\varphi}(U) = -\eta_\varphi(\widetilde D).  $$
 	    This finishes the proof.
 \end{proof}

\begin{proposition}\label{convergenceatinf-hypb}
	Under the same assumptions of Theorem $\ref{thm:main}$, let $U_t$ be the representative of $\rho(\widetilde D)$ given by 
	\[ U_t=e^{2\pi iF_t(\widetilde D)} \textup{ with }
	F_t(x)=\frac{1}{\sqrt \pi}\int_\infty^{x/t}e^{-s^2}ds  \textup{ for } t>0 \textup{ and } U_0 = 1.\]  Then the path $U_t$ defines an invertible element of $B_{L,0}(\mathcal SG)^+$. Furthermore, if $\varphi$ in  $C^{2m}(\mathbb CG,\langle h\rangle))$ has polynomial growth,  then  the following integral
	$$\int_0^\infty \widetilde{\varphi\#\tr}(\dot{U}_tU_t^{-1}\otimes(U_t\otimes U_t^{-1})^{\otimes m}) dt$$
	converges absolutely.
\end{proposition}
\begin{proof}
	Since $U_t(x)-1$ admits an analytic continuation to an entire function, it follows from Lemma $\ref{lemma:f(D)inB}$ that  $U_t = U_t(\widetilde D)\in B(\mathcal SG)^+$ for each $t>0$ and the path $U_t$ is smooth with respect to the norm $\|\cdot \|_B$ on $(0, \infty)$. It remains to show that $U_t$ is continuous  at $t=0$ with respect to the norm $\|\cdot \|_B$. 
	
	Since $\widetilde D$ is invertible, let $\sigma > 0 $ be the spectral gap of $\widetilde D$ at zero.
	Then the spectral radius of $e^{-\widetilde D^2}$ as an element in $B(\mathcal SG)$ is $e^{-\sigma^2}$, since $B(\mathcal SG)$ is a smooth dense subalgebra of $C^*(\widetilde M)^G$. Recall that  $B(\mathcal SG)$ is a Banach algebra with respect to the norm\footnote{Rescale the norm $\|\cdot\|_B$ if necessary.} $\|\cdot\|_B$ (cf. Proposition $\ref{estimationondelta}$), that is, 
	$$\|A_1A_2\|_{B}\leqslant \|A_1\|_B\|A_2\|_B\text{ , for any }A_1,A_2\in B(\mathcal SG).$$
	By the spectral radius formula
	$$\lim_{n\to\infty}(\|(e^{-\widetilde D^2})^n\|_B)^{\frac 1 n}=e^{-\sigma^2},$$
	there exists $C_1>0$ such that 
	$$\|e^{-\frac{1}{t}\widetilde D^2}\|_B\leqslant C_1e^{-\frac{1}{2t}\sigma^2}$$
	for all sufficiently small $t>0$.
	It follows that there exists $C >0 $ such that
	\begin{equation}\label{6.3}
	\begin{split}
	\|\dot{U}_tU_t^{-1}\|_B=&
	\big\|-2\sqrt\pi i\frac{\widetilde D}{t^2} e^{-\widetilde D^2/t^2}\big\|_B\\
	\leqslant& \frac{1}{t^2}\big\|2\sqrt\pi \widetilde D e^{-\widetilde D^2}\big\|_B\cdot\big\| e^{-(1/t^2-1)\widetilde D^2}\big\|_B\\
	\leqslant & C\frac{1}{t^2} e^{-\frac{1}{2}\sigma^2/t^2}.
	\end{split}
	\end{equation}
	By the definition of $U_t$, we have 
	$$U_t-1=\exp\left(\int_0^t \dot{U}_sU_s^{-1}ds\right)-1.$$
	In fact, the integral on the right hand side converges in $B(\mathcal SG)$, thanks to the inequality in line \eqref{6.3}. In particular, it follows that 
	\begin{equation}\label{6.4}
	\begin{split}
	\|U_t-1\|_B=&\Big\|\sum_{n=1}^\infty\frac{1}{n!}\left(\int_0^t \dot{U}_sU_s^{-1}ds\right)^n\Big\|_B	\\
	\leqslant& \sum_{n=1}^\infty \frac{1}{n!}\left(\frac{C}{t^2}e^{-\frac{\sigma^2}{2t^2}}\right)^n=	
	\exp\left( \frac{C}{t^2}e^{-\frac{\sigma^2}{2t^2}}\right) -1,
	\end{split}
	\end{equation}
	which goes to zero as $t$ goes to zero. It follows that the path $U_t$, with $t\in [0, \infty)$,  gives an invertible element in $B_{L,0}(\mathcal SG)^+$. Furthermore, by Proposition \ref{convergenceatzerohigher} and Proposition  \ref{highertrace}, the following integral
		\[ \int_0^\infty \widetilde{\varphi\#\tr}(\dot{U}_tU_t^{-1}\otimes(U_t\otimes U_t^{-1})^{\otimes m}) dt\]
		converges absolutely. This finishes the proof. 
\end{proof}

\section{Delocalized higher Atiyah-Patodi-Singer index theorem}\label{Application}
In this section, we apply the results from previous sections to prove a delocalized higher Atiyah-Patodi-Singer index theorem.

Let us first review the Connes-Chern character map in our context. We shall only discuss the even dimensional case; the odd case is similar.  Let $G$ be a discrete group, and $[\alpha]\in HC^{2m}(\mathbb CG)$. The $[\alpha]$-component of the Connes-Chern character of an idempotent $p \in \mathcal SG$ is given by  
\begin{equation}\label{ConnesChern}
\ch_{[\alpha]}(p) = 
\frac{(2m)!}{m!}\varphi\#\tr(p^{\otimes 2m+1}),
\end{equation}
where $\varphi$ is a cyclic cocycle representative of $[\alpha]$. It has been implied that  $\ch_{[\varphi]}(p)$ is independent of the choice of representative of $[\alpha]$. Indeed, for a cyclic coboundary $b\psi$, we have 
$$b\psi\#\tr(p^{\otimes 2m+1})=\psi(p^{\otimes 2m})=0.$$
The last equality follows from the fact $\psi$ is cyclic, which in particular implies that $\psi(p^{\otimes 2m}) = - \psi(p^{\otimes 2m})$.

If $G$ is hyperbolic, then by Proposition \ref{highertrace} the formula in line $\eqref{ConnesChern}$ continues to make sense for idempotents in  $B(\mathcal SG)$, as long as  $\varphi$ has polynomial growth. In fact, in this case, $\ch_{[\alpha]}$ defines a Connes-Chern character map at the level of $K$-theory: 
\[ \ch_{[\alpha]}\colon K_0(B(\mathcal SG)) \to \mathbb C.  \]
Indeed, suppose $[p_0] = [p_1]\in K_0(B(\mathcal SG))$.  Let $p_t$ be a piecewise smooth path of idempotents in $B(\mathcal SG)^+$ connecting $p_0$ and $p_1$. Suppose $\varphi$ has polynomial growth. Then a routine calculation shows that
\[ 
\frac{d}{dt}\varphi\#\tr\big(p_t^{\otimes 2m+1}\big)=(2m+1)(b\varphi\#\tr)\Big((\dot{p_t}p_t-p_t\dot{p_t})\otimes p_t^{\otimes 2m+1}\Big) = 0,
\] 
since $\varphi$ is a cyclic cocycle. It follows immediately that 
\[ \varphi\#\tr(p_0^{\otimes 2m+1})=\varphi\#\tr(p_1^{\otimes 2m+1}). \] Therefore, the map $\ch_{[\alpha]}\colon K_0(B(\mathcal SG)) \to \mathbb C$ is well-defined.

By Theorem $\ref{growthofcocycle}$, every cyclic cohomology class of a hyperbolic group has a representative with polynomial growth. Hence, to summarize the above, we have the following proposition. 

\begin{proposition}
	Suppose that $G$ is a hyperbolic group and  $\left\langle h\right\rangle $ is a conjugacy class of $G$. For each class $[\alpha]\in HC^{2m}(\mathbb C, \langle h \rangle)$,  the Connes-Chern character map 
	\begin{equation}
	\ch_{[\alpha]}\colon K_0(C^*_r(G))\to \mathbb C
	\end{equation}
	given by the formula $\eqref{ConnesChern}$ is well-defined. 
Moreover, we have $\ch_{S[\alpha]} = \ch_{[\alpha]}$, where $S\colon HC^{2m}(\mathbb CG,\left\langle h\right\rangle )\to HC^{2m+2}(\mathbb CG,\left\langle h\right\rangle )$ is Connes' periodicity map. 
\end{proposition}
\begin{proof}
The formula for Connes's periodicity map is given in Definition $\ref{def:periodic}$. A straightforward computation shows that  
	$$(S\varphi\#\tr)(p^{\otimes 2m+3})=\frac{1}{2(2m+1)}\varphi\#\tr(p^{\otimes 2m+1}),$$
	from which the second statement of the proposition  immediately follows.	
\end{proof}

The Connes-Chern character $\ch_{[\alpha]}\colon K_0(C_r^*(G)) \to \mathbb C$ above and the delocalized Connes-Chern character map $\tau_{[\alpha]} \colon K_i(C^*_{L,0}(\widetilde M)^G)\to \mathbb C$ from Theorem $\ref{thm:main}$ are related as follows. 

\begin{proposition}\label{boundary}
	Suppose $G$ is hyperbolic and $\left\langle h\right\rangle $ is a nontrivial conjugacy class of $G$. Given $[\alpha]\in HC^{2m}(\mathbb CG,\langle h\rangle)$,  we have the following commutative diagram:
  $$\begin{CD}
  K_0(C_r^*(G)) @>\ch_{[\alpha]}>> \mathbb C\\
  @V\partial VV  @VV{\times (-2)}V\\
  K_1(C_{L,0}^*(\widetilde M)^G) @>\tau_{[\alpha]}>> \mathbb C
  \end{CD}$$
	where $\partial\colon K_0(C_r^*(G))\to K_1(C_{L,0}^*(\widetilde M)^G)$ is the connecting map in the six-term K-theoretical exact sequence for the short exact sequence:
	$$0\to C^*_{L,0}(\widetilde M)^G\to C_L^*(\widetilde M)^G\to C^*(\widetilde M)^G\to 0$$
\end{proposition}
\begin{proof}
	Each element of  $K_0(C_r^*(G))$ is represented by the formal difference of two idempotents in $B(\mathcal SG)^+$. For notational simplicity, let us carry out the computation for an  idempotent $p$ in $B(\mathcal SG)$. 
	
	 Recall that $\partial [p]$ is defined as follows:  let $\{a_t\}_{t\in [0, \infty)}$ be the following lift of $p$ in $B_L(\widetilde M)^G$:
	\[ a_t = \begin{cases}
	   (1-t) p & \textup{ if } 0\leqslant t \leqslant 1, \\
	   0 & \textup{ if }  t \geqslant 1.
	 \end{cases} \]
	 Then we have 
	\[  \partial p \coloneqq u \textup{ with }   u(t) = e^{2\pi i a_t} \textup{ for } t\in [0, \infty). \]
Note that for $0\leqslant t \leqslant 1$, we have \[ u_t=e^{2\pi i (1-t)p}=1+(e^{2\pi i (1-t)}-1)p. \]
It follows that 
	\begin{align*}
	&\int_0^\infty  \widetilde{\varphi\#\tr}(\dot u_tu_t^{-1}\otimes(u_t\otimes u_t^{-1})^{\otimes m}) dt\\
	=&\int_0^1  \widetilde{\varphi\#\tr}(\dot u_tu_t^{-1}\otimes(u_t\otimes u_t^{-1})^{\otimes m}) dt\\
	= &-\varphi\#\tr(p^{\otimes 2m+1})\int_0^1(2\pi i)(e^{2\pi i(1-t)}-1)^m(e^{-2\pi i(1-t)}-1)^mdt\\=&
	-2\pi i\frac{(2m)!}{(m!)^2}\varphi\#\tr(p^{\otimes 2m+1})=
	-\frac{2\pi i}{m!}\ch_{[\alpha]}([p]).
	\end{align*}
It follows that 
	$$\ch_{[\alpha]}([p])=-\frac{1}{2} \tau_{[\alpha]}(\partial [p]).$$
	This finishes the proof. 
\end{proof}

 Let $W$ be a compact $n$-dimensional spin manifold with boundary $M=\partial W$. Suppose $W$ is equipped with a Riemannian metric which has product structure near $M$ and in additional has positive scalar curvature on $M$. Let $\widetilde W$ be the universal covering of $W$ equipped with the metric lift from $W$. Denote $\pi_1(W)$ by $G$. The associated Dirac operator $\widetilde D_W$ naturally defines a higher index in $K_n(C^*(\widetilde W)^G)$, denoted by $\ind_G(\widetilde D_W)$ as in \cite[Section 3]{Xiepos}. Denote the lift of $M$ with respect to the covering map by $\widetilde M=\partial \widetilde W$. The associated Dirac operator $\widetilde D_M$ naturally defines a higher rho invariant $\rho(\widetilde D_M)$ in $K_{n-1}(C_{L,0}^*(\widetilde M)^G)$. The image of $\rho(\widetilde D_M)$ under the natural homomorphism $K_{n-1}(C_{L,0}^*(\widetilde M)^G)\to K_{n-1}(C_{L,0}^*(\widetilde W)^G)$ will still be denoted by $\rho(\widetilde D_M)$.

We denote by $\partial:K_{n}(C^*(\widetilde W)^G)\to K_{n-1}(C_{L,0}^*(\widetilde W)^G)$ the connecting map in the six-term K-theoretical exact sequence for the short exact sequence:
$$0\to C^*_{L,0}(\widetilde W)^G\to C_L^*(\widetilde W)^G\to C^*(\widetilde W)^G\to 0.$$
By \cite[Theorem 1.14]{Piazzarho} and \cite[Theorem A]{Xiepos}, we have 
\begin{equation}\label{K-indexformula}
\partial (\ind_G(\widetilde D_W)) =\rho(\widetilde D_M) \textup{ in } K_{n-1}(C_{L,0}^*(\widetilde W)^G).
\end{equation}
This together with Proposition $\ref{boundary}$ implies  the following delocalized Atiyah-Patodi-Singer index theorem. 
\begin{theorem}\label{application}
	Let $W$ be a compact even-dimensional spin manifold with boundary $M$. Suppose $W$ is equipped with a Riemannian metric which has product structure near $M$ and in additional has positive scalar curvature on $M$. Suppose $G=\pi_1(W)$ is hyperbolic and $\left\langle h\right\rangle $ is a non-trivial conjugacy class of $G$.  Then for any $[\alpha]\in HC^{2m}(\mathbb CG,\left\langle h\right\rangle )$, we have  
	\begin{equation}
	\ch_{[\alpha]}(\ind_G(\widetilde D_W))=\frac{1}{2}\eta_{[\alpha]}(\widetilde D_M).
	\end{equation}
\end{theorem}
\begin{proof}
	Observe that Proposition $\ref{boundary}$ still holds if we replace $\widetilde M$ by $\widetilde W$. In particular, we have the following commutative diagram: 
	 $$\begin{CD}
	K_0(C_r^*(G)) @>\ch_{[\alpha]}>> \mathbb C\\
	@V\partial VV  @VV{\times (-2)}V\\
	K_1(C_{L,0}^*(\widetilde W)^G) @>\tau_{[\alpha]}>> \mathbb C
	\end{CD}$$
	Now the theorem follows immediately from Theorem $\ref{thm:main}$ and  the equality: $$\partial (\ind_G(\widetilde D_W))=\rho(\widetilde D_M)\textup{ in } K_{n-1}(C_{L,0}^*(\widetilde W)^G). $$
\end{proof}

By using Theorem $\ref{thm:main}$, we have derived Theorem $\ref{application}$ as a consequence of a  $K$-theoretic counterpart. This  is possible only because we have realized  $\eta_{[\alpha]}(\widetilde D_M)$ as the pairing between the cyclic cocycle $[\alpha]$ and the $C^\ast$-algebraic secondary invariant $\rho(\widetilde D_M)$ in $K_1(C_{L,0}^\ast(\widetilde W)^G)$.

Alternatively, one can also derive Theorem  $\ref{application}$ from a version of  higher Atiyah-Patodi-Singer index theorem due to Leichtnam and Piazza \cite[Theorem 4.1]{MR1708639} and Wahl \cite[Theorem 9.4 \& 11.1]{Charlotte}. This version of  higher Atiyah-Patodi-Singer index theorem is stated in terms of noncommutative differential forms on a smooth dense subalgebra of $C_r^\ast(G)$; or noncommutative differential forms on a certain class of smooth dense subalgebras (if exist) of general $C^\ast$-algebras (not just group $C^\ast$-algebras) in Wahl's version.  In the case of Gromov's hyperbolic groups, one can choose such a smooth dense subalgebra to be Puschnigg's smooth dense subalgebra $B(\mathbb CG)$. For hyperbolic groups, every cyclic cohomology class of $\mathbb CG$ continuously extends to a  cyclic cohomology class of $B(\mathbb CG)$, cf. Section $\ref{SmoothDenseSubalgebras}$ and Section $ \ref{CyclicCohomologyGrowth}$. Now Theorem $\ref{application}$ follows by pairing the higher Atiyah-Patodi-Singer index formula of Leichtnam-Piazza and Wahl with the delocalized cyclic cocycles of $\mathbb C\Gamma$. 

One can also try to pair the higher Atiyah-Patodi-Singer index formula of Leichtnam-Piazza and Wahl with group cocycles of $\Gamma$, or equivalently cyclic cocycles in  $HC^{\ast}(\mathbb C\Gamma,\left\langle 1\right\rangle )$, where $\langle 1\rangle$ stands for the conjugacy class of the identity element of $\Gamma$. In this case, for fundamental groups with property RD,  Gorokhovsky, Moriyoshi and Piazza proved a higher Atiyah-Patodi-Singer index theorem for group cocycles with polynomial growth \cite[Theorem 7.2]{MR3500822}.

\section{Delocalized higher eta invariant and its relation to Lott's  higher eta invariant}\label{Identification}

 In this section,  we shall establish the relation between our definition of the delocalized higher eta invariant (cf. Definition $\ref{higherdelocalizedeta}$) and Lott's  higher eta invariant \cite[Section 4.4 \& 4.6]{Lott1}. In particular, we prove that our definition of the delocalized higher eta invariant is equal to  Lott's  higher eta invariant up to a constant $\displaystyle \frac{1}{\sqrt{\pi}}$. The main techniques used in this section are from Connes' papers \cite{AC88,MR1116414}.

Let $M$ be a closed manifold and $\widetilde M$ the universal covering over $M$. Suppose $D$ is a first order self-adjoint elliptic differential operator acting on a vector bundle $E$ over $M$ and $\widetilde D$ the lift of $D$ to $\widetilde M$.
 Suppose that $\left\langle h\right\rangle$ is a nontrivial conjugacy class of $G = \pi_1(M)$ and $\widetilde D$ is invertible. Throughout this section,  we assume that  $G$ has polynomial growth.

Let  $\mathscr B$ be the following dense subalgebra of $C_r^*(G)$: 
\[ \mathscr B=\big\{f: G\to \mathbb C  \mid \sum_{g\in G}(1+\ell(g))^{2k}|f(g)|^2<\infty \textup{ for all } k \in \mathbb N \}, \]
where $\ell$ is a word-length function on $G$. $\mathscr B$ is Fr\'echet locally $m$-convex algebra. Moreover, since $G$ has polynomial growth, $\mathscr B$ is a smooth dense subalgebra of $C_r^\ast(G)$.  The universal graded differential algebra of $\mathscr B$ is 
\[ \Omega_\ast(\mathscr B) = \bigoplus_{k=0}^\infty \Omega_k(\mathscr B) \]
where as a vector space, $\Omega_k(\mathscr B) = \mathscr B \otimes (\mathscr B)^{\otimes k}$. As $\mathscr B$ is a Fr\'echet algebra, we consider the completion of $\Omega_\ast(\mathscr B)$, which will still be denoted by $\Omega_\ast(\mathscr B)$. 

Let $\mathfrak E = (\widetilde M\times_G \mathscr B)\otimes E$ be the associated $\mathscr B$-vector bundle and $C^\infty (M; \mathfrak E)$ be its space of smooth sections. Now suppose $\psi$ is a smooth function on $\widetilde M$ with comapct support  such that 
\[ \sum_{g\in G}g^*\psi=1. \] Then we have a superconnection  $\nabla\colon C^\infty(M; \mathfrak E)\to C^\infty(M; \mathfrak E\otimes_{\mathscr B} \Omega_1(\mathscr B))$ given by 
$$\nabla(f)=\sum_{g\in G}(\psi \cdot  g^*f)\otimes_{\mathscr B} dg.$$  
\begin{definition}[{\cite[Section 4.4 \& 4.6]{Lott1}}]\label{def:Lotteta} For each $\beta>0$, 
	Lott's higher eta invariant $\widetilde \eta(\widetilde D)$ is defined by the formula
	\[ \widetilde  \eta(\widetilde D, \beta)= \beta^{1/2}\int_0^\infty \STR(\widetilde De^{-\beta(t\widetilde D +\nabla)^2})dt.  \]
	Here we follow the superconnection formalism, and $\STR$ is the corresponding supertrace, cf. \cite[Proposition 22]{Lott1}. 
\end{definition}

We recall the following periodic version of Lott's higher eta invariant. 
 
 \begin{definition}[{\cite[Section 4.6]{Lott1}}]
 Define $\tilde \eta(\widetilde D) \in\Omega_*(\mathscr B)$ to be  
 	$$\tilde\eta(\widetilde D)=\int_0^\infty e^{-\beta}\tilde\eta(\widetilde D, \beta^2)d\beta.$$
 \end{definition}
Similar estimates as those in Section $\ref{HigherEtaInvariants}$ show that, under the assumption $G$ has polynomial growth, the above integral converges in $\Omega_\ast(\mathscr B)$, hence $\tilde \eta(\widetilde D)$ is well-defined.

Let us write 
 $$\tilde\eta(\widetilde D)=\sum_{m\geqslant 0}\tilde\eta_k(\widetilde D)=\sum_{m\geqslant 0}\int_0^\infty e^{-\beta}\tilde\eta_{2m}(\widetilde D, \beta^2)d\beta,$$
 where $\tilde\eta_{2m}(\widetilde D)$ and $\tilde\eta_{2m}(\widetilde D, \beta^2)$ are the $2m$-th components of $\tilde\eta (\widetilde D)$  and $\tilde\eta(\widetilde D, \beta^2)$  in $\Omega_{2m}(\mathscr B)$ respectively. 
 
  For each $m\geq 0$, only a finite number of terms in Duhamel expansion of $\widetilde De^{-\beta^2(t\widetilde D +\nabla)^2}$ will contribute to $\tilde \eta_{2m}(\widetilde D, \beta^2)$.  Suppose $\varphi\in C^{2m}(G, \langle h\rangle)$ is a cyclic cocycle with  polynomial growth. Without loss of generality, we assume that $\varphi$ is normalized, that is, 
  \[ \varphi(g_0, g_1, \cdots, g_{2k}) =0 \textup{ if } g_i=1 \textup{ for some } i\geq 1.  \]  Let us consider the paring  $\langle\varphi,  \tilde\eta_{2m}(\widetilde D, \beta^2) \rangle.$ Observe that,  since  $\tilde\eta_{2m}(\widetilde D, \beta^2)$  is paired with $\varphi$, we can relax the smoothness condition on $\psi$  in the definition of the connection $\nabla$ above and  choose  $\psi$ to be the characteristic function of a fundamental domain of $\widetilde M$ under the action $G$. More precisely, for such a choice of $\psi$, we should treat the summands $t\nabla \widetilde D$ and $t\widetilde D\nabla$ in the supercommutator  $[\nabla, t\widetilde D] = t(\nabla \widetilde D + \widetilde D\nabla)$ separately so that we avoid taking the differential of $\psi$.   As a consequence, the term $\nabla^2$ does not contribute to the pairing $\langle\varphi,  \tilde\eta_{2m}(\widetilde D, \beta^2) \rangle$, since $\varphi$ is normalized. To summarize, we have 
\begin{align*}
& \langle\varphi,\,  \tilde\eta_{2m}(\widetilde D, \beta^2) \rangle \\
=&\int_0^\infty \beta^{2m+1} \int_{(\sum_{j=0}^{2m} s_j)=\beta}\Big\langle \varphi, \STR\big(
  \widetilde De^{-s_0\beta t^2\widetilde D^2}[\nabla,t\widetilde D]
  e^{-s_1\beta t^2\widetilde D^2} [\nabla,t\widetilde D]\cdots \\
  & \hspace{7cm} \times [\nabla,t\widetilde D]
  e^{-s_{2m}\beta t^2\widetilde D^2}\big) \Big\rangle ds_0ds_1\cdots ds_{2m}dt\\[2ex]
=&\int_0^\infty \beta^{m+\frac{1}{2}}  \int_{(\sum_{j=0}^{2m} s_j)=\beta }\Big\langle \varphi, \STR\big(
\widetilde De^{-s_0t^2\widetilde D^2}[\nabla,t\widetilde D]
e^{-s_1t^2\widetilde D^2} [\nabla,t\widetilde D]\cdots \\
& \hspace{7cm} \times [\nabla,t\widetilde D]
e^{-s_{2m}t^2\widetilde D^2}\big) \Big\rangle ds_0ds_1\cdots ds_{2m}dt,\\
  \end{align*}
where the second equality follows from the change of variables $t \mapsto \sqrt{\beta} t$. Here $[\nabla, t\widetilde D]$ is the supercommutator (i.e. graded-commutator) of $\nabla$ and $t\widetilde D$.

 Let $f_0(x)=\widetilde De^{-xt^2\widetilde D^2}$ and 
 $f_j(x)=[\nabla,t\widetilde D]
 e^{-xt^2\widetilde D^2}$ for $j>0$. From the above calculation, we see that
$$\langle \varphi, \tilde\eta_{2m}(\widetilde D, \beta^2) \rangle = \beta^{m+\frac{1}{2}}  \int_0^\infty \big\langle \varphi, \STR\big(f_0 \ast f_1\ast\cdots \ast f_{2m}(\beta)\big)\big\rangle dt,$$
 where $\ast$ stands for the convolution:
 \[ (f\ast h)(\beta) = \int_{0}^\beta f(x)h(\beta-x) dx.  \]
Recall that  the Laplace transform 
\[ f\mapsto \mathcal L(f)(s) = \int_0^\infty e^{-s\beta}f(\beta) d\beta \] 
converts convolutions of functions into pointwise products of functions. Let us define 
\[ H(s) \coloneqq  \big\langle \varphi, \STR\textstyle{\big(\prod_{j=0}^{2m}\mathcal L(f_j)(s)\big)} \big\rangle,  \]
which is the Laplace transform of 
\[  \int_0^\infty \big\langle \varphi, \STR\big(f_0 \ast f_1\ast\cdots \ast f_{2m}(\beta)\big)\big\rangle dt. \]
Recall that  
 $$\mathcal L(e^{-a\beta})(s)=\frac{1}{a+s}.$$
Therefore, we have 
 \begin{align*}
  H(s) =& \int_0^\infty \big\langle \varphi, \STR{\textstyle \big(
 \frac{\widetilde D}{t^2\widetilde D^2+s}[\nabla,t\widetilde D]\frac{1}{t^2\widetilde D^2+s}\cdots [\nabla,t\widetilde D]\frac{1}{t^2\widetilde D^2+s}\big) }\big\rangle dt \\
 = &  s^{-m-\frac{1}{2}}\int_0^\infty \big\langle \varphi, \STR{\textstyle \big(
 \frac{\widetilde D}{t^2\widetilde D^2+1}[\nabla,t\widetilde D]\frac{1}{t^2\widetilde D^2+1}\cdots [\nabla,t\widetilde D]\frac{1}{t^2\widetilde D^2+1}\big)} \big\rangle dt \\
= & s^{-m-\frac{1}{2}} H(1)
 \end{align*}
 where the second equality follows from the change of variables $t\to t/\sqrt s$.
Apply the inverse Laplace transform to $H(s)$ and we obtain
\begin{align*}
&\int_0^\infty \big\langle \varphi, \STR\big(f_0 \ast f_1\ast\cdots \ast f_{2m}(\beta)\big)\big\rangle dt = \frac{\beta^{m-\frac{1}{2}}}{\Gamma(m+\frac{1}{2})} H(1)
\end{align*}
where $\Gamma(m+\frac{1}{2}) = \frac{1}{4^m}\frac{(2m)!}{ m!}  \sqrt \pi$. It follows that 
 \begin{align*}
 &\langle \varphi, \tilde\eta_{2m}(\widetilde D, \beta^2) \rangle =  \frac{\beta^{2m}}{\Gamma(m+\frac{1}{2})} H(1).
 \end{align*}
Hence we have 
\begin{align*}
\langle \varphi, \tilde\eta_{2m}(\widetilde D) \rangle =& \int_{0}^\infty \langle \varphi,  e^{-\beta}\tilde\eta_{2m}(\widetilde D, \beta^2) \rangle d\beta \\
= & 
\frac{4^m m!}{\sqrt \pi}
\int_0^\infty \big\langle \varphi, \STR{\textstyle \big(
	\frac{\widetilde D}{t^2\widetilde D^2+1}[\nabla,t\widetilde D]\frac{1}{t^2\widetilde D^2+1}\cdots [\nabla,t\widetilde D]\frac{1}{t^2\widetilde D^2+1}\big)} \big\rangle dt
\end{align*}
We shall identify this formula with our formula for delocalized higher eta invariant in  Definition $\ref{higherdelocalizedeta}$.  To this end, let us first define 
 $$w_t(x)\coloneqq \frac{tx-i}{tx+i}.$$
 Note that the path $w_{1/t}(\widetilde D)$ is a representative of the higher rho invariant $\rho(\widetilde D)$. 
A direct computation shows that 
 $$\dot w_t(\widetilde D)w_t(\widetilde D)^{-1}=\frac{2i\widetilde{D}}{t^2\widetilde D^2+1}$$
 and
 $$[w_t(\widetilde D),\nabla]\cdot[w_t(\widetilde D)^{-1},\nabla]
 =4(t\widetilde D+i)^{-1}[\nabla,t\widetilde D](t^2\widetilde D^2+1)^{-1}[\nabla,t\widetilde D](t\widetilde D-i)^{-1},$$
where $[w_t(\widetilde D),\nabla]$ and $[w_t(\widetilde D)^{-1},\nabla]$ are the usual ungraded commutator\footnote{In the superconnection formalism of Definition $\ref{def:Lotteta}$ above, the relevant bundle has been ``doubled". In other words, the superconnection is formulated on a superbundle, which is the direct sum of two copies of the original bundle together with the obvious $\mathbb Z/2$-grading. On the other hand,  the operator $w_t(\widetilde D)$ is  defined on the original bundle, instead of the superbundle.}. 
For notational simplicity, let us write $w_t$ in place of $w_t(\widetilde D)$. The above computation implies that   \begin{align*}
 &\tr\big(\dot w_tw_t^{-1} ([w_t,\nabla][w_t^{-1},\nabla])^m\big)\\
 =&4^m\tr\textstyle \big(\frac{2i\widetilde D}{t^2\widetilde D^2+1}\big(\frac{1}{t\widetilde D+i}[\nabla,t\widetilde D]\frac{1}{t^2\widetilde D^2+1}[\nabla,t\widetilde D]\frac{1}{t\widetilde D-i} \big)^m \big)\\
 =&4^m\tr\textstyle \big(\frac{2i\widetilde D}{t^2\widetilde D^2+1}
 [\nabla,t\widetilde D]\frac{1}{t^2\widetilde D^2+1}\cdots 
 [\nabla,t\widetilde D]\frac{1}{t^2\widetilde D^2+1}\big).
 \end{align*}
 Since $\varphi$ is normalized, again $\nabla^2$ does not contribute to the pairing 
 \[ \langle\varphi, \tr\big(\dot w_tw_t^{-1} ([w_t,\nabla][w_t^{-1},\nabla])^m\big)\rangle.  \]
It follows that 
 \begin{align*}
 &\langle\varphi, \tr\big(\dot w_tw_t^{-1} ([w_t,\nabla][w_t^{-1},\nabla])^m\big)\rangle\\
 =&\langle \varphi, \tr\big(\dot w_tw_t^{-1}(\nabla w_t\nabla w_t^{-1})^m\big)\rangle +\langle \varphi, \tr\big(\dot w_tw_t^{-1}( w_t\nabla w_t^{-1}\nabla)^m\big)\rangle,
  \end{align*}
 when $m\geqslant 1$. By the definition of the connection $\nabla$ and the definition of the trace, we have
\begin{align*}
\langle \varphi, \tr\big(\dot w_tw_t^{-1}(\nabla w_t\nabla w_t^{-1})^m\big)\rangle = \widetilde{\varphi\#\tr}\big(\dot w_tw_t^{-1}\otimes
 (w_t\otimes w_t^{-1})^{\otimes m}\big).
\end{align*}
On the other hand, we have 
\begin{align*}
&\langle \varphi, \tr\big(\dot w_tw_t^{-1}( w_t\nabla w_t^{-1}\nabla)^m\big)\rangle = \langle \varphi,\tr\big(\dot w_t(\nabla w_t^{-1}\nabla w_t)^{m-1} \nabla w_t^{-1} \nabla \big)\rangle = 0.
\end{align*}
Therefore, for all  $m\geqslant 0$, we have
\begin{align*}
&\int_0^\infty \langle \varphi, \tr\big(\dot w_tw_t^{-1} ([w_t,\nabla][w_t^{-1},\nabla])^m\big) \rangle dt \\
= & \int_0^\infty  \widetilde{\varphi\#\tr}\big(\dot w_tw_t^{-1}\otimes 
(w_t\otimes w_t^{-1})^{\otimes m}\big)dt\\
= &   \frac{\pi i}{m!} \tau_{\varphi}(w). 
\end{align*}
where $\tau_\varphi$ is the map from Definition \ref{etaforhyperbolic}. Now  Lemma $\ref{lm:conest}$ below proves the convergence of $\tau_{\varphi}(w)$. To summarize, we have established the following precise relation between our definition of the delocalized higher eta invariant (cf. Definition $\ref{higherdelocalizedeta}$) and Lott's  higher eta invariant \cite[Section 4.4 \& 4.6]{Lott1}. 
\begin{proposition}\label{8.3}
Suppose $D$ is a first order self-adjoint elliptic differential operator acting on a vector bundle $E$ over $M$ and $\widetilde D$ the lift of $D$ to $\widetilde M$.
Assume that $\left\langle h\right\rangle$ is a nontrivial conjugacy class of $G = \pi_1(M)$, $\widetilde D$ is invertible and $G=\pi_1(M)$ has polynomial growth. Then we have  	\[ \tau_{\varphi}(\rho(\widetilde D))=    \frac{1}{\sqrt\pi} \langle \varphi, \tilde\eta_{2m}(\widetilde D)\rangle. \]
\end{proposition}

In the remaining part of this section, we prove Lemma $\ref{lm:conest}$ below and hence complete the proof of the above proposition. The proof of Lemma $\ref{lm:conest}$ uses the technical assumption that $G=\pi_1(M)$ has polynomial growth. It remains an open question how to identify our formulation of higher eta invariants with Lott's higher eta invariant in general. 
\begin{lemma}\label{lm:conest}
	With the same notation as above, the following integral 
	\begin{equation}\label{eq:poly}
	\int_0^\infty  \widetilde{\varphi\#\tr}\big(\dot w_tw_t^{-1}\otimes 
	(w_t\otimes w_t^{-1})^{\otimes m}\big)dt
	\end{equation} 
	converges absolutely. 
\end{lemma}
\begin{proof}
let $a_t(x) = w_t(x) -1 = -2i(tx+i)^{-1}$. For each $t >0$, the Fourier transform of $a_t$ is 
 $$\widehat{a}_t(\xi)=-4\pi\frac{1}{t} e^{-\xi/t}\theta(\xi),$$
 where $\theta$ is the characteristic function of the interval $(0,\infty)$. The function $\widehat{a}_t$ and all its derivatives  are smooth away from $\xi=0$ and decay exponentially as $|\xi|\to\infty$. It follows from the proof of Lemma \ref{functionalcal} that the Schwartz kernel of  $w_t(\widetilde D)-1$ is smooth away from the diagonal of $\widetilde M\times\widetilde M$.  The same holds for $w_t^{-1}-1$ and $\dot w_tw_t^{-1}$. Similar arguments as in the proof of Proposition \ref{convergenceofTRhigher} and Proposition \ref{convergenceatzerohigher} show that  $\widetilde{\varphi\#\tr}(\dot w_tw_t^{-1}\otimes
 (w_t\otimes w_t^{-1})^{\otimes m})$ is finite for each $t>0$ and furthermore the integral in line $\eqref{eq:poly}$ converges absolutely for small $t$.  
 
 Now we prove  the integral in line $\eqref{eq:poly}$ converges absolutely for large $t$. Since $G$ acts freely and cocompactly on $\widetilde M$, there exists a constant $\varepsilon>0$ such that $\dist(x,gx)>\varepsilon$ for all $x\in \widetilde M$ and all $g \neq e\in G$. 
Fix a point $x_0\in \widetilde M$. For $x\in \widetilde M$, let $\nu(x)$ be a smooth approximation of the distance from $x$ to $x_0$. More precisely, let  $\nu$ be a smooth function on $\widetilde M$ satisfying the following:
 \begin{enumerate}
 	\item  $\dist(x,x_0)\leqslant \nu(x)\leqslant 2\dist(x,x_0)$
 	if $\dist(x,x_0)\geqslant \varepsilon$,
 	\item  $\nu$ has uniformly bounded derivatives up to order $N$ with $N$ sufficiently large. 
 \end{enumerate}  
Let $\delta$ be the unbound derivation on $C^*(\widetilde M)^G$ defined by  \[ \delta(T)\coloneqq [T, \nu] = T\circ \nu - \nu \circ T\] for $T\in C^*(\widetilde M)^G$. If $T$ admits a distributional Schwartz kernel which is smooth away from the diagonal of $\widetilde  M\times \widetilde M$, then we have 
 $$\delta^k(T)(x,y)=T(x,y)(\nu(x)-\nu(y))^k,$$
 for all $(x,y)\in\widetilde M\times\widetilde M$.
 
 Denote $A_t= w_t(\widetilde D)-1$. We have
 $$\delta(A_t)=2it\frac{1}{t\widetilde D+i}[\widetilde D,\nu]\frac{1}{t\widetilde D+i}.$$ 
Since $\nu$ has uniformly bounded derivatives and $\widetilde D$ is invertible, there exists $C_1>0$ such that 
 $$\|\widetilde D\delta(A_t)\widetilde D\|_{op}\leqslant \frac{C_1}{t},$$
 where $\|\cdot\|_{op}$ denotes the operator norm.
 By induction, we see that there exists $C_2 > 0 $ such that
 $$\sup_{k+j\leqslant N+1}\|\widetilde D^k\delta^N(A_t)\widetilde D^{j}\|_{op}\leqslant \frac{C_2}{t}.$$
Let $K_t(x,y)$ be the distributional Schwartz kernel of $A_t$. Then the Schwartz kernel of $\delta^N(A_t)$ is 
\[ K_t(x,y)(\nu(y)-\nu(x))^{N}\]
for all $x, y \in \widetilde M$.  It follows from Lemma \ref{kernelupperbound} that  there exists $C_2>0$ such that 
 $$|K_t(x,y)|\cdot |\nu(y)-\nu(x)|^{N}\leqslant \frac{C_2}{t}$$
 for all $x, y\in\widetilde M$, where $|K_t(x,y)|$ is the norm of the matrix $K_t(x, y)$. 
 In particular, if $(x,y)=(x_0,gx_0)$ with $g\neq e$, then we have  
 $$|K_t(x_0,gx_0)|\leqslant \frac{C_2}{t\cdot(\dist(x_0,gx_0))^N}.$$ 
 Now for each $x\in \widetilde M$, use a smooth approximation of the distance function centered at $x$ and apply the same estimates above. Since the action of $G$ on $\widetilde M$ is cocompact, we may choose $C_2$ so that 
  $$|K_t(x,gx)|\leqslant \frac{C_2}{t\cdot(\dist(x,gx))^N}$$
  for all $x\in \mathcal F$ and all $g\neq e\in G$, where $\mathcal F$ is a fundamental domain of $\widetilde M$ under the action of $G$. Similar estimates hold for the Schwartz kernels of $w_t^{-1}-1$ and $\dot w_tw_t^{-1}$.

By assumption,  $G$ has polynomial growth and $\varphi$ is a delocalized cyclic cocycle with polynomial growth. A straightforward computation shows that there exists a constant $C >0$ such that 
$$\Big|\widetilde{\varphi\#\tr}(\dot w_tw_t^{-1}\otimes
(w_t\otimes w_t^{-1})^{\otimes m})\Big|\leqslant \frac{C}{t^{2m+2}},$$
which implies that the integral in line $\eqref{eq:poly}$ converges absolutely for large $t$. This finishes the proof.

\end{proof}


\begin{thebibliography}{10}
	
	\bibitem{A-P-S75a}
	M.~F. Atiyah, V.~K. Patodi, and I.~M. Singer.
	\newblock Spectral asymmetry and {R}iemannian geometry. {I}.
	\newblock {\em Math. Proc. Cambridge Philos. Soc.}, 77:43--69, 1975.
	
	\bibitem{A-P-S75b}
	M.~F. Atiyah, V.~K. Patodi, and I.~M. Singer.
	\newblock Spectral asymmetry and {R}iemannian geometry. {II}.
	\newblock {\em Math. Proc. Cambridge Philos. Soc.}, 78(3):405--432, 1975.
	
	\bibitem{A-P-S76}
	M.~F. Atiyah, V.~K. Patodi, and I.~M. Singer.
	\newblock Spectral asymmetry and {R}iemannian geometry. {III}.
	\newblock {\em Math. Proc. Cambridge Philos. Soc.}, 79(1):71--99, 1976.
	
	\bibitem{MR679698}
	Paul Baum and Ronald~G. Douglas.
	\newblock {$K$}\ homology and index theory.
	\newblock In {\em Operator algebras and applications, {P}art {I} ({K}ingston,
		{O}nt., 1980)}, volume~38 of {\em Proc. Sympos. Pure Math.}, pages 117--173.
	Amer. Math. Soc., Providence, R.I., 1982.
	
	\bibitem{MR3296587}
	Moulay-Tahar Benameur and Indrava Roy.
	\newblock The {H}igson-{R}oe exact sequence and {$\ell^2$} eta invariants.
	\newblock {\em J. Funct. Anal.}, 268(4):974--1031, 2015.
	
	\bibitem{MR0444836}
	P.~Bernat, N.~Conze, M.~Duflo, M.~L\'{e}vy-Nahas, M.~Ra\"{\i}s, P.~Renouard,
	and M.~Vergne.
	\newblock {\em Repr\'{e}sentations des groupes de {L}ie r\'{e}solubles}.
	\newblock Dunod, Paris, 1972.
	\newblock Monographies de la Soci\'{e}t\'{e} Math\'{e}matique de France, No. 4.
	
	\bibitem{JBJC89}
	Jean-Michel Bismut and Jeff Cheeger.
	\newblock {$\eta$}-invariants and their adiabatic limits.
	\newblock {\em J. Amer. Math. Soc.}, 2(1):33--70, 1989.
	
	\bibitem{MR1744486}
	Martin~R. Bridson and Andr\'{e} Haefliger.
	\newblock {\em Metric spaces of non-positive curvature}, volume 319 of {\em
		Grundlehren der Mathematischen Wissenschaften [Fundamental Principles of
		Mathematical Sciences]}.
	\newblock Springer-Verlag, Berlin, 1999.
	
	\bibitem{MR814144}
	Dan Burghelea.
	\newblock The cyclic homology of the group rings.
	\newblock {\em Comment. Math. Helv.}, 60(3):354--365, 1985.
	
	\bibitem{MR806699}
	Jeff Cheeger and Mikhael Gromov.
	\newblock Bounds on the von {N}eumann dimension of {$L^2$}-cohomology and the
	{G}auss-{B}onnet theorem for open manifolds.
	\newblock {\em J. Differential Geom.}, 21(1):1--34, 1985.
	
	\bibitem{AC88}
	A.~Connes.
	\newblock Entire cyclic cohomology of {B}anach algebras and characters of
	{$\theta$}-summable {F}redholm modules.
	\newblock {\em $K$-Theory}, 1(6):519--548, 1988.
	
	\bibitem{Connes}
	Alain Connes.
	\newblock Noncommutative differential geometry.
	\newblock {\em Inst. Hautes \'{E}tudes Sci. Publ. Math.}, (62):257--360, 1985.
	
	\bibitem{MR1116414}
	Alain Connes.
	\newblock On the {C}hern character of {$\theta$} summable {F}redholm modules.
	\newblock {\em Comm. Math. Phys.}, 139(1):171--181, 1991.
	
	\bibitem{ConnesNCG}
	Alain Connes.
	\newblock {\em Noncommutative geometry}.
	\newblock Academic Press, Inc., San Diego, CA, 1994.
	
	\bibitem{CGM90}
	Alain Connes, Mikha{\"{\i}}l Gromov, and Henri Moscovici.
	\newblock Conjecture de {N}ovikov et fibr{\'e}s presque plats.
	\newblock {\em C. R. Acad. Sci. Paris S{\'e}r. I Math.}, 310(5):273--277, 1990.
	
	\bibitem{CM90}
	Alain Connes and Henri Moscovici.
	\newblock Cyclic cohomology, the {N}ovikov conjecture and hyperbolic groups.
	\newblock {\em Topology}, 29(3):345--388, 1990.
	
	\bibitem{Harpe}
	Pierre de~la Harpe.
	\newblock Groupes hyperboliques, alg\`ebres d'op\'{e}rateurs et un
	th\'{e}or\`eme de {J}olissaint.
	\newblock {\em C. R. Acad. Sci. Paris S\'{e}r. I Math.}, 307(14):771--774,
	1988.
	
	\bibitem{MR3563244}
	Robin~J. Deeley and Magnus Goffeng.
	\newblock Realizing the analytic surgery group of {H}igson and {R}oe
	geometrically part {III}: higher invariants.
	\newblock {\em Math. Ann.}, 366(3-4):1513--1559, 2016.
	
	\bibitem{MR511246}
	Harold Donnelly.
	\newblock Eta invariants for {$G$}-spaces.
	\newblock {\em Indiana Univ. Math. J.}, 27(6):889--918, 1978.
	
	\bibitem{MR836727}
	Ezra Getzler.
	\newblock A short proof of the local {A}tiyah-{S}inger index theorem.
	\newblock {\em Topology}, 25(1):111--117, 1986.
	
	\bibitem{EG93b}
	Ezra Getzler.
	\newblock Cyclic homology and the {A}tiyah-{P}atodi-{S}inger index theorem.
	\newblock In {\em Index theory and operator algebras ({B}oulder, {CO}, 1991)},
	volume 148 of {\em Contemp. Math.}, pages 19--45. Amer. Math. Soc.,
	Providence, RI, 1993.
	
	\bibitem{MR3500822}
	Alexander Gorokhovsky, Hitoshi Moriyoshi, and Paolo Piazza.
	\newblock A note on the higher {A}tiyah-{P}atodi-{S}inger index theorem on
	{G}alois coverings.
	\newblock {\em J. Noncommut. Geom.}, 10(1):265--306, 2016.
	
	\bibitem{Gromov}
	M.~Gromov.
	\newblock Hyperbolic groups.
	\newblock In {\em Essays in group theory}, volume~8 of {\em Math. Sci. Res.
		Inst. Publ.}, pages 75--263. Springer, New York, 1987.
	
	\bibitem{Higson1}
	Nigel Higson and John Roe.
	\newblock Mapping surgery to analysis. {I}. {A}nalytic signatures.
	\newblock {\em $K$-Theory}, 33(4):277--299, 2005.
	
	\bibitem{Higson2}
	Nigel Higson and John Roe.
	\newblock Mapping surgery to analysis. {II}. {G}eometric signatures.
	\newblock {\em $K$-Theory}, 33(4):301--324, 2005.
	
	\bibitem{Higson3}
	Nigel Higson and John Roe.
	\newblock Mapping surgery to analysis. {III}. {E}xact sequences.
	\newblock {\em $K$-Theory}, 33(4):325--346, 2005.
	
	\bibitem{JLO88}
	Arthur Jaffe, Andrzej Lesniewski, and Konrad Osterwalder.
	\newblock Quantum {$K$}-theory. {I}. {T}he {C}hern character.
	\newblock {\em Comm. Math. Phys.}, 118(1):1--14, 1988.
	
	\bibitem{Sheagan}
	Sheagan John.
	\newblock Higher rho invariants, cyclic cohomology and groups of polynomial
	growth.
	\newblock thesis in prepration.
	
	\bibitem{Jolissaint}
	Paul Jolissaint.
	\newblock Rapidly decreasing functions in reduced {$C^*$}-algebras of groups.
	\newblock {\em Trans. Amer. Math. Soc.}, 317(1):167--196, 1990.
	
	\bibitem{GK88}
	Gennadi Kasparov.
	\newblock Equivariant {$KK$}-theory and the {N}ovikov conjecture.
	\newblock {\em Invent. Math.}, 91(1):147--201, 1988.
	
	\bibitem{MR1998480}
	Gennadi Kasparov and Georges Skandalis.
	\newblock Groups acting properly on ``bolic'' spaces and the {N}ovikov
	conjecture.
	\newblock {\em Ann. of Math. (2)}, 158(1):165--206, 2003.
	
	\bibitem{Keswani}
	Navin Keswani.
	\newblock Geometric {$K$}-homology and controlled paths.
	\newblock {\em New York J. Math.}, 5:53--81, 1999.
	
	\bibitem{Laff}
	Vincent Lafforgue.
	\newblock A proof of property ({RD}) for cocompact lattices of {${\rm
			SL}(3,\bold R)$} and {${\rm SL}(3,\bold C)$}.
	\newblock {\em J. Lie Theory}, 10(2):255--267, 2000.
	
	\bibitem{L1}
	Vincent Lafforgue.
	\newblock {$K$}-th\'{e}orie bivariante pour les alg\`ebres de {B}anach et
	conjecture de {B}aum-{C}onnes.
	\newblock {\em Invent. Math.}, 149(1):1--95, 2002.
	
	\bibitem{MR2874956}
	Vincent Lafforgue.
	\newblock La conjecture de {B}aum-{C}onnes \`a coefficients pour les groupes
	hyperboliques.
	\newblock {\em J. Noncommut. Geom.}, 6(1):1--197, 2012.
	
	\bibitem{Leichtnam}
	Eric Leichtnam and Paolo Piazza.
	\newblock The {$b$}-pseudodifferential calculus on {G}alois coverings and a
	higher {A}tiyah-{P}atodi-{S}inger index theorem.
	\newblock {\em M\'{e}m. Soc. Math. Fr. (N.S.)}, (68):iv+121, 1997.
	
	\bibitem{MR1708639}
	{{\'E}}ric Leichtnam and Paolo Piazza.
	\newblock Homotopy invariance of twisted higher signatures on manifolds with
	boundary.
	\newblock {\em Bull. Soc. Math. France}, 127(2):307--331, 1999.
	
	\bibitem{Lott1}
	John Lott.
	\newblock Higher eta-invariants.
	\newblock {\em $K$-Theory}, 6(3):191--233, 1992.
	
	\bibitem{Lott}
	John Lott.
	\newblock Delocalized {$L^2$}-invariants.
	\newblock {\em J. Funct. Anal.}, 169(1):1--31, 1999.
	
	\bibitem{MR1618772}
	John Lott.
	\newblock Diffeomorphisms and noncommutative analytic torsion.
	\newblock {\em Mem. Amer. Math. Soc.}, 141(673):viii+56, 1999.
	
	\bibitem{RM93}
	Richard~B. Melrose.
	\newblock {\em The {A}tiyah-{P}atodi-{S}inger index theorem}, volume~4 of {\em
		Research Notes in Mathematics}.
	\newblock A K Peters Ltd., Wellesley, MA, 1993.
	
	\bibitem{Meyer}
	Ralf Meyer.
	\newblock Combable groups have group cohomology of polynomial growth.
	\newblock {\em Q. J. Math.}, 57(2):241--261, 2006.
	
	\bibitem{Mineyev}
	I.~Mineyev.
	\newblock Straightening and bounded cohomology of hyperbolic groups.
	\newblock {\em Geom. Funct. Anal.}, 11(4):807--839, 2001.
	
	\bibitem{MR1914618}
	Igor Mineyev and Guoliang Yu.
	\newblock The {B}aum-{C}onnes conjecture for hyperbolic groups.
	\newblock {\em Invent. Math.}, 149(1):97--122, 2002.
	
	\bibitem{Nistor}
	V.~Nistor.
	\newblock Group cohomology and the cyclic cohomology of crossed products.
	\newblock {\em Invent. Math.}, 99(2):411--424, 1990.
	
	\bibitem{Piazza}
	Paolo Piazza and Thomas Schick.
	\newblock Groups with torsion, bordism and rho invariants.
	\newblock {\em Pacific J. Math.}, 232(2):355--378, 2007.
	
	\bibitem{Piazzarho}
	Paolo Piazza and Thomas Schick.
	\newblock Rho-classes, index theory and {S}tolz' positive scalar curvature
	sequence.
	\newblock {\em J. Topol.}, 7(4):965--1004, 2014.
	
	\bibitem{Puschnigg}
	Michael Puschnigg.
	\newblock New holomorphically closed subalgebras of {$C^*$}-algebras of
	hyperbolic groups.
	\newblock {\em Geom. Funct. Anal.}, 20(1):243--259, 2010.
	
	\bibitem{MR790678}
	Daniel Quillen.
	\newblock Superconnections and the {C}hern character.
	\newblock {\em Topology}, 24(1):89--95, 1985.
	
	\bibitem{Roecoarse}
	John Roe.
	\newblock Coarse cohomology and index theory on complete {R}iemannian
	manifolds.
	\newblock {\em Mem. Amer. Math. Soc.}, 104(497):x+90, 1993.
	
	\bibitem{Roe}
	John Roe.
	\newblock {\em Index theory, coarse geometry, and topology of manifolds},
	volume~90 of {\em CBMS Regional Conference Series in Mathematics}.
	\newblock Published for the Conference Board of the Mathematical Sciences,
	Washington, DC; by the American Mathematical Society, Providence, RI, 1996.
	
	\bibitem{Charlotte}
	Charlotte Wahl.
	\newblock The {A}tiyah-{P}atodi-{S}inger index theorem for {D}irac operators
	over {$C^\ast$}-algebras.
	\newblock {\em Asian J. Math.}, 17(2):265--319, 2013.
	
	\bibitem{MR1707352}
	Shmuel Weinberger.
	\newblock Higher {$\rho$}-invariants.
	\newblock In {\em Tel {A}viv {T}opology {C}onference: {R}othenberg
		{F}estschrift (1998)}, volume 231 of {\em Contemp. Math.}, pages 315--320.
	Amer. Math. Soc., Providence, RI, 1999.
	
	\bibitem{Weinberger:2016dq}
	Shmuel Weinberger, Zhizhang Xie, and Guoliang Yu.
	\newblock Additivity of higher rho invariants and nonrigidity of topological
	manifolds.
	\newblock arXiv:1608.03661, 2016.
	
	\bibitem{MR3416114}
	Shmuel Weinberger and Guoliang Yu.
	\newblock Finite part of operator {$K$}-theory for groups finitely embeddable
	into {H}ilbert space and the degree of nonrigidity of manifolds.
	\newblock {\em Geom. Topol.}, 19(5):2767--2799, 2015.
	
	\bibitem{Xiepos}
	Zhizhang Xie and Guoliang Yu.
	\newblock Positive scalar curvature, higher rho invariants and localization
	algebras.
	\newblock {\em Adv. Math.}, 262:823--866, 2014.
	
	\bibitem{MR3590536}
	Zhizhang Xie and Guoliang Yu.
	\newblock Higher rho invariants and the moduli space of positive scalar
	curvature metrics.
	\newblock {\em Adv. Math.}, 307:1046--1069, 2017.
	
	\bibitem{Xie}
	Zhizhang Xie and Guoliang Yu.
	\newblock Delocalized eta invariants, algebraicity, and {$K$}-theory of group
	{$C^*$}-algebras.
	\newblock {\em arXiv:1805.07617}, 2018.
	
	\bibitem{Yulocalization}
	Guoliang Yu.
	\newblock Localization algebras and the coarse {B}aum-{C}onnes conjecture.
	\newblock {\em $K$-Theory}, 11(4):307--318, 1997.
	
	\bibitem{MR3551834}
	Rudolf Zeidler.
	\newblock Positive scalar curvature and product formulas for secondary index
	invariants.
	\newblock {\em J. Topol.}, 9(3):687--724, 2016.
	
\end{thebibliography}
\end{document}